\documentclass[12pt]{amsart}
\usepackage[T1]{fontenc}
\usepackage[utf8]{inputenc}
\usepackage[english]{babel}
\usepackage{lmodern,amsmath,amssymb,amsthm,mathtools}
\usepackage{microtype,needspace}
\usepackage{mathrsfs,esint,enumitem,booktabs,longtable,array,aliascnt}
\usepackage{xcolor}
\usepackage[colorlinks=true,linkcolor=blue!55!black,citecolor=green!40!black,urlcolor=blue!55!black]{hyperref}
\usepackage[nameinlink,noabbrev]{cleveref}
\numberwithin{equation}{section}
\newtheorem{theorem}{Theorem}[section]
\newtheorem{introtheorem}{Theorem}

\newtheorem{restatedtheorem}{Theorem}

\crefname{introtheorem}{Theorem}{Theorems}
\Crefname{introtheorem}{Theorem}{Theorems}
\crefname{restatedtheorem}{Theorem}{Theorems}
\Crefname{restatedtheorem}{Theorem}{Theorems}

\newaliascnt{proposition}{theorem}
\newtheorem{proposition}[proposition]{Proposition}
\aliascntresetthe{proposition}
\newaliascnt{lemma}{theorem}
\newtheorem{lemma}[lemma]{Lemma}
\aliascntresetthe{lemma}
\newaliascnt{corollary}{theorem}
\newtheorem{corollary}[corollary]{Corollary}
\aliascntresetthe{corollary}
\theoremstyle{definition}
\newaliascnt{definition}{theorem}

\aliascntresetthe{definition}
\theoremstyle{remark}
\newaliascnt{remark}{theorem}

\aliascntresetthe{remark}
\newcommand{\R}{\mathbb{R}}
\newcommand{\Sph}{\mathbb{S}}
\newcommand{\eps}{\varepsilon}
\newcommand{\dd}{\,\mathrm{d}}

\newcommand{\ip}[2]{\left\langle#1,#2\right\rangle}

\allowdisplaybreaks[2]
\title[The flat singular Escobar-Yamabe equation]{Global Classification of the Flat Singular Escobar-Yamabe Equation}
\author{Felipe G. Bogo}
\address{Instituto de Matem\'atica e Estat\'istica, Universidade Federal
Fluminense, Niter\'oi, RJ, Brasil, 24210-201.}
\email{felipebogo@id.uff.br}
\hypersetup{pdfauthor={Felipe G. Bogo}}
\date{}

\begin{document}

\begin{abstract}
In every dimension $n\geq3$, we classify positive solutions of
the flat Escobar--Yamabe equation on the Euclidean half-space
with a possible isolated boundary singularity at the origin.
Solutions with a removable singularity are boundary bubbles.
Up to dilation, the nonremovable solutions are the cylindrical
solution and a family of periodic Delaunay solutions uniquely
parametrized by the necksize $\eps\in(0,1)$, the minimum of the
normalized boundary trace. In each fundamental period $T_\eps$,
the trace has exactly two critical points, a nondegenerate
minimum and a nondegenerate maximum. The period and the solution
centered at the trace minimum and rescaled to period one depend
smoothly on $\eps$. As $\eps\to0$,
$T_\eps=\frac{4}{n-2}\log(1/\eps)+O(1)$.
As $\eps\to1$, the solutions converge to the cylindrical solution
and their periods tend to that of its linearized oscillation.
We use uniqueness in the odd and even classes for linear Robin
equations to obtain reflection symmetry and continuation in
the necksize.
\end{abstract}
\maketitle
\enlargethispage*{0pt}
\vspace{-24pt}

\setcounter{tocdepth}{1}
\tableofcontents
\vspace{-24pt}

\section{Introduction}\label{sec:introduction}

Let $n\geq3$. We study positive solutions of the singular
Escobar--Yamabe problem
\begingroup
\predisplaypenalty=0\relax
\begin{equation}\label{eq:intro-flat-equation}
 \begin{cases}
 \Delta U=0 & \text{for }x_n>0,\\
 \partial_{x_n} U + \dfrac{n-2}{2}u^{\frac n{n-2}} = 0
   & \text{on }\partial\R^n_+\setminus\{0\},
 \end{cases}
\end{equation}
\endgroup
where $\R^n_+=\{x=(\bar x,x_n)\in\R^n:x_n\geq0\}$ and
$u(\bar x)=U(\bar x,0)$. We assume that $U$ is smooth on
$\R^n_+\setminus\{0\}$. The conformal metric $g=U^{4/(n-2)}\delta$,
where $\delta$ is the Euclidean metric, is scalar flat and has
boundary mean curvature $h_g=1$ away from $0$. We take $h_g$ to
be the average of the principal curvatures with respect to the
outward unit normal.

The relation with the fractional Yamabe problem is described by
the conformal extension results of Chang and
Gonz\'alez~\cite{ChangGonzalez} and Case and Chang~\cite{CaseChang}.
At order $s=1/2$, the Poisson extension of a solution of
\[
 (-\Delta_{\R^{n-1}})^{1/2}u=\frac{n-2}{2}u^{\frac n{n-2}}
 \qquad\text{in }\R^{n-1}\setminus\{0\}
\]
solves~\eqref{eq:intro-flat-equation} by the extension formula of
Caffarelli and Silvestre~\cite{CS}.

We say that the singularity at $0$ is removable if the trace $u$
extends continuously across $0$. In the nonremovable case,
Caffarelli, Jin, Sire and Xiong~\cite[Theorem~1.1]{CJXS} prove that
\[
 c_1|\bar x|^{-\frac{n-2}{2}}
 \leq u(\bar x)
 \leq c_2|\bar x|^{-\frac{n-2}{2}}
\]
for some $c_1,c_2>0$ and all nonzero $\bar x$ sufficiently close to
$0$. If the singularity is removable, the weak equation extends
across $0$ by the argument in
Section~\ref{subsec:global-kelvin-reduction}. The half-space
Liouville theorem of Jin, Li and Xiong~\cite[Theorem~1.8]{JLX}
then gives the boundary bubbles in~\eqref{eq:intro-flat-bubble}.
The classification for regular solutions is due to
Li and Zhu~\cite[Theorem~1.2]{LiZhu}.

The case without boundary was considered by Caffarelli, Gidas and
Spruck~\cite[Theorem~8.1(a) and Corollary~8.2(b)]{CGS}. Positive
solutions of
\begin{equation}\label{eq:intro-classical-equation}
 \Delta w+\frac{n(n-2)}4 w^{\frac{n+2}{n-2}}=0
 \qquad\text{in }\R^n\setminus\{0\}
\end{equation}
are standard bubbles if the singularity is removable and are radial
if it is nonremovable. For a radial solution, the Emden--Fowler
transformation is
\[
 w(r)=r^{-\frac{n-2}{2}}v(t),\qquad r=|x|,\qquad t=-\log r.
\]
This gives the ordinary differential equation
\begin{equation}\label{eq:intro-cylindrical-ode}
 v''+\frac{n(n-2)}4v^{\frac{n+2}{n-2}}
       -\frac{(n-2)^2}4v=0.
\end{equation}
The Hamiltonian energy
\begin{equation}\label{eq:intro-hamiltonian}
 H(v,v')=(v')^2-\frac{(n-2)^2}4v^2
                 +\frac{(n-2)^2}4v^{\frac{2n}{n-2}}
\end{equation}
is constant along solutions $v$ of~\eqref{eq:intro-cylindrical-ode}.
Its periodic orbits give the Fowler solutions, also called
Delaunay-type solutions after Delaunay's surfaces of
revolution~\cite{Delaunay}.

\begin{proposition}[Mazzeo and Pacard~{\cite[Proposition~1]{MP}}]
\label{prop:intro-mazzeo-pacard-prop-1}
For every energy $H$ in the interval $ (-(n-2)((n-2)/n)^{n/2}/2,\,0)$ there is a unique positive bounded solution $v$
of~\eqref{eq:intro-cylindrical-ode} with that energy and with
$v'(0)=0$ and $v''(0)>0$. This solution is periodic, takes values in
$(0,1)$ and can be parametrized by its minimum
\[
 \eps=v(0)\in\left(0,\left(\frac{n-2}{n}\right)^{\frac{n-2}{4}}\right).
\]
At the left endpoint of the energy interval the only positive
bounded solution is the constant
\[
 v\equiv\left(\frac{n-2}{n}\right)^{\frac{n-2}{4}}.
\]
At energy zero, a nonnegative bounded solution is either zero or
\[
 v(t)=(\cosh(t-t_0))^{-\frac{n-2}{2}}
 \qquad\text{for some }t_0\in\R.
\]
\end{proposition}

For the singular Escobar--Yamabe problem, Caffarelli, Jin, Sire
and Xiong~\cite[Theorem~1.3]{CJXS} show that every solution with a
nonremovable singularity is cylindrically symmetric, that is,
\[
 U(\bar x,x_n)=U(|\bar x|,x_n).
\]
The Emden--Fowler transform can depend on both $t$ and the
angle with the $x_n$-axis, so the reduced equation is a partial
differential equation. DelaTorre and
Gonz\'alez~\cite{DG} studied the fractional equation, its Hamiltonian
and the linearization at the cylindrical solution. DelaTorre,
del Pino, Gonz\'alez and Wei~\cite[Theorem~1.1]{DDPGW} constructed
periodic solutions by a variational argument. At order $1/2$,
their theorem applies when the half-space has dimension $n\geq4$.
Cruz-Bl\'azquez, DelaTorre and Ruiz~\cite[Section~1]{CDR} asked
whether every solution of the fractional equation whose
Emden--Fowler transform is bounded above and away from zero
is periodic in $t$. Andrade, DelaTorre, do~\'O,
Ratzkin and Wei~\cite[Theorem~1]{Andrade} classified the singular
solutions for fractional orders sufficiently close to $1$.

The normalized Emden--Fowler transform $v$ in
Section~\ref{subsec:notation} solves~\eqref{eq:boundary-cylinder}.
Its boundary coefficient $c$ and the normalization factor $b_n$ are
\[
 c=2\left(\frac{\Gamma(n/4)}{\Gamma((n-2)/4)}\right)^2,
 \qquad b_n=\left(\frac{2c}{n-2}\right)^{\frac{n-2}{2}}.
\]
For cylindrically symmetric solutions we write $\tilde v(t)$
for the trace on $\R\times\partial\Sph^{n-1}_+$. The cylindrical
solution $C_n$ has trace $1$ in this normalization. We call the
minimum of $\tilde v$ the \emph{necksize}.

\begin{introtheorem}\label{intro:necksize-family}
Let $n\geq3$. For every $\eps\in(0,1)$, there is a unique positive
smooth cylindrically symmetric solution $v_\eps$ of
\eqref{eq:boundary-cylinder} with bounded trace and
\[
 \tilde v_\eps(0)=\min_{t\in\R}\tilde v_\eps(t)=\eps.
\]
The solution is even in $t$ and periodic. Its fundamental period
$T_\eps$ is twice the first positive critical point of $\tilde v_\eps$.
The trace $\tilde v_\eps$ has one minimum and one maximum per period,
with $\tilde v_\eps''(0)>0$, $\tilde v_\eps''(T_\eps/2)<0$ and
$\tilde v_\eps'(t)>0$ for $0<t<T_\eps/2$.

Every positive smooth cylindrically symmetric solution $v$
of~\eqref{eq:boundary-cylinder} with trace bounded above and away
from zero is either $C_n$ or has the form
$v(t,\theta)=v_\eps(t+\tau,\theta)$ for a unique $\eps\in(0,1)$
and a phase $\tau\in\R$ determined modulo $T_\eps$. We set
$v_1=C_n$.

With $v_\eps^{(1)}(s,\theta)=v_\eps(T_\eps s,\theta)$, the map
$\eps\mapsto(T_\eps,v_\eps^{(1)})$ is smooth from $(0,1)$ into
\[
 (0,\infty)\times
 C^{k,\gamma}_{1,\mathrm{sym,even}}(\R\times\Sph^{n-1}_+)
\]
for every integer $k\geq0$ and $0<\gamma<1$.
\end{introtheorem}

\medskip
\noindent\emph{Idea of the proof of Theorem~\ref{intro:necksize-family}.}\enspace
Proposition~\ref{s3:rigidity} gives reflection symmetry about
each critical point of $\tilde v$ and determines an even
solution by $\tilde v(0)$. The differences of solutions and
the Jacobi fields used in that proposition satisfy linear
equations with a coefficient $b(t)$ given by \eqref{eq:s3-secant}.
For $n\geq5$, we estimate the trace by decomposing the solution
along the first eigenfunction and its orthogonal complement.
For $n=3,4$, we obtain a contraction for the trace using
the eigenfunctions for a fixed value of the boundary coefficient.

Testing~\eqref{eq:boundary-cylinder} against $C_n$ gives
\eqref{eq:s3-moment-equation}, which yields two distinct critical
points of $\tilde v$ for every noncylindrical solution under
consideration. Composing the reflections about these points
gives a period. Odd uniqueness for the Jacobi field $\partial_t v$
gives nondegeneracy of the critical points. Even uniqueness gives
uniqueness for a prescribed necksize after translating a minimum to
$t=0$.

Minimization on a compact cylinder of sufficiently large period
gives one noncylindrical solution. After rescaling to period one,
we treat the period as an unknown and prescribe $\tilde v(0)$.
Even uniqueness and the Fredholm alternative give invertibility
of the linearization of this augmented problem. The implicit
function theorem and compactness then give, respectively, openness
and relative closedness of the set of attained necksizes in $(0,1)$.

As $\eps\to1$, $v_\eps$ converges to $C_n$ and $T_\eps$ tends to
the period determined by the linearization there
(Proposition~\ref{s3:endpoint-periods}). We prove that $\eps\to0$ implies,
\[
 T_\eps=\frac{4}{n-2}\log(1/\eps)+O(1)
\]
in Corollary~\ref{s3:logarithmic-period}.

\begin{introtheorem}\label{intro:flat-classification}
Let $n\geq3$ and let $U$ be a positive solution of
\eqref{eq:intro-flat-equation}, smooth on
$\R^n_+\setminus\{0\}$. Exactly one of the following holds.
\begin{enumerate}[label=\textup{(\roman*)}]
\item The singularity at $0$ is removable, and $U$ is the bubble
\begin{equation}\label{eq:intro-flat-bubble}
 U(\bar x,x_n)=
 \left(\frac{2\lambda}{|\bar x-a|^2+(x_n+\lambda)^2}\right)^{\frac{n-2}{2}}
\end{equation}
for unique parameters $a\in\R^{n-1}$ and $\lambda>0$.

\item The singularity at $0$ is nonremovable. There exist a unique
necksize $\eps\in(0,1]$ and a phase $\tau\in\R$ such that,
for $r>0$ and $\theta\in\Sph^{n-1}_+$,
\[
 U(r\theta)= b_nr^{-(n-2)/2}\,
                 v_\eps\bigl(-\log r+\tau,\theta\bigr),
\]
where $v_\eps$ is given by Theorem~\ref{intro:necksize-family}.
For $0<\eps<1$, the phase is unique modulo $T_\eps$. For
$\eps=1$, all phases give the same cylindrical solution.
\end{enumerate}
Conversely, both formulas give positive solutions of
\eqref{eq:intro-flat-equation}. The solutions in \textup{(i)}
extend smoothly across $0$, and those in \textup{(ii)} have a
nonremovable singularity there.
\end{introtheorem}

\medskip
\noindent\emph{Idea of the proof of Theorem~\ref{intro:flat-classification}.}\enspace
If $u$ extends continuously across $0$, the weak equation for
$U$ extends across $0$ and the regular half-space classification
gives~\eqref{eq:intro-flat-bubble}. In the nonremovable case, we
apply the local two-sided bounds at $0$ to $U$ and its Kelvin
transform $U_{0,1}$ from~\eqref{eq:flat-kelvin}. The trace of
$U_{0,1}$ is also singular at $0$, since Kelvin inversion sends
bubbles to bubbles. This gives two-sided bounds for
$|x|^{(n-2)/2}U(x)$ on $\R^n_+\setminus\{0\}$. Cylindrical symmetry
and Theorem~\ref{intro:necksize-family} then give the formula in
\textup{(ii)}.

In the companion paper~\cite{PartII}, we use these results to study
the Jacobi fields of the Delaunay family and the asymptotic behavior
of solutions near an isolated boundary singularity. For a smooth
background metric, two-sided bounds on the boundary restriction
give convergence to a fixed Delaunay solution. The corresponding
problem without boundary was studied by Korevaar, Mazzeo, Pacard
and Schoen~\cite{KMPS} and by Marques~\cite{Marques}.

The classification also supplies models for constructing singular
solutions of the boundary Yamabe problem on manifolds with boundary.
In the classical Yamabe problem, Silva Santos~\cite{SilvaSantos}
used gluing to construct complete conformal metrics with Delaunay-type
ends on punctured closed manifolds under nondegeneracy and curvature
hypotheses. For the fractional Yamabe equation, Ao, DelaTorre,
Gonz\'alez and Wei~\cite[Theorem~1.1]{AoDGW} used Delaunay solutions
in a gluing construction with prescribed isolated singularities.

\subsection*{Structure of the paper}

In Section~\ref{sec:preliminaries}, Section~\ref{subsec:notation}
introduces the normalized Emden--Fowler transform and Kelvin
inversion. Section~\ref{subsec:estimates-compactness} gives uniform
regularity and compactness on the cylinder, and
Section~\ref{subsec:linearized-spaces} specifies the linearized
equations and periodic function spaces.

In Section~\ref{sec:delaunay-family},
Section~\ref{subsec:delaunay-rigidity} proves rigidity and
periodicity. Section~\ref{subsec:delaunay-local-chart} constructs
the family and proves its smooth dependence on the necksize.
Section~\ref{subsec:delaunay-hamiltonian} records the Hamiltonian
and determines the limiting behavior of the period.

In Section~\ref{sec:global-classification},
Section~\ref{subsec:global-kelvin-reduction} treats removable
singularities and obtains the global two-sided bounds in the
nonremovable case. Section~\ref{subsec:global-classification-proof}
uses these results and Theorem~\ref{intro:necksize-family} to
complete the proof of Theorem~\ref{intro:flat-classification}.

In Appendix~\ref{sec:flat-foundations},
Section~\ref{subsec:flat-extension-bounds} controls the extension
on $\R\times\Sph^{n-1}_+$ by bounds on its trace.
Section~\ref{subsec:flat-trace-bounds} bounds the first derivative
of the trace and proves second-order inequalities for powers of the
trace, which are used to estimate the boundary coefficients. Section~\ref{subsec:flat-spectral-estimates}
gives the eigenvalue bounds and the resolvent identity.
In Appendix~\ref{flatref:section},
Section~\ref{app:weighted-identity} estimates a weighted integral
of the trace, and Section~\ref{app:even-jacobi} uses
this estimate to prove uniqueness in each parity for $n\geq5$.
Appendix~\ref{nfour:section} obtains the corresponding uniqueness
statements for $n=3,4$ from integral equations for the trace and
estimates for their kernels.

\subsection*{Relation to concurrent work}

While this manuscript was in preparation, Jin, Liao and Xiong
submitted their preprint~\cite{JinLiaoXiong} to arXiv on
November~5, 2026. Their Theorems~1.1 and~2.1 give the same
half-space classification. They parametrize the linear solution
space by two projections of the Cauchy data in Proposition~2.10
and classify the nonlinear solutions by their Hamiltonian energy.
Their Proposition~A.2 also gives the necksize parametrization by
showing that the necksize is a real-analytic function of the
energy with a real-analytic inverse. Their Lemma~A.1 proves
uniqueness for even solutions from $Y(0,L)$ under the additional
bound $|b'|\leq b$.

The scalar resolvent identity in the first line
of~\eqref{flatbase:robin-resolvent-identity} also appears
in~\cite[(2.20)]{JinLiaoXiong}. The weighted
identity~\eqref{flatref:exact-dissipation} is analogous
to~\cite[(A.4)]{JinLiaoXiong}. To obtain the weighted
estimate, both papers decompose the solution along the first
Robin eigenfunction and its orthogonal complement and integrate
the corresponding identity against $e^{-(n-1)t}$.
Appendices~\ref{sec:flat-foundations} and~\ref{flatref:section}
give the derivations used in this paper.

We use uniqueness from $Y(0,L)$ in the even class and from
$\partial_tY(0,L)$ in the odd class. We deduce periodicity
by composing reflections about trace critical points and construct
the family by continuation directly in the necksize.

Corollary~\ref{s3:logarithmic-period} gives the rate of divergence
of the period at small necksize,
$T_\eps=\frac{4}{n-2}\log(1/\eps)+O(1)$.
Proposition~\ref{flatbase:pressure} gives second-order inequalities
for powers of the trace in dimensions $n\geq5$. For $n=3,4$,
Theorem~\ref{nfour:linear-uniqueness} establishes uniqueness in
each parity from the corresponding boundary datum, assuming only
that $b$ is smooth and even and $0\leq b\leq n/2$. These estimates
and the different proof of the classification are the reasons for
presenting our results here.

\subsection*{Statement on the use of AI}

The need for Proposition~\ref{s3:rigidity} was identified before
we used AI. The model GPT-Sol 6 was instrumental in finding the
method of proof of that proposition. We also used the model to
revise the manuscript. All mathematical statements, proofs, and citations were independently checked
by the author, who takes full responsibility for the contents of the article. The final manuscript was
written and revised by the author.

\subsection*{Acknowledgments}

This work forms part of my Ph.D.\ thesis at Universidade Federal
Fluminense. I am grateful to my advisor, S\'ergio Almaraz, for his
guidance and advice throughout this work. I also thank
Xuezhang Chen for his warm welcome at Nanjing University and his
generous support during my stay.

This study was financed in part by the Coordena\c{c}\~ao de
Aperfei\c{c}oamento de Pessoal de N\'{\i}vel Superior - Brasil
(CAPES) - Finance Code 001. I received support under process
No.~88887.815246/2023-00 and a scholarship through the Programa
Institucional de Doutorado Sandu\'{\i}che no Exterior (PDSE),
process No.~88881.127484/2025-01.

\newpage

\section{Preliminaries}\label{sec:preliminaries}

Section~\ref{subsec:notation} gives the cylindrical formulation
of~\eqref{eq:intro-flat-equation} and describes the action of
dilations and Kelvin inversion. Section~\ref{subsec:estimates-compactness}
establishes uniform regularity and compactness for solutions
bounded above and away from zero. In Section~\ref{subsec:linearized-spaces},
we introduce the linearized equations and the estimates for
Jacobi fields with polynomial growth used in the rigidity argument.
We also define the periodic H\"older spaces for the local
parametrization by necksize.

\subsection{Notation}\label{subsec:notation}
\label{subsec:flat-cylinder}\label{subsec:flat-stationary}

We use the half-space and trace conventions of
Section~\ref{sec:introduction}. We introduce the normalized
Emden--Fowler variables, write the interior equation and the
boundary condition on the cylinder, and express dilations as
translations in $t$. Cylindrical symmetry means invariance under
rotations of $\bar x$ with $x_n$ fixed. For cylindrically symmetric
functions, we give the angular form of the Laplacian, the conormal
derivative at the equator, and the regularity condition at the
pole. We then record the Kelvin transform and its action in $t$
when the center of inversion is the origin.

Write $x=r\theta$ with $r>0$ and
$\theta\in\Sph^{n-1}_+=\{\theta\in\Sph^{n-1}:\theta_n\geq0\}$.
In the normalization of the introduction, the Emden--Fowler
variables are
\[
 t=-\log r,\qquad
 v(t,\theta)=b_n^{-1}r^{\frac{n-2}{2}}U(r\theta).
\]
We write $\tilde v$ for the restriction of $v$ to
$\R\times\partial\Sph^{n-1}_+$. Thus
$\tilde v(t,\vartheta)=v(t,(\vartheta,0))$ for
$\vartheta\in\Sph^{n-2}$. The dilation
$U_\lambda(x)=\lambda^{(n-2)/2}U(\lambda x)$ corresponds to
$v_\lambda(t,\theta)=v(t-\log\lambda,\theta)$.

In these variables, \eqref{eq:intro-flat-equation} becomes
\begin{equation}\label{eq:boundary-cylinder}
 \begin{cases}
 v_{tt}+\Delta_{\Sph^{n-1}_+}v-\dfrac{(n-2)^2}{4}v=0
       &\text{in }\R\times\Sph^{n-1}_+,\\
 \partial_{\nu_{\Sph^{n-1}_+}}v=c\,\tilde v^{\frac n{n-2}}
       &\text{on }\R\times\partial\Sph^{n-1}_+.
 \end{cases}
\end{equation}
Here $\Delta_{\Sph^{n-1}_+}$ is the Laplace--Beltrami operator
of the round metric, and $\nu_{\Sph^{n-1}_+}$ is the outward unit
conormal along $\partial\Sph^{n-1}_+$. The cylinder carries the
product metric with volume element $\dd\mu=\dd t\,\dd\theta$.

For $\theta\in\Sph^{n-1}_+$, let $\beta=\arccos\theta_n$ and
$L=\pi/2$. Then
\[
 \theta=(\sin\beta\,\vartheta,\cos\beta),\qquad
 0\leq\beta\leq L,\qquad \vartheta\in\Sph^{n-2}.
\]
For cylindrically symmetric functions, we write $v=v(t,\beta)$ and
$\tilde v(t)=v(t,L)$. For $0<\beta<L$, the Laplace--Beltrami
operator on such functions is
\[
 \Delta_{\Sph^{n-1}_+}v
   =\omega^{-1}(\omega v_\beta)_\beta
   =v_{\beta\beta}+(n-2)\cot\beta\,v_\beta,
 \qquad \omega(\beta)=\sin^{n-2}\beta.
\]
On $\R\times\partial\Sph^{n-1}_+$ we have
$\partial_{\nu_{\Sph^{n-1}_+}}v=v_\beta(t,L)$.
For each fixed $t$ and $\vartheta$, the restriction of $v$
to $s\mapsto(\sin s\,\vartheta,\cos s)$ is smooth and even near
$s=0$ by cylindrical symmetry. Hence $v_\beta(t,0)=0$.

\begin{samepage}
For $a=(\bar a,0)$ and $\rho>0$, the Kelvin transform is
\begin{equation}\label{eq:flat-kelvin}
 U_{a,\rho}(x)=
 \left(\frac{\rho}{|x-a|}\right)^{n-2}
 U\left(a+\frac{\rho^2(x-a)}{|x-a|^2}\right).
\end{equation}
It preserves~\eqref{eq:intro-flat-equation} on the inverted
domain. For $a=0$, the map $r\mapsto \rho^2/r$ gives the
cylindrical action $t\mapsto-t-2\log\rho$.
\end{samepage}

\subsection{Uniform regularity and compactness of solutions}
\label{subsec:estimates-compactness}

We define uniformly local H\"older norms on the cylinder and
its boundary. For solutions with common positive lower and upper
bounds, we first obtain a H\"older estimate in half-space
coordinates and then use Schauder estimates on the cylinder
to control higher derivatives. The estimates are uniform under
translations in $t$ and give convergence of a subsequence on
every finite cylinder. They also apply when the boundary
coefficient is $1$.

For an integer $k\geq0$ and $0<\gamma<1$, the space
$C^{k,\gamma}_{\mathrm{uloc}}(\R\times\Sph^{n-1}_+)$ consists of
locally $C^{k,\gamma}$ functions for which
\[
 \|f\|_{C^{k,\gamma}_{\mathrm{uloc}}(\R\times\Sph^{n-1}_+)}
 =\sup_{\tau\in\R}
    \|f\|_{C^{k,\gamma}([\tau-1,\tau+1]\times\Sph^{n-1}_+)}
 <\infty.
\]
The same definition applies on $\R\times\partial\Sph^{n-1}_+$.

\begin{proposition}\label{prop:flat-uniform-estimates}
Let $v$ be a positive smooth solution of~\eqref{eq:boundary-cylinder}
with $m\leq v\leq M$ on $\R\times\Sph^{n-1}_+$ for some
$m,M>0$. For every integer $k\geq0$ and
$0<\gamma<1$, there is a constant $C=C(n,m,M,k,\gamma)$ such that
\begin{equation}\label{eq:cylindrical-uniform-regularity}
 \|v\|_{C^{k,\gamma}_{\mathrm{uloc}}(\R\times\Sph^{n-1}_+)}
 \leq C.
\end{equation}
Every sequence of such solutions with the same bounds has a
subsequence converging in $C^{k,\gamma}$ on each finite cylinder,
for every $k\geq0$ and $0<\gamma<1$, to a solution with the same
bounds. The same conclusions hold when the boundary coefficient
is $1$.
\end{proposition}

\begin{proof}
Let $U$ be the solution of~\eqref{eq:intro-flat-equation}
corresponding to $v$.
For $\tau\in\R$, the dilation in Section~\ref{subsec:notation} gives
\[
 U_{e^{-\tau}}(r\theta)
 =b_n r^{-\frac{n-2}{2}}v(\tau-\log r,\theta).
\]
Consider the half-annulus
\[
 A^+=\{x\in\R^n_+:e^{-3}<|x|<e^3\}.
\]
The bound $v\leq M$ controls $U_{e^{-\tau}}$ on $A^+$,
including its trace on $A^+\cap\partial\R^n_+$, uniformly in $\tau$.
Choose finitely many cylinders $B_\rho^{n-1}(\bar y)\times[0,\rho)$
contained in $A^+$ whose cylinders of half the radius cover
$\{e^{-1}\leq |x|\leq e\}$ near $\partial\R^n_+$.
By~\eqref{eq:intro-flat-equation}, $-\partial_{x_n}U_{e^{-\tau}}$
is uniformly bounded in $L^n$ on each base. Since $n$ exceeds
the boundary dimension $n-1$, the H\"older estimate of Jin, Li
and Xiong~\cite[Proposition~2.6(iii)]{JLX} at order $1/2$ applies.
Combining it with the interior estimates for harmonic functions and
returning to cylindrical variables gives
\[
 \|v\|_{C^{0,\alpha}([\tau-1,\tau+1]\times\Sph^{n-1}_+)}
 \leq C(n,M)
\]
for some $0<\alpha<1$ independent of $\tau$.

For $j\geq0$, applying the local Schauder estimates for the
Neumann problem~\cite[Theorem~9.3]{ADN} to~\eqref{eq:boundary-cylinder}
with conormal derivative $c\tilde v^{n/(n-2)}$ on
$\R\times\partial\Sph^{n-1}_+$ gives
\[
 \|v\|_{C^{j+1,\alpha}([\tau-1,\tau+1]\times\Sph^{n-1}_+)}
 \leq C\left(
 M+\|c\tilde v^{\frac n{n-2}}\|_
 {C^{j,\alpha}([\tau-2,\tau+2]\times\partial\Sph^{n-1}_+)}
 \right).
\]
The constant $C$ is independent of $\tau$ by translation
invariance. To bound the H\"older quotients on $[\tau-2,\tau+2]$,
we use the interval of radius $1$ centered at the midpoint of the
two time coordinates when they differ by at most $2$. For the
remaining pairs, the supremum norms control the H\"older
quotients. Since $s\mapsto c s^{n/(n-2)}$ is smooth on $[m,M]$,
the chain rule bounds the right-hand side uniformly in $\tau$
whenever $v$ has a uniform $C^{j,\alpha}_{\mathrm{uloc}}$ bound.
Induction gives these bounds for all $j\geq0$. The bound in
$C^{k+1,\alpha}_{\mathrm{uloc}}$ implies
\eqref{eq:cylindrical-uniform-regularity} for every $0<\gamma<1$.
When the boundary coefficient is $1$, we apply the estimates
to $c^{-(n-2)/2}v$.

For a sequence with the same bounds, diagonal extraction on
$[-j,j]\times\Sph^{n-1}_+$ gives convergence in $C^j$ there for
every integer $j\geq1$. Taking $j\geq k+1$ large enough to cover
a given finite cylinder gives convergence in $C^{k,\gamma}$
on that cylinder. The limit satisfies~\eqref{eq:boundary-cylinder}
and the same pointwise bounds.
\end{proof}

\subsection{The linearized equations and periodic function spaces}
\label{subsec:linearized-spaces}

We linearize the interior equation and the boundary condition
and define Jacobi fields with polynomial growth. Under uniform
bounds on the linearized boundary coefficient and its derivatives,
local estimates give H\"older bounds with the same order of
growth as the field. We then define periodic H\"older spaces on
the cylinder and its boundary and specify the subspaces with
cylindrical symmetry and prescribed parity in $t$. Along an
even, cylindrically symmetric periodic solution, the linearized
operators preserve the period and the symmetry conditions.

Let $v$ be a solution of~\eqref{eq:boundary-cylinder}. For a smooth
function $z$ on $\R\times\Sph^{n-1}_+$, write $\tilde z$ for
its restriction to $\R\times\partial\Sph^{n-1}_+$. The linearized interior operator $\mathcal A$ and
the boundary operator $\mathcal B$ are
\[
 \begin{aligned}
 \mathcal A z
 &= -\partial_t^2z-\Delta_{\Sph^{n-1}_+}z
     +\dfrac{(n-2)^2}{4}z,\\
 \mathcal B z
 &= \partial_{\nu_{\Sph^{n-1}_+}}z
    -\dfrac n{n-2}c\,\tilde v^{\frac2{n-2}}\tilde z .
 \end{aligned}
\]
When the choice of $v$ needs to be explicit, we write
$\mathcal B_v$ for $\mathcal B$.
A smooth function $\Phi$ is a Jacobi field along $v$ if
\begin{equation}\label{eq:jacobi-local-pair}
 \begin{cases}
 \mathcal A\Phi=0 &\text{in }\R\times\Sph^{n-1}_+,\\
 \mathcal B\Phi=0 &\text{on }\R\times\partial\Sph^{n-1}_+.
 \end{cases}
\end{equation}

A Jacobi field $\Phi$ has polynomial growth of order at most
$q\geq0$ if there is a constant $C>0$ such that
\begin{equation}\label{eq:jacobi-polynomial-space}
 |\Phi(t,\theta)|\leq C(1+t^2)^{q/2}
 \qquad\text{for all }(t,\theta)\in\R\times\Sph^{n-1}_+.
\end{equation}
Assume that the coefficient $\frac n{n-2}c\,\tilde v^{2/(n-2)}$
in $\mathcal B$ and all its derivatives are uniformly bounded
on $\R\times\partial\Sph^{n-1}_+$. These bounds follow from
Proposition~\ref{prop:flat-uniform-estimates} when $m\leq v\leq M$
for positive constants $m$ and $M$.
For every integer $k\geq0$ and $0<\gamma<1$, the local
Schauder estimates for~\eqref{eq:jacobi-local-pair} give
\[
 \|\Phi\|_{C^{k,\gamma}([\tau-1,\tau+1]\times\Sph^{n-1}_+)}
 \leq C_{k,\gamma}
 \|\Phi\|_{C^0([\tau-2,\tau+2]\times\Sph^{n-1}_+)},
\]
where $C_{k,\gamma}$ is independent of $\tau$ and depends
only on $n,k,\gamma$ and the bounds for the coefficient of $\tilde z$
in $\mathcal B$ and its derivatives. Since $1+t^2\leq9(1+\tau^2)$
for $|t-\tau|\leq2$, \eqref{eq:jacobi-polynomial-space} implies
\[
 \|\Phi\|_{C^{k,\gamma}([\tau-1,\tau+1]\times\Sph^{n-1}_+)}
 \leq C'(1+\tau^2)^{q/2},
\]
where $C'$ is independent of $\tau$.

For $T>0$, we denote by
$C^{k,\gamma}_T(\R\times\Sph^{n-1}_+)$ the subspace of
$C^{k,\gamma}_{\mathrm{uloc}}(\R\times\Sph^{n-1}_+)$ consisting
of $T$-periodic functions,
\[
 f(t+T,\theta)=f(t,\theta)
       \qquad\text{for every }(t,\theta)\in\R\times\Sph^{n-1}_+,
\]
equipped with the norm
$\|f\|_{C^{k,\gamma}([0,T]\times\Sph^{n-1}_+)}$.
The same convention applies on $\R\times\partial\Sph^{n-1}_+$.

The subscripts $\mathrm{sym}$, $\mathrm{even}$ and $\mathrm{odd}$
indicate cylindrical symmetry, evenness in $t$ and oddness in $t$,
respectively. On $\R\times\partial\Sph^{n-1}_+$, the subscript
$\mathrm{sym}$ means independence of $\vartheta\in\Sph^{n-2}$. We combine these
subscripts with $T$ when imposing the corresponding conditions on
$T$-periodic functions. Each condition defines a closed subspace.
When $v$ is $T$-periodic, cylindrically symmetric and even in $t$,
the operators $\mathcal A$ and $\mathcal B$ preserve the period
and the indicated symmetries.

\section{The Delaunay family}
\label{sec:delaunay-family}

We prove Theorem~\ref{intro:necksize-family} for cylindrically
symmetric solutions with positive two-sided trace bounds.
Reflection symmetry, periodicity and uniqueness for solutions
with a prescribed trace minimum at $t=0$ are established in
Section~\ref{subsec:delaunay-rigidity}. The local parametrization
and compactness give the family for every necksize in $(0,1)$.
We conclude with the Hamiltonian and the behavior of the
fundamental period at both ends of the family.

\subsection{Rigidity and periodicity}
\label{subsec:delaunay-normalization}
\label{subsec:delaunay-rigidity}
\label{subsec:delaunay-periodicity}

The two-sided trace bounds control the solution on the whole
cylinder. They allow us to apply linear uniqueness to differences
of solutions and to Jacobi fields.

Let $C_n$ be the solution of the Dirichlet problem for
$-\Delta_{\Sph^{n-1}_+}+(n-2)^2/4$ with boundary value $1$.
The Dirichlet form is coercive, and the maximum principle gives
$C_n>0$. Uniqueness gives cylindrical symmetry. The constant-mode
Dirichlet-to-Neumann symbol at zero frequency
in~\cite[Theorem~1.1]{DG}, with boundary dimension $n-1$ and order
$1/2$, equals $c$ (see~\eqref{eq:s3-symbol}). Thus, with $\beta$,
$L$ and $\omega$ as in Section~\ref{subsec:notation}, $C_n$ satisfies
\begin{equation}\label{eq:s3-stationary-angular}
 \begin{cases}
 (\omega C_n')'=\dfrac{(n-2)^2}{4}\omega C_n,\\[3pt]
 C_n(L)=1,\qquad C_n'(L)=c.
 \end{cases}
\end{equation}
Smoothness and cylindrical symmetry give $C_n'(0)=0$, as in
Section~\ref{subsec:notation}.

Let $\mathscr B_n$ be the class of positive
cylindrically symmetric solutions of~\eqref{eq:boundary-cylinder}
that are smooth up to $\R\times\partial\Sph^{n-1}_+$ and satisfy,
for some $m,M>0$,
\begin{equation}\label{eq:s3-bilateral-class}
 m\leq\tilde v(t)\leq M
 \qquad\text{for every }t\in\R.
\end{equation}
For $v\in\mathscr B_n$, set
\begin{equation}\label{eq:flat-normalization-table}
 \widehat v=c^{\frac{n-2}{2}}v,\qquad
 \widehat v(t,L)=c^{\frac{n-2}{2}}\tilde v(t).
\end{equation}
Then~\eqref{eq:boundary-cylinder} gives
\begin{equation}\label{eq:s3-rp-auxiliary-pair}
 \begin{cases}
 \mathcal A\widehat v=0
     & \text{in }\R\times\Sph^{n-1}_+,\\
 \partial_\beta\widehat v(t,L)
     =\widehat v(t,L)^{\frac n{n-2}}
     & \text{for }t\in\R.
 \end{cases}
\end{equation}
Applying Proposition~\ref{flatbase:selection} to $\widehat v$
and dividing by the factor in~\eqref{eq:flat-normalization-table}
gives
\begin{equation}\label{eq:s3-canonical-comparison}
 mC_n(\beta)\leq v(t,\beta)\leq MC_n(\beta)
 \qquad\text{for every }(t,\beta)\in\R\times[0,L].
\end{equation}

Let $v_1$ and $v_2$ be bounded solutions of $\mathcal A v_i=0$ for $i=1,2$
in $\R\times\Sph^{n-1}_+$ whose traces $\tilde v_1$ and
$\tilde v_2$ depend only on $t$. If $\tilde v_1\leq\tilde v_2$
on $\R$, the maximum principle argument in the proof of
Proposition~\ref{flatbase:selection}, applied to $v_1-v_2$,
gives $v_1\leq v_2$ in $\R\times\Sph^{n-1}_+$.
Set $\eta=\|\tilde v_1-\tilde v_2\|_{L^\infty(\R)}$.
The traces of $v_2\pm\eta C_n$ are $\tilde v_2\pm\eta$, so
$v_2-\eta C_n\leq v_1\leq v_2+\eta C_n$. Equivalently,
\begin{equation}\label{eq:s3-norm-linear-comparison}
 |v_1(t,\theta)-v_2(t,\theta)|
 \leq\|\tilde v_1-\tilde v_2\|_{L^\infty(\R)}C_n(\theta).
\end{equation}

If $h$ is a period of $\tilde v$, then $v(t+h,\theta)$ and
$v(t,\theta)$ are bounded solutions with the same trace, so
they coincide by~\eqref{eq:s3-norm-linear-comparison}. The converse
follows by restriction to $\R\times\partial\Sph^{n-1}_+$. Thus
\begin{equation}\label{eq:s3-period-group-definition}
 h\text{ is a period of }v
 \quad\text{if and only if}\quad
 h\text{ is a period of }\tilde v.
\end{equation}

If $\tilde v$ is constant, say $\tilde v\equiv a>0$, then
\eqref{eq:s3-norm-linear-comparison} gives $v=aC_n$.
The boundary condition in~\eqref{eq:boundary-cylinder} and
the value of $C_n'(L)$ in~\eqref{eq:s3-stationary-angular} give
\[
 ac=aC_n'(L)=\partial_\beta v(t,L)=c\,a^{\frac n{n-2}}.
\]
Dividing by $ac$ gives $a^{2/(n-2)}=1$, so $a=1$.
In particular, $C_n$ is the unique cylindrical solution in
$\mathscr B_n$.

Proposition~\ref{flatbase:trace-gradient}, applied to
$\widehat v=c^{\frac{n-2}{2}}v$, bounds the trace uniformly over
$\mathscr B_n$. Indeed, the nonnegative derivative term
in~\eqref{flatbase:gradient} gives
$c^2\tilde v(t)^{\frac4{n-2}}\leq(n-2)^2/4$. Thus
\begin{equation}\label{eq:s3-universal-upper}
 \tilde v(t)\leq M_n
 :=\left(\frac{n-2}{2c}\right)^{\frac{n-2}{2}}
 \qquad\text{for every }t\in\R.
\end{equation}

Fix $\delta>0$. For every $v\in\mathscr B_n$ with
$\inf_\R\tilde v\geq\delta$, the bounds in
\eqref{eq:s3-canonical-comparison} hold with $m=\delta$ and
$M=M_n$. Since $C_n$ is positive on $\Sph^{n-1}_+$,
these bounds control $v$ above and away from zero in terms
of $n$ and $\delta$. By Proposition~\ref{prop:flat-uniform-estimates},
for every integer $k\geq0$ and $0<\gamma<1$ there is a constant
$C=C(n,\delta,k,\gamma)$ such that
\begin{equation}\label{eq:s3-norm-uniform-estimates}
 \sup_{\tau\in\R}
 \|v\|_{C^{k,\gamma}([\tau-1,\tau+1]\times\Sph^{n-1}_+)}
 \leq C.
\end{equation}

Let $v_1,v_2\in\mathscr B_n$ and set $w=\widehat v_1-\widehat v_2$.
By the fundamental theorem of calculus, the difference of the
boundary nonlinearities is $b(t)w(t,L)$, where
\begin{equation}\label{eq:s3-secant}
 b(t)=\frac n{n-2}\int_0^1
 \bigl((1-s)\widehat v_1(t,L)+s\widehat v_2(t,L)\bigr)^{
       \frac2{n-2}}\,\dd s.
\end{equation}
This formula also defines $b(t)$ where the traces agree.
Subtracting the equations for $\widehat v_1$ and $\widehat v_2$ gives
\begin{equation}\label{eq:s3-rp-difference-equation}
 \begin{cases}
 \mathcal A w=0 & \text{in }\R\times\Sph^{n-1}_+,\\
 \partial_\beta w(t,L)=b(t)w(t,L)
     & \text{on }\R\times\partial\Sph^{n-1}_+.
 \end{cases}
\end{equation}

When $v_1=v_2=v$, the coefficient in~\eqref{eq:s3-secant}
is $\frac n{n-2}\widehat v(t,L)^{2/(n-2)}
=\frac n{n-2}c\tilde v(t)^{2/(n-2)}$. Thus a cylindrically
symmetric Jacobi field $\Phi$ along $v$ satisfies
\begin{equation}\label{eq:s3-rp-jacobi-growth}
 \begin{cases}
 \mathcal A\Phi=0 &\text{in }\R\times\Sph^{n-1}_+,\\
 \mathcal B\Phi
   =\partial_\beta\Phi(t,L)
      -\dfrac n{n-2}\widehat v(t,L)^{\frac2{n-2}}\Phi(t,L)=0
     &\text{on }\R\times\partial\Sph^{n-1}_+.
 \end{cases}
\end{equation}

\begin{proposition}\label{s3:rigidity}
Let $n\geq3$.
\begin{enumerate}[label=\textup{(\roman*)}]
\item Let $v\in\mathscr B_n$. If $\tilde v'(t_0)=0$ at some
$t_0\in\R$, then
\[
 v(t_0+t,\theta)=v(t_0-t,\theta)
 \qquad\text{for every }(t,\theta)\in\R\times\Sph^{n-1}_+.
\]

\item Let $v_1,v_2\in\mathscr B_n$ be even in $t$ and satisfy
$\tilde v_1(0)=\tilde v_2(0)$. Then $v_1=v_2$ on
$\R\times\Sph^{n-1}_+$.

\item Let $v\in\mathscr B_n$ be even in $t$, and let $\Phi$ be a
cylindrically symmetric Jacobi field along $v$ with polynomial
growth in the sense of~\eqref{eq:jacobi-polynomial-space}.
If $\Phi$ is odd in $t$ and $\partial_t\Phi(0,L)=0$, then
$\Phi\equiv0$. If $\Phi$ is even in $t$ and $\Phi(0,L)=0$, then
$\Phi\equiv0$.
\end{enumerate}
\end{proposition}

\begin{proof}
The trace bound in~\eqref{eq:s3-universal-upper} gives
$\widehat v(t,L)\leq((n-2)/2)^{(n-2)/2}$. Thus the boundary
coefficients in~\eqref{eq:s3-secant}
and~\eqref{eq:s3-rp-jacobi-growth} satisfy
\begin{equation}\label{eq:s3-rp-coefficient-bound}
 0<b(t)\leq\frac n2\qquad(t\in\R).
\end{equation}
For $n\geq5$, Corollary~\ref{flatbase:secant} also gives
$|b'|\leq b$ for the coefficients in the three cases below.
For the reflected difference and for the Jacobi equation along
an even solution, \eqref{flatbase:history} gives
$|b'(t)|/b(t)\leq(n-2)t/2$ for $t\geq0$, since $n\geq5$
and $\tanh s\leq s$ for $s\geq0$.

For \textup{(i)}, set
\[
 w_{t_0}(t,\theta)
   =\widehat v(t_0+t,\theta)-\widehat v(t_0-t,\theta).
\]
The two functions solve~\eqref{eq:s3-rp-auxiliary-pair} and
have the same trace bounds as $\widehat v$. The change $t\mapsto-t$
interchanges their traces, and the change $s\mapsto1-s$
in~\eqref{eq:s3-secant} shows that $b$ is even. Since
$\partial_t\widehat v(t_0,L)=0$, we have
\begin{equation}\label{eq:s3-rp-reflected-data}
 w_{t_0}(-t,\theta)=-w_{t_0}(t,\theta),\qquad
 \partial_t w_{t_0}(0,L)=2\partial_t\widehat v(t_0,L)=0.
\end{equation}
By~\eqref{eq:s3-norm-uniform-estimates}, $w_{t_0}$ and its first
derivatives are bounded. The equation for $w_{t_0}$ and
\eqref{eq:s3-rp-reflected-data} give $w_{t_0}=0$ by
Lemma~\ref{flatref:parity-uniqueness}\textup{(i)} if $n\geq5$,
and by Theorem~\ref{nfour:linear-uniqueness}\textup{(i)} if $n=3,4$.
Hence $v$ is symmetric about $t=t_0$.

For \textup{(ii)}, $w=\widehat v_1-\widehat v_2$ and $b$
are even in $t$. The two traces agree at $t=0$, so $w(0,L)=0$.
The function $w$ and its first derivatives are bounded by
\eqref{eq:s3-norm-uniform-estimates}. Applying
Lemma~\ref{flatref:parity-uniqueness}\textup{(ii)} if $n\geq5$
and Theorem~\ref{nfour:linear-uniqueness}\textup{(ii)} if $n=3,4$
gives $w=0$, so $v_1=v_2$.

For \textup{(iii)}, $b$ is even because $v$ is even.
The positive lower trace bound and~\eqref{eq:s3-norm-uniform-estimates}
give uniform bounds for $b$ and its derivatives. The estimates in
Section~\ref{subsec:linearized-spaces} therefore give the same order
of polynomial growth for the local H\"older norms of $\Phi$.
We apply Lemma~\ref{flatref:parity-uniqueness} if $n\geq5$
and Theorem~\ref{nfour:linear-uniqueness} if $n=3,4$,
using the corresponding part for the parity of $\Phi$.
\end{proof}

To locate critical points of the trace, we integrate the
interior equation against $C_n$. For $v\in\mathscr B_n$, set
\begin{equation}\label{eq:s3-moment}
 I(t)=\int_0^L\omega(\beta)C_n(\beta)v(t,\beta)\,\dd\beta.
\end{equation}
If $M$ is an upper bound for $\tilde v$, then
\eqref{eq:s3-canonical-comparison} gives
\[
 0<I(t)\leq M\int_0^L\omega(\beta)C_n(\beta)^2\,\dd\beta
 \qquad\text{for all }t\in\R.
\]
Differentiating under the integral and using the interior equation
in~\eqref{eq:boundary-cylinder}, we obtain
\[
 \begin{aligned}
 I''(t)
 &=\int_0^L C_n
 \left(\frac{(n-2)^2}{4}\omega v
       -\partial_\beta(\omega\partial_\beta v)\right)\,\dd\beta\\
 &=\int_0^L
 \left[v(\omega C_n')'
       -C_n\partial_\beta(\omega\partial_\beta v)\right]\,\dd\beta\\
 &=\bigl[\omega(C_n'v-C_n\partial_\beta v)\bigr]_{\beta=0}^{\beta=L},
 \end{aligned}
\]
where the second equality follows from~\eqref{eq:s3-stationary-angular}
and the last from integration by parts. The endpoint term at
$\beta=0$ vanishes since $\omega(0)=0$.
At $\beta=L$, we have $\omega(L)=C_n(L)=1$ and $C_n'(L)=c$.
The boundary condition in~\eqref{eq:boundary-cylinder} thus gives
\begin{equation}\label{eq:s3-moment-equation}
 I''(t)=c\bigl(\tilde v(t)-\tilde v(t)^{\frac n{n-2}}\bigr).
\end{equation}

\begin{lemma}\label{s3:two-critical}
Let $v\in\mathscr B_n\setminus\{C_n\}$. Then $\tilde v$ takes
values both below and above $1$ and has at least two distinct
critical points.
\end{lemma}

\begin{proof}
If $\tilde v\geq1$ on $\R$, then $I$ is concave by
\eqref{eq:s3-moment-equation}. If $\tilde v\leq1$ on $\R$,
then $I$ is convex. In either case, boundedness on $\R$
forces $I$ to be constant. Since $c\tilde v>0$,
\eqref{eq:s3-moment-equation} gives $\tilde v^{2/(n-2)}\equiv1$.
The trace is therefore identically $1$, and
\eqref{eq:s3-norm-linear-comparison} gives $v=C_n$.
For $v\ne C_n$, there must be $t_-,t_+\in\R$ such that
$\tilde v(t_-)<1<\tilde v(t_+)$. The corresponding conclusion
for the fractional Yamabe equation under positive two-sided
bounds is proved in~\cite[Proposition~2.5]{CDR}.

Suppose that $\tilde v$ is monotone on $[R,\infty)$.
The bounds in~\eqref{eq:s3-bilateral-class} give a limit
$\ell\in(0,\infty)$ as $t\to+\infty$. By
\eqref{eq:s3-moment-equation},
\[
 \lim_{t\to+\infty}I''(t)
 =c\bigl(\ell-\ell^{\frac n{n-2}}\bigr).
\]
If this limit were nonzero, there would be $\delta>0$ and
$\sigma\in\{-1,1\}$ such that $\sigma I''(t)\geq\delta$ for
all $t\geq R$, after increasing $R$. The integral formula
for $I(t)$ based at $R$ then gives
\[
 \begin{aligned}
 \sigma I(t)
 &=\sigma I(R)+(t-R)\sigma I'(R)
   +\int_R^t(t-s)\sigma I''(s)\,\dd s\\
 &\geq \sigma I(R)+(t-R)\sigma I'(R)
   +\frac{\delta}{2}(t-R)^2.
 \end{aligned}
\]
The right-hand side tends to $+\infty$, whereas $I$ is bounded.
Thus $\ell=\ell^{\frac n{n-2}}$, and
$\ell>0$ gives $\ell=1$. If $\tilde v$ is monotone near
$-\infty$, the same argument applied to $I(-t)$ and $\tilde v(-t)$
gives $\tilde v(t)\to1$ as $t\to-\infty$.

If $\tilde v$ had fewer than two critical points, its derivative
would have a fixed strict sign on each tail, and $\tilde v$
would tend to $1$ at both ends. Choose $R$ large enough that
$t_-$ and $t_+$ lie in $(-R,R)$ and the values $\tilde v(\pm R)$
lie strictly between $\tilde v(t_-)$ and $\tilde v(t_+)$. The
minimum and maximum on $[-R,R]$ are then attained at distinct
interior points. Both points are critical, a contradiction.
\end{proof}

The following proposition proves the periodicity conjectured
in~\cite[Section~1]{CDR} at order $1/2$ for solutions in
$\mathscr B_n$ and describes the critical points of their traces.

\begin{proposition}\label{s3:bilateral-profiles}
Let $n\geq3$ and $v\in\mathscr B_n$. If $v\ne C_n$, then $v$ is
periodic in $t$ and has the following properties.
\begin{enumerate}[label=\textup{(\roman*)}]
\item At each critical point $t_0$ of $\tilde v$, we have
$\tilde v''(t_0)\ne0$, and $v$ is symmetric about $t=t_0$.

\item Let $t_0$ be a minimum point of $\tilde v$, and let $t_1$
be the first critical point to its right. Then $T=2(t_1-t_0)$
is a period of $v$. The derivative $\tilde v'$ is positive
between $t_0$ and $t_1$ and negative between $t_1$ and $t_0+T$.

\item The next minimum after $t_0$ occurs at $t_0+T$, and $t_1$
is the unique maximum point of the trace between these minima.
The periods of $v$ and $\tilde v$ are exactly the integer
multiples of $T$.
\end{enumerate}
\end{proposition}

\begin{proof}
For $v\ne C_n$, Lemma~\ref{s3:two-critical} gives two critical
points $\tau_0<\tau_1$ of $\tilde v$. By
Proposition~\ref{s3:rigidity}\textup{(i)}, reflection about
either point preserves $v$. Reflecting first about $\tau_0$
and then about $\tau_1$ gives
\begin{equation}\label{eq:s3-rp-reflection-period}
 v\bigl(t+2(\tau_1-\tau_0),\theta\bigr)=v(t,\theta)
 \qquad\text{for every }(t,\theta)\in\R\times\Sph^{n-1}_+.
\end{equation}

At any critical point $t_0$, the function
$v_{t_0}(t,\theta)=v(t_0+t,\theta)$ is even in $t$. Its derivative
$\Phi=\partial_t v_{t_0}$ is odd and satisfies
\[
 \begin{cases}
 \mathcal A\Phi=0 &\text{in }\R\times\Sph^{n-1}_+,\\
 \mathcal B_{v_{t_0}}\Phi=0
      &\text{on }\R\times\partial\Sph^{n-1}_+,
 \end{cases}
\]
by differentiating~\eqref{eq:boundary-cylinder}. It is bounded
by~\eqref{eq:s3-norm-uniform-estimates} and has polynomial
growth with $q=0$. Since $\Phi(t,L)=\tilde v'(t_0+t)$, a vanishing
second derivative $\tilde v''(t_0)$ would give
$\partial_t\Phi(0,L)=0$. Proposition~\ref{s3:rigidity}\textup{(iii)}
would then give $\Phi\equiv0$, making $v$ independent of $t$
and hence equal to $C_n$. Thus $\tilde v''(t_0)\ne0$.
In particular, every zero of $\tilde v'$ is simple and isolated.

By periodicity, $\tilde v$ attains its minimum at some $t_0$,
and adding a positive period to $t_0$ gives another minimum.
The set of critical points is closed, and $t_0$ is isolated
in it. Hence there is a first critical point $t_1$ to the
right of $t_0$.

Since $\tilde v''(t_0)>0$, the derivative is positive to the
right of $t_0$ for $t$ sufficiently close to $t_0$. It stays
positive between $t_0$ and $t_1$ because there is no critical
point between them. Differentiating
the reflection identity about $t_1$ gives
\[
 \tilde v'(2t_1-t)=-\tilde v'(t).
\]
Consequently the derivative is negative between $t_1$ and
$2t_1-t_0$. The reflections about $t_0$ and $t_1$ give the
period $T=2(t_1-t_0)$ as in~\eqref{eq:s3-rp-reflection-period},
and $2t_1-t_0=t_0+T$.

The trace thus increases from its minimum at $t_0$ to its
maximum at $t_1$ and returns to its minimum at $t_0+T$.
By periodicity, the minimum points are $t_0+jT$ for
integers $j$. A period of $\tilde v$ sends $t_0$ to one of these
points, so it is an integer multiple of $T$. Conversely, every
integer multiple of $T$ is a period of $v$ and of $\tilde v$.
\end{proof}

If the trace of a noncylindrical solution has a minimum at $t_0$,
we use $v(t_0+t,\theta)$ to place this minimum at $t=0$.
By Proposition~\ref{s3:bilateral-profiles}\textup{(iii)},
the difference between any two minimum points is a period
of $v$, so all choices give the same centered solution.
The lower bound in~\eqref{eq:s3-bilateral-class} and the trace
value below $1$ from Lemma~\ref{s3:two-critical} place the
minimum in $(0,1)$. We call it the necksize and write
\begin{equation}\label{eq:s3-necksize}
 \eps=\tilde v(0)=\min_{t\in\R}\tilde v(t).
\end{equation}
By Proposition~\ref{s3:bilateral-profiles}\textup{(i)--(ii)},
the centered solution is even in $t$, and its trace increases
strictly from its minimum at $0$ to its maximum at $T/2$.

\subsection{Parametrization by necksize}
\label{subsec:delaunay-local-chart}
\label{subsec:delaunay-global-family}

We use the minimum of the trace to parametrize centered
solutions, first locally and then on all of $(0,1)$.

Fix a centered noncylindrical solution $v_0\in\mathscr B_n$ with
necksize $\eps_0\in(0,1)$ and fundamental period $T$ as in
Proposition~\ref{s3:bilateral-profiles}. Let $M_T$ be the
compact cylinder obtained from $[0,T]\times\Sph^{n-1}_+$ by
identifying $(0,\theta)$ with $(T,\theta)$ for each
$\theta\in\Sph^{n-1}_+$. We equip $M_T$ with the product metric
and use $t$ for the cylindrical coordinate on its periodic lift.
For $0<\gamma<1$, set
\begin{equation}\label{eq:s3-chart-spaces}
 \begin{aligned}
 \mathscr X_\gamma
   &=C^{2,\gamma}_{T,\mathrm{sym,even}}(\R\times\Sph^{n-1}_+),\\
 \mathscr X_\gamma^+
   &=\{v\in\mathscr X_\gamma:\min_{[0,T]}\tilde v>0\},\\
 \mathscr Y_\gamma^0
   &=C^{0,\gamma}_{T,\mathrm{sym,even}}(\R\times\Sph^{n-1}_+)
     \times C^{1,\gamma}_{T,\mathrm{sym,even}}
        (\R\times\partial\Sph^{n-1}_+),\\
 \mathscr Y_\gamma
   &=\mathscr Y_\gamma^0\times\R.
 \end{aligned}
\end{equation}
The set $\mathscr X_\gamma^+$ is open in $\mathscr X_\gamma$
because restriction to the compact boundary $\partial M_T$ is
continuous in the uniform norm.

The factor $\lambda>0$ allows the period to vary from $T$ to
$\lambda T$. On $M_T$, we define
$\mathfrak F:\mathscr X_\gamma^+\times(0,\infty)\times\R
\to\mathscr Y_\gamma$ by
\begin{equation}\label{eq:s3-augmented-map}
 \mathfrak F(v,\lambda,\eps)=
 \begin{pmatrix}
 -\lambda^{-2}\partial_t^2v-\Delta_{\Sph^{n-1}_+}v
   +\dfrac{(n-2)^2}{4}v\\[3pt]
 \partial_\beta v(\cdot,L)-c\tilde v^{\frac n{n-2}}\\
 \tilde v(0)-\eps
 \end{pmatrix}.
\end{equation}
In particular, $\mathfrak F(v_0,1,\eps_0)=0$.

\begin{proposition}\label{s3:local-neck-chart}
Let $n\geq3$ and $0<\gamma<1$, and let $v_0$, $\eps_0$ and $T$
be as above. There exist an open interval $I\Subset(0,1)$
containing $\eps_0$, an open neighborhood $\mathcal N$ of $(v_0,1)$
in $\mathscr X_\gamma^+\times(0,\infty)$, and a smooth map
$\eps\mapsto(v_\eps^{(T)},\lambda_\eps)$ from $I$ into $\mathcal N$.
At $\eps=\eps_0$ it takes the value $(v_0,1)$, and its graph is
given by
\begin{equation}\label{eq:s3-lc-zero-set}
 \mathfrak F^{-1}(0)\cap(\mathcal N\times I)
 =\{(v_\eps^{(T)},\lambda_\eps,\eps):\eps\in I\}.
\end{equation}
For each $\eps\in I$, define
\begin{equation}\label{eq:s3-lc-physical-chart}
 \begin{aligned}
 v_\eps(t,\theta)&=v_\eps^{(T)}(t/\lambda_\eps,\theta),\\
 T_\eps&=\lambda_\eps T.
 \end{aligned}
\end{equation}
Then $v_\eps\in\mathscr B_n$ is centered and noncylindrical, with
necksize $\eps$. Its trace first returns to its minimum at
$t=T_\eps$, which is the fundamental period of $v_\eps$.

For every integer $k\geq0$, the map
$\eps\mapsto(v_\eps^{(T)},\lambda_\eps)$ is smooth with values
in $C^{k,\gamma}_{T,\mathrm{sym,even}}(\R\times\Sph^{n-1}_+)
\times(0,\infty)$.
\end{proposition}

\begin{proof}
Restriction to $\R\times\partial\Sph^{n-1}_+$ defines a bounded
linear map
\[
 \mathscr X_\gamma\ni v\longmapsto\tilde v
 \in C^{1,\gamma}_{T,\mathrm{sym,even}}
       (\R\times\partial\Sph^{n-1}_+).
\]
For $v$ in a sufficiently small neighborhood of any point of
$\mathscr X_\gamma^+$, the trace values lie in a compact
subinterval of $(0,\infty)$. Taylor's formula and the
multiplication estimate in the boundary H\"older space
therefore give smoothness of
$v\mapsto c\,\tilde v^{n/(n-2)}$. The remaining terms of
$\mathfrak F$ are linear in $(v,\eps)$ with coefficients smooth
in $\lambda>0$, so $\mathfrak F$ is smooth.

Let $\mathcal A$ be the interior operator and $\mathcal B_{v_0}$
the linearized boundary operator at $v_0$, as in
Section~\ref{subsec:linearized-spaces}. The derivative of
$\mathfrak F$ with respect to $(v,\lambda)$ at $(v_0,1,\eps_0)$
is the bounded linear operator
$\mathcal D:\mathscr X_\gamma\times\R\to\mathscr Y_\gamma$ given by
\begin{equation}\label{eq:s3-lc-linearization}
 \mathcal D(z,\mu)
 =
 \bigl(\mathcal Az+2\mu\partial_t^2v_0,\
        \mathcal B_{v_0}z,\ z(0,L)\bigr).
\end{equation}
The kernel of $\mathcal D$ consists of the pairs $(z,\mu)$
satisfying
\begin{equation}\label{eq:s3-lc-kernel-system}
 \begin{aligned}
 \mathcal Az&=-2\mu\partial_t^2v_0,\\
 \mathcal B_{v_0}z&=0,\\
 z(0,L)&=0.
 \end{aligned}
\end{equation}
Differentiating~\eqref{eq:boundary-cylinder} in $t$ shows that
$\partial_t v_0$ is a Jacobi field along $v_0$. The product rule
gives
\begin{equation}\label{eq:s3-lc-time-commutator}
 \begin{aligned}
 \mathcal A(t\partial_t v_0)&=-2\partial_t^2v_0,\\
 \mathcal B_{v_0}(t\partial_t v_0)&=0.
 \end{aligned}
\end{equation}
Thus, for $(z,\mu)\in\ker\mathcal D$, the function
$\Phi=z-\mu t\partial_t v_0$ satisfies the Jacobi equation along
$v_0$. The positivity and smoothness of $\tilde v_0$ give smooth
coefficients in the boundary condition. Starting from
$\Phi\in C^{2,\gamma}_{\mathrm{loc}}$, successive Schauder
estimates give smoothness up to $\R\times\partial\Sph^{n-1}_+$.

Both terms in $\Phi$ are cylindrically symmetric. The function
$z$ is even, and $t\partial_t v_0$ is the product of two odd
functions, so $\Phi$ is even. The correction vanishes at $t=0$,
giving $\Phi(0,L)=z(0,L)=0$. Since $z$ and $\partial_t v_0$ are
continuous and $T$-periodic, they are bounded. Consequently,
\[
 |\Phi(t,\theta)|
 \leq \|z\|_{C^0}+|\mu|\,|t|\,\|\partial_t v_0\|_{C^0}
 \leq C(1+t^2)^{1/2}.
\]
Thus $\Phi$ has polynomial growth of order at most $1$ in
the sense of~\eqref{eq:jacobi-polynomial-space}, and
Proposition~\ref{s3:rigidity}\textup{(iii)} gives $\Phi=0$.
Using $z=\mu t\partial_t v_0$ and periodicity, we obtain
\[
 0=z(t+T,\theta)-z(t,\theta)=\mu T\partial_t v_0(t,\theta).
\]
Since $v_0$ is noncylindrical, $\partial_t v_0$ is nonzero
somewhere. Hence $\mu=0$ and $z=0$.

Consider the Neumann operator
\[
 \mathcal P:\mathscr X_\gamma\to\mathscr Y_\gamma^0,
 \qquad
 \mathcal Pz=(\mathcal Az,\partial_\beta z(\cdot,L)).
\]
For $(f,g)\in\mathscr Y_\gamma^0$, periodicity allows us to
regard $f$ and $g$ as functions on $M_T$ and $\partial M_T$,
respectively. On $M_T$, the associated bilinear form is
\[
 (z,\varphi)\longmapsto
 \int_{M_T}\left(
  \langle\nabla z,\nabla\varphi\rangle
  +\frac{(n-2)^2}{4}z\varphi
 \right)\dd V .
\]
It is bounded and coercive on $H^1(M_T)$ because its diagonal
terms control the $L^2$ norms of $\nabla z$ and $z$. The
interior and boundary data define a bounded linear functional
on $H^1(M_T)$ by the trace theorem. Lax--Milgram gives a unique
weak solution.

Schauder regularity for the Neumann problem gives
$z\in C^{2,\gamma}(M_T)$, since $f\in C^{0,\gamma}(M_T)$ and
$g\in C^{1,\gamma}(\partial M_T)$. The reflection $t\mapsto -t$
and the rotations defining cylindrical symmetry preserve
$\mathcal A$, the normal derivative and both data. By uniqueness,
the solution is invariant under these isometries. Its periodic
lift therefore belongs to $\mathscr X_\gamma$. Thus $\mathcal P$
is bijective, and the bounded inverse theorem makes it an
isomorphism.

The operator $\mathcal P\oplus\operatorname{Id}_{\R}$ is an
isomorphism from $\mathscr X_\gamma\times\R$ to
$\mathscr Y_\gamma^0\times\R=\mathscr Y_\gamma$. By
\eqref{eq:s3-lc-linearization}, the difference between
$\mathcal D$ and this isomorphism is
\[
 (z,\mu)\longmapsto
 \left(
  2\mu\partial_t^2v_0,\
  -\frac n{n-2}c\,\tilde v_0^{\frac2{n-2}}\tilde z,\
  z(0,L)-\mu
 \right).
\]
The first and last components have finite rank. The boundary
component is compact because restriction and multiplication
by the smooth coefficient are bounded into
$C^{2,\gamma}(\partial M_T)$, which embeds compactly into
$C^{1,\gamma}(\partial M_T)$. Hence the difference is compact.
By stability of the Fredholm index under compact
perturbations~\cite[\S~6.4, item~1(c)]{Brezis}, $\mathcal D$
is Fredholm of index zero. Since its kernel is trivial, its
range has codimension zero in $\mathscr Y_\gamma$, so $\mathcal D$
is surjective. The bounded inverse theorem makes $\mathcal D$
an isomorphism. The implicit function theorem
gives~\eqref{eq:s3-lc-zero-set}.

Since $\min_{M_T}v_0>0$, we may shrink $\mathcal N$ and $I$
so that each $v_\eps^{(T)}$ is positive on $M_T$ and close
to $v_0$ in $C^2(M_T)$.

By Proposition~\ref{s3:bilateral-profiles}\textup{(i)--(ii)},
$\tilde v_0''(0)>0$, $\tilde v_0''(T/2)<0$, and $\tilde v_0'$
is positive on $(0,T/2)$. For $v\in\mathscr X_\gamma$,
evenness and periodicity give $\tilde v'(0)=\tilde v'(T/2)=0$.
When $v$ is sufficiently close to $v_0$ in $C^2(M_T)$, the
signs of $\tilde v_0''$ persist near $0$ and $T/2$. Integrating
from these two points gives $\tilde v'>0$ near both ends of
$(0,T/2)$. On the compact interval between those neighborhoods,
$\tilde v_0'$ has a positive lower bound, so $\tilde v'>0$
there as well. Reflection about $T/2$ gives strict decrease on
$(T/2,T)$. Applied to $v_\eps^{(T)}$, this shows that the trace
has its minimum at $0$, where it equals $\eps$ by
\eqref{eq:s3-augmented-map}.

The function $v_\eps$ defined in~\eqref{eq:s3-lc-physical-chart}
solves~\eqref{eq:boundary-cylinder} because
$\mathfrak F(v_\eps^{(T)},\lambda_\eps,\eps)=0$. Since
$\tilde v_\eps\geq\eps>0$, the boundary term
$c\,\tilde v_\eps^{n/(n-2)}$ has the H\"older regularity of
$\tilde v_\eps$. Successive applications of the boundary
Schauder estimates give smoothness up to
$\R\times\partial\Sph^{n-1}_+$. The function $v_\eps$ is positive
and cylindrically symmetric,
and its periodic trace is bounded with minimum $\eps>0$.
Thus $v_\eps\in\mathscr B_n$.

The trace attains its minimum at $0$ and is strictly greater
than $\eps$ on $(0,T_\eps)$. Since $v_\eps$ is $T_\eps$-periodic,
this shows that it is centered and noncylindrical and that
$T_\eps$ is its fundamental period.

It remains to prove smooth dependence in the higher H\"older
norms on $M_T$. Let $\mathcal D_\eps$ be the derivative of
$\mathfrak F$ with respect to $(v,\lambda)$ at
$(v_\eps^{(T)},\lambda_\eps,\eps)$. Smoothness of $\mathfrak F$
and continuity of the family in $\mathscr X_\gamma\times(0,\infty)$
give $\mathcal D_\eps\to\mathcal D$ in operator norm as
$\eps\to\eps_0$. We may therefore restrict $I$ so that
\[
 \bigl\|\mathcal D^{-1}(\mathcal D_\eps-\mathcal D)\bigr\|<1
 \qquad\text{for every }\eps\in I.
\]
The operator $\mathcal D^{-1}\mathcal D_\eps$ is the identity
plus $\mathcal D^{-1}(\mathcal D_\eps-\mathcal D)$, so it is
invertible by the Neumann series. It follows that
$\mathcal D_\eps:\mathscr X_\gamma\times\R\to\mathscr Y_\gamma$
is an isomorphism throughout $I$.

Fix $\eps_*\in I$ and $k\geq2$. For
$v\in C^{k,\gamma}_{T,\mathrm{sym,even}}(\R\times\Sph^{n-1}_+)$
with positive trace, the formula~\eqref{eq:s3-augmented-map}
defines a smooth map with values in
\[
 C^{k-2,\gamma}_{T,\mathrm{sym,even}}(\R\times\Sph^{n-1}_+)
 \times
 C^{k-1,\gamma}_{T,\mathrm{sym,even}}
       (\R\times\partial\Sph^{n-1}_+)
 \times\R.
\]
The smoothness follows from the same composition argument.
Let $f\in C^{k-2,\gamma}(M_T)$ and
$g\in C^{k-1,\gamma}(\partial M_T)$ be even and cylindrically
symmetric, and let $a\in\R$. Since $(f,g,a)\in\mathscr Y_\gamma$,
the isomorphism $\mathcal D_{\eps_*}$ gives a unique pair
$(z,\mu)\in\mathscr X_\gamma\times\R$ satisfying
$\mathcal D_{\eps_*}(z,\mu)=(f,g,a)$. Surjectivity in the
higher H\"older spaces follows once $z\in C^{k,\gamma}(M_T)$.
The interior equation
for $z$ is uniformly elliptic, and its right-hand side differs
from $f$ by a multiple of the smooth function
$\partial_t^2v_{\eps_*}^{(T)}$. The boundary condition expresses
the normal derivative as $g$ plus the trace of $z$ multiplied
by a smooth function. Thus, if $z\in C^{j,\gamma}(M_T)$ for
$2\leq j<k$, the interior right-hand side lies in
$C^{j-1,\gamma}(M_T)$ and the normal derivative lies in
$C^{j,\gamma}(\partial M_T)$. Schauder regularity then gives
$z\in C^{j+1,\gamma}(M_T)$. Starting from $z\in C^{2,\gamma}(M_T)$
gives $z\in C^{k,\gamma}(M_T)$ by iteration. Injectivity
follows from uniqueness in $\mathscr X_\gamma\times\R$, and
the bounded inverse theorem gives an isomorphism in the
higher H\"older spaces.

The implicit function theorem now gives a smooth family in
$C^{k,\gamma}_{T,\mathrm{sym,even}}(\R\times\Sph^{n-1}_+)
\times(0,\infty)$ near $\eps_*$. It passes through
$(v_{\eps_*}^{(T)},\lambda_{\eps_*})$ and, by continuous
embedding into $\mathscr X_\gamma\times(0,\infty)$, lies in
$\mathcal N$ for $\eps$ close to $\eps_*$. The graph
identity~\eqref{eq:s3-lc-zero-set} identifies the two families.
Since $\eps_*$ was arbitrary, this proves smooth dependence on
$I$ for $k\geq2$. The continuous inclusions from $C^{2,\gamma}$
into $C^{k,\gamma}$ give the cases $k=0,1$.
\end{proof}

The local parametrization gives openness of the set of
realized necksizes, and the uniform estimates give relative
closedness in $(0,1)$.

\begin{restatedtheorem}\label{s3:global-family}
Let $n\geq3$. For every $\eps\in(0,1)$ there is a unique solution
$v_\eps\in\mathscr B_n$ satisfying
\begin{equation}\label{eq:s3-centered-data}
 \tilde v_\eps(0)=\min_{t\in\R}\tilde v_\eps(t)=\eps.
\end{equation}
The solution $v_\eps$ is noncylindrical, even in $t$, and
periodic. Denote its fundamental period by $T_\eps$.

\begin{enumerate}[label=\textup{(\roman*)}]
\item Let $M_n$ be the upper bound in~\eqref{eq:s3-universal-upper}.
The trace attains its maximum at $T_\eps/2$, and this maximum
is strictly greater than $1$. For every $t\in\R$,
\begin{equation}\label{eq:s3-family-bounds}
 \eps\leq\tilde v_\eps(t)\leq\tilde v_\eps(T_\eps/2)\leq M_n.
\end{equation}
The corresponding bounds on the cylinder are
\[
 \eps C_n(\theta)\leq v_\eps(t,\theta)\leq M_n C_n(\theta)
 \qquad\text{for all }(t,\theta)\in\R\times\Sph^{n-1}_+.
\]

\item If $t_1$ is the first positive critical point of
$\tilde v_\eps$, then
\begin{equation}\label{eq:s3-fundamental-period}
 T_\eps=2t_1.
\end{equation}
The periods of $v_\eps$ and $\tilde v_\eps$ are precisely the
integer multiples of $T_\eps$. The critical points of
$\tilde v_\eps$ consist of the minimum points $t_j^{\min}$ and
the maximum points $t_j^{\max}$, given by
\begin{equation}\label{eq:s3-critical-point-sets}
 \begin{aligned}
 t_j^{\min}&=jT_\eps,\\
 t_j^{\max}&=\left(j+\tfrac12\right)T_\eps,
 \end{aligned}
 \qquad j\in\mathbb Z.
\end{equation}
We have $\tilde v_\eps''(0)>0$ and $\tilde v_\eps''(T_\eps/2)<0$.
Moreover,
\begin{equation}\label{eq:s3-half-period-monotonicity}
 \tilde v_\eps'(t)=-\tilde v_\eps'(T_\eps-t)>0
 \qquad\text{for }0<t<T_\eps/2.
\end{equation}

\item The only cylindrical solution in $\mathscr B_n$ is
$C_n$, and we set $v_1=C_n$. Every noncylindrical solution
$v\in\mathscr B_n$ has the form
\begin{equation}\label{eq:s3-orbit-classification}
 v(t,\theta)=v_\eps(t+\tau,\theta)
\end{equation}
for a unique $\eps\in(0,1)$ and some $\tau\in\R$. Two choices
of $\tau$ give the same solution exactly when their difference
is an integer multiple of $T_\eps$.

\item Define the $1$-periodic function $v_\eps^{(1)}$ by
\[
 v_\eps^{(1)}(s,\theta)=v_\eps(T_\eps s,\theta).
\]
For every integer $k\geq0$ and $0<\gamma<1$, the map
\begin{equation}\label{eq:s3-fixed-period-smoothness}
 \begin{aligned}
 (0,1)&\longrightarrow
 (0,\infty)\times
 C^{k,\gamma}_{1,\mathrm{sym,even}}(\R\times\Sph^{n-1}_+),\\
 \eps&\longmapsto(T_\eps,v_\eps^{(1)})
 \end{aligned}
\end{equation}
is smooth.
\end{enumerate}
\end{restatedtheorem}

\begin{proof}
Let
\begin{equation}\label{eq:s3-gc-realized-set}
 E=\bigl\{\eps\in(0,1):
 \tilde v(0)=\min_{t\in\R}\tilde v(t)=\eps
 \text{ for some }v\in\mathscr B_n\bigr\}.
\end{equation}
To show that $E$ is nonempty, we construct a periodic solution
on $M_T$, with $T>0$ to be chosen. Set $q=2(n-1)/(n-2)$
and minimize
\[
 Q_T(V)=
 \frac{\displaystyle\int_{M_T}
     \left(|\nabla V|^2+\frac{(n-2)^2}{4}V^2\right)\dd\mu}
 {\displaystyle\left(\int_{\partial M_T}
                           |V|^q\dd S\right)^{2/q}}
\]
over cylindrically symmetric functions in $H^1(M_T)$ with
nonzero trace. Here $\dd\mu$ is the product volume element
and $\dd S$ is the induced boundary measure. In this class
the trace depends only on $t$. The trace theorem gives a trace in
$H^{1/2}$ of the circle of length $T$, which embeds compactly
in $L^q$ for every finite $q$. Normalizing the denominator
to $1$ gives a bounded minimizing sequence in $H^1(M_T)$.
Its weak limit retains the trace norm by compactness and
minimizes $Q_T$ by lower semicontinuity.

Taking the absolute value gives a nonnegative minimizer.
Its Euler--Lagrange equation is the interior equation
in~\eqref{eq:boundary-cylinder} with
$\partial_\beta V=\lambda\tilde V^{q-1}$ on $\partial M_T$,
where $\lambda>0$ is the minimum of $Q_T$ under the trace
normalization. Averaging any test function over the rotations
defining cylindrical symmetry gives the weak equation for all
test functions. The trace belongs to $L^r(\partial M_T)$ for every finite $r$,
and so does $\lambda\tilde V^{q-1}$. For $r>n$, duality gives
$L^r(\partial M_T)\hookrightarrow W^{-1/r,r}(\partial M_T)$,
and Neumann regularity gives $V\in W^{1,r}(M_T)$.
By Sobolev embedding, $\tilde V$ is H\"older continuous.
Since $q>2$, the normal derivative $\lambda\tilde V^{q-1}$
has the same H\"older exponent. The Schauder estimates then
give $V\in C^{1,\gamma}(M_T)$ for some $0<\gamma<1$. The strong maximum principle and
the boundary point lemma give positivity on $M_T$, since
a zero on the boundary would have zero normal derivative.
The boundary condition then gives smoothness up to $\partial M_T$.
Multiplication by $(\lambda/c)^{1/(q-2)}$ makes its boundary
condition the one in~\eqref{eq:boundary-cylinder}.

Among functions independent of $t$, the least value of
$Q_T$ is $Q_T(C_n)$. Indeed, for a fixed constant trace,
its multiple of $C_n$ minimizes the numerator by the
interior equation and integration by parts. For
$V_s=C_n(1+s\cos(2\pi t/T))$, expansion of the two
integrals gives
\[
 \frac{Q_T(V_s)}{Q_T(C_n)}
 =1+\frac{s^2}{2}
   \left[
    \frac{(2\pi/T)^2}{c}\int_0^L\omega C_n^2\dd\beta
    -(q-2)
   \right]+O(s^3).
\]
Here we used~\eqref{eq:s3-stationary-angular} to compute
$\int_0^L\omega\bigl((C_n')^2+\frac{(n-2)^2}{4}C_n^2\bigr)
\dd\beta=c$. The coefficient of $s^2$ is negative for
sufficiently large $T$. For such $T$ the minimizing solution
depends on $t$ and its periodic lift belongs to $\mathscr B_n$.
Its trace is nonconstant, by uniqueness of the bounded
extension in~\eqref{eq:s3-norm-linear-comparison}.

Integrating~\eqref{eq:s3-moment-equation} over a period gives
$\int_0^T(\tilde v-\tilde v^{n/(n-2)})\dd t=0$ for this
normalized solution $v$. Since the trace is positive and
nonconstant, it takes values below and above $1$. Choosing
the origin of $t$ at a minimum thus gives a necksize
$\eps_0\in E$.

Now fix $\eps_0\in E$ and choose $v_0\in\mathscr B_n$ whose
trace attains its minimum $\eps_0$ at $0$. Since $\eps_0<1$
and the only cylindrical solution has trace $1$, $v_0$ is
noncylindrical. Proposition~\ref{s3:local-neck-chart}, applied
at $v_0$, gives an open interval $I\subset(0,1)$ containing
$\eps_0$. For each $\eps\in I$, it gives a solution in
$\mathscr B_n$ whose trace attains its minimum $\eps$ at $0$.
Thus $I\subset E$, proving that $E$ is open in $(0,1)$.

Let $(\eps_j)$ be a sequence in $E$ converging to
$\eps_\infty\in(0,1)$. Choose $v_j\in\mathscr B_n$ whose trace
attains its minimum $\eps_j$ at $0$. For all sufficiently
large $j$, we have $\eps_j\geq\eps_\infty/2$, so
$\tilde v_j(t)\geq\eps_\infty/2$ for every $t\in\R$. The
uniform upper bound~\eqref{eq:s3-universal-upper} gives
$\tilde v_j(t)\leq M_n$. Applying the extension bounds
in~\eqref{eq:s3-canonical-comparison} with these two constants
gives
\begin{equation}\label{eq:s3-lc-uniform-extension-bounds}
 \tfrac12\eps_\infty C_n(\theta)
 \leq v_j(t,\theta)\leq M_n C_n(\theta)
 \qquad\text{for all }(t,\theta)\in\R\times\Sph^{n-1}_+.
\end{equation}
The uniform estimates in~\eqref{eq:s3-norm-uniform-estimates}
and a diagonal argument give convergence, along a subsequence,
in $C^k([-R,R]\times\Sph^{n-1}_+)$ for every $k\geq0$ and
$R>0$ to a solution $v_\infty$ of~\eqref{eq:boundary-cylinder}.
The bounds and cylindrical
symmetry pass to the limit, so $v_\infty\in\mathscr B_n$.
Since $\tilde v_j(t)\geq\tilde v_j(0)=\eps_j$ for all $t$,
passing to the limit gives
$\tilde v_\infty(t)\geq\tilde v_\infty(0)=\eps_\infty$ for all
$t$. Hence $\eps_\infty\in E$, and $E$ is relatively closed
in $(0,1)$. By connectedness,
\begin{equation}\label{eq:s3-gc-open-closed-conclusion}
 E=(0,1).
\end{equation}

For any solution satisfying~\eqref{eq:s3-centered-data}, the
trace has derivative zero at $0$, so
Proposition~\ref{s3:rigidity}\textup{(i)} gives evenness in $t$.
Two such solutions for the same $\eps$ have the same trace
value at $0$ and hence agree by
Proposition~\ref{s3:rigidity}\textup{(ii)}. Denote this unique
solution by $v_\eps$. It is noncylindrical because $\eps<1$.

Proposition~\ref{s3:bilateral-profiles} gives the periods and
the critical points in \textup{(ii)}, including the signs of
the second derivatives at $0$ and $T_\eps/2$. It also gives
$\tilde v_\eps'(t)>0$ for $0<t<T_\eps/2$. Evenness and
periodicity imply $\tilde v_\eps(T_\eps-t)=\tilde v_\eps(t)$.
Differentiating this identity gives
\eqref{eq:s3-half-period-monotonicity}.

For \textup{(i)}, the maximum at $T_\eps/2$ is greater than
$1$ by Lemma~\ref{s3:two-critical} and at most $M_n$ by
\eqref{eq:s3-universal-upper}. Applying
\eqref{eq:s3-canonical-comparison} with trace bounds $\eps$
and $M_n$ gives the bounds for $v_\eps$ on the cylinder.

Section~\ref{subsec:delaunay-normalization} identifies $C_n$
as the only cylindrical solution. For a noncylindrical solution
$v\in\mathscr B_n$, choose $t_0$ at which its trace attains
its minimum and set $\eps=\tilde v(t_0)$. The function
$(t,\theta)\mapsto v(t+t_0,\theta)$ has trace minimum $\eps$
at $0$, so by uniqueness it equals $v_\eps$. It follows that
$v(t,\theta)=v_\eps(t-t_0,\theta)$, which gives
\eqref{eq:s3-orbit-classification} with $\tau=-t_0$. The
parameter $\eps$ is determined by the minimum of $\tilde v$.
If two values $\tau_1$ and $\tau_2$ give the same solution,
the substitution $s=t+\tau_2$ shows that $\tau_1-\tau_2$ is a
period of $v_\eps$. Conversely, any such period gives equality.
By \textup{(ii)}, these periods are precisely the integer
multiples of $T_\eps$. This proves \textup{(iii)}.

For \textup{(iv)}, fix $\eps_0\in(0,1)$ and set
$T_0=T_{\eps_0}$. Proposition~\ref{s3:local-neck-chart},
applied at $v_{\eps_0}$, gives the smooth map
$\eps\mapsto(v_\eps^{(T_0)},\lambda_\eps)$ near $\eps_0$.
The function $v_\eps^{(T_0)}(t/\lambda_\eps,\theta)$ has trace
minimum $\eps$ at $0$. By uniqueness it equals $v_\eps$,
and~\eqref{eq:s3-lc-physical-chart} gives
$T_\eps=\lambda_\eps T_0$. Consequently,
\begin{equation}\label{eq:s3-gc-rescaled-local-chart}
 v_\eps^{(1)}(s,\theta)
 =v_\eps(T_\eps s,\theta)
 =v_\eps^{(T_0)}(T_0s,\theta).
\end{equation}
Since $T_0$ is fixed, the substitution $t=T_0s$ gives a bounded
linear map from the periodic H\"older space of period $T_0$
to the corresponding space of period $1$. The local map
therefore gives smooth dependence of $(T_\eps,v_\eps^{(1)})$
near $\eps_0$ in~\eqref{eq:s3-fixed-period-smoothness}. Since
$\eps_0$ was arbitrary, this proves \textup{(iv)}.
\end{proof}

Smooth dependence in Theorem~\ref{s3:global-family}\textup{(iv)}
gives compactness on the cylinder of period $1$ for necksizes
in a compact subinterval of $(0,1)$.

\begin{lemma}\label{s3:period-compactness}
Let $0<\eps_-\leq\eps_+<1$ and let $(\eps_j)$ be a sequence in
$[\eps_-,\eps_+]$. Write $v_j=v_{\eps_j}$ and $T_j=T_{\eps_j}$,
as in Theorem~\ref{s3:global-family}. After passing to a
subsequence, $\eps_j\to\eps_\infty$ for some
$\eps_\infty\in[\eps_-,\eps_+]$. With $v_\infty=v_{\eps_\infty}$
and $T_\infty=T_{\eps_\infty}$, we have
\begin{equation}\label{eq:s3-lc-compactness-conclusions}
 \begin{gathered}
 T_j\longrightarrow T_\infty,\\
 v_j\longrightarrow v_\infty
 \quad\text{in }C^{k,\gamma}_{\mathrm{loc}}
       (\R\times\Sph^{n-1}_+).
 \end{gathered}
\end{equation}
On the cylinder of period $T_\infty$, we also have
\begin{equation}\label{eq:s3-lc-rescaled-compactness}
 \begin{gathered}
 v_{\eps_j}^{(1)}(t/T_\infty,\theta)
 =v_j\bigl((T_j/T_\infty)t,\theta\bigr)
 \longrightarrow v_\infty(t,\theta)\\
 \text{in }C^{k,\gamma}_{T_\infty,\mathrm{sym,even}}
       (\R\times\Sph^{n-1}_+).
 \end{gathered}
\end{equation}
Both function convergences hold along the same subsequence
for every integer $k\geq0$ and every $0<\gamma<1$.
\end{lemma}

\begin{proof}
By compactness of $[\eps_-,\eps_+]$, we may pass to a
subsequence for which $\eps_j\to\eps_\infty$ with
$\eps_\infty\in[\eps_-,\eps_+]$. The continuity of the map
in~\eqref{eq:s3-fixed-period-smoothness} gives
\[
 T_j\longrightarrow T_\infty
 \qquad\text{and}\qquad
 v_{\eps_j}^{(1)}\longrightarrow v_{\eps_\infty}^{(1)}
 \quad\text{in }C^{k,\gamma}_{1,\mathrm{sym,even}}
       (\R\times\Sph^{n-1}_+)
\]
for every integer $k\geq0$ and every $0<\gamma<1$, along the
subsequence chosen using only $(\eps_j)$. Composition with the
fixed dilation $s=t/T_\infty$ preserves the symmetries and is
bounded between the periodic H\"older spaces of periods $1$
and $T_\infty$. Each $t$-derivative contributes a factor of
$T_\infty^{-1}$, and the H\"older quotients are controlled by
the Lipschitz constant of the dilation in product coordinates.
Since $v_{\eps_\infty}^{(1)}(t/T_\infty,\theta)=v_\infty(t,\theta)$,
this gives~\eqref{eq:s3-lc-rescaled-compactness}.

For the local convergence, fix $R>0$ and write
$K_R=[-R,R]\times\Sph^{n-1}_+$. For $|t|\leq R$ and large $j$,
both $t/T_j$ and $t/T_\infty$ belong to
$[-2R/T_\infty,2R/T_\infty]$. The derivatives through order
$k+1$ of $v_{\eps_j}^{(1)}$ converge uniformly on the cylinder
over this interval, and $t/T_j\to t/T_\infty$ uniformly.
In $v_j(t,\theta)=v_{\eps_j}^{(1)}(t/T_j,\theta)$, each
$t$-derivative contributes a factor $T_j^{-1}$. The chain rule
and uniform continuity of the limiting derivatives give
$v_j\to v_\infty$ in $C^{k+1}(K_R)$. On $K_R$, the mean value
estimate in product charts bounds the H\"older seminorms of
the derivatives of order $k$ by the $C^{k+1}$ norm. This proves
the local convergence in~\eqref{eq:s3-lc-compactness-conclusions}.
\end{proof}
\stepcounter{theorem}

\subsection{The Hamiltonian and the limiting periods}
\label{subsec:delaunay-hamiltonian}

We record the conserved Hamiltonian of the Delaunay family
and evaluate it at $C_n$ and at a trace minimum. The even
Jacobi field at $C_n$ determines the limit of the fundamental
period as $\eps\to1$. As $\eps\to0$, the period grows
logarithmically.

For $v\in\mathscr B_n$, define the Hamiltonian divided by
the area of $\Sph^{n-2}$ by
\begin{equation}\label{eq:s3-hamiltonian}
 \begin{aligned}
 H_v(t)
 ={}&\frac12\int_0^L\omega(\beta)
       \left(v_t(t,\beta)^2-v_\beta(t,\beta)^2
                   -\frac{(n-2)^2}{4}v(t,\beta)^2\right)\dd\beta\\
    &+\frac{c(n-2)}{2(n-1)}
                \tilde v(t)^{\frac{2(n-1)}{n-2}}.
 \end{aligned}
\end{equation}
This is the cylindrical form of the Hamiltonian of DelaTorre
and Gonz\'alez~\cite[Theorem~1.2]{DG} at order $1/2$.
Differentiating in $t$ and integrating by parts in $\beta$ gives
\begin{equation}\label{eq:s3-he-energy-derivative}
 \begin{aligned}
 H_v'(t)
 ={}&\int_0^L\omega v_t
       \left(v_{tt}+\omega^{-1}(\omega v_\beta)_\beta
                         -\frac{(n-2)^2}{4}v\right)\dd\beta\\
    &-[\omega v_\beta v_t]_0^L
       +c\,\tilde v(t)^{\frac n{n-2}}\tilde v'(t).
 \end{aligned}
\end{equation}
The integral vanishes by the interior equation
in~\eqref{eq:boundary-cylinder}, and the term at $\beta=0$
vanishes since $\omega(0)=0$. At $\beta=L$, we have
$\omega(L)=1$ and $v_t(t,L)=\tilde v'(t)$, so the boundary
condition in~\eqref{eq:boundary-cylinder} gives
\begin{equation}\label{eq:s3-hamiltonian-conservation}
 H_v'(t)
 =\bigl(c\,\tilde v(t)^{\frac n{n-2}}-v_\beta(t,L)\bigr)
      \tilde v'(t)=0.
\end{equation}
We denote this constant by $H(v)$. For the extension $b_n v$,
the Hamiltonian with integrals over $\Sph^{n-1}_+$ and its
boundary is $|\Sph^{n-2}|b_n^2H(v)$. The factor $b_n^2$ in the
boundary term follows from $\frac{n-2}{2}b_n^{2/(n-2)}=c$.

Multiplying the equation for $C_n$ in
\eqref{eq:s3-stationary-angular} by $C_n$ and integrating by
parts gives
\begin{equation}\label{eq:s3-he-stationary-dirichlet-energy}
 \int_0^L\omega
       \left((C_n')^2+\frac{(n-2)^2}{4}C_n^2\right)\dd\beta
 =[\omega C_n C_n']_0^L=c,
\end{equation}
where we used $C_n(L)=1$ and $C_n'(L)=c$. Hence
\eqref{eq:s3-hamiltonian} gives
\begin{equation}\label{eq:s3-stationary-energy}
 H(C_n)=-\frac{c}{2(n-1)}.
\end{equation}

For $0<\eps<1$, write
\begin{equation}\label{eq:s3-he-neck-decomposition}
 v_\eps(0,\beta)=\eps C_n(\beta)+R_\eps(\beta),
 \qquad R_\eps(L)=0.
\end{equation}
By~\eqref{eq:s3-stationary-angular} and the product rule,
\begin{equation}\label{eq:s3-he-neck-cross-term}
 \begin{aligned}
 \int_0^L\omega
   \left(C_n'R_\eps'+\frac{(n-2)^2}{4}C_nR_\eps\right)\dd\beta
 &=\int_0^L(\omega C_n'R_\eps)'\dd\beta\\
 &=[\omega C_n'R_\eps]_0^L=0.
 \end{aligned}
\end{equation}
The boundary terms vanish because $\omega(0)=0$ and
$R_\eps(L)=0$. By evenness in $t$, $\partial_t v_\eps(0,\theta)=0$
for $\theta\in\Sph^{n-1}_+$. Thus~\eqref{eq:s3-hamiltonian}
at $t=0$ and identities~\eqref{eq:s3-he-stationary-dirichlet-energy}
and~\eqref{eq:s3-he-neck-cross-term} give
\begin{equation}\label{eq:s3-neck-energy}
 \begin{aligned}
 H(v_\eps)
 ={}&-\frac{c\eps^2}{2}
       \left(1-\frac{n-2}{n-1}\eps^{\frac2{n-2}}\right)\\
    &-\frac12\int_0^L\omega
          \left((R_\eps')^2+\frac{(n-2)^2}{4}R_\eps^2\right)
          \dd\beta<0.
 \end{aligned}
\end{equation}
The factor in parentheses is positive for $0<\eps<1$.

The limiting period at $C_n$ is determined by the
Dirichlet-to-Neumann operator associated with the interior
equation in~\eqref{eq:boundary-cylinder}.

For $\xi\in\R$, let $A_\xi$ be the smooth cylindrically
symmetric solution of
\begin{equation}\label{eq:s3-he-angular-symbol-problem}
 (\omega A_\xi')'
 =\left(\frac{(n-2)^2}{4}+\xi^2\right)\omega A_\xi,
 \qquad A_\xi(L)=1.
\end{equation}
Existence, uniqueness and smooth dependence on $\xi$ follow
from the Dirichlet problem for
$-\Delta_{\Sph^{n-1}_+}+(n-2)^2/4+\xi^2$ on
$\Sph^{n-1}_+$ with boundary value $1$ on
$\partial\Sph^{n-1}_+$.

The bounded extension of $e^{i\xi t}$ that solves the interior
equation in~\eqref{eq:boundary-cylinder} is
$e^{i\xi t}A_\xi(\beta)$. By the conormal convention of
Section~\ref{subsec:notation},
\begin{equation}\label{eq:s3-he-mode-conormal}
 \left.\partial_\beta
       \bigl(e^{i\xi t}A_\xi(\beta)\bigr)\right|_{\beta=L}
 =A_\xi'(L)e^{i\xi t}.
\end{equation}
The Dirichlet-to-Neumann operator therefore sends $e^{i\xi t}$
to $\Theta_n(\xi)e^{i\xi t}$, where
$\Theta_n(\xi)=A_\xi'(L)$. In boundary dimension $n-1$ and at
order $1/2$, the formula of DelaTorre and
Gonz\'alez~\cite[Theorem~1.1]{DG} for functions independent of
the spherical variable gives
\begin{equation}\label{eq:s3-symbol}
 \Theta_n(\xi)
 =2\left|\frac{\Gamma(\frac n4+\frac{i\xi}2)}
                 {\Gamma(\frac{n-2}{4}+\frac{i\xi}2)}\right|^2,
 \qquad \Theta_n(0)=c.
\end{equation}
Here $A_0=C_n$, which also gives the stated value at $\xi=0$.

The equation in~\eqref{eq:s3-he-angular-symbol-problem}
depends on $\xi$ through $\xi^2$, so uniqueness gives
$A_{-\xi}=A_\xi$ and $\Theta_n(-\xi)=\Theta_n(\xi)$.
Green's identity for $A_\xi$ and $A_\eta$ gives
\[
 \Theta_n(\xi)-\Theta_n(\eta)
 =(\xi^2-\eta^2)\int_0^L\omega A_\xi A_\eta\,\dd\beta.
\]
Here $A_\xi(L)=A_\eta(L)=\omega(L)=1$ and the term at
$\beta=0$ vanishes since $\omega(0)=0$. Dividing by
$\xi-\eta$ and letting $\eta\to\xi$, we obtain
\begin{equation}\label{eq:s3-he-symbol-derivative}
 \Theta_n'(\xi)=2\xi\int_0^L\omega A_\xi^2\,\dd\beta>0
 \qquad\text{for }\xi>0.
\end{equation}

Applying Stirling's formula to the quotient in
\eqref{eq:s3-symbol} gives $\Theta_n(\xi)/|\xi|\to1$ as
$|\xi|\to\infty$. Since $\Theta_n$ is strictly increasing on
$(0,\infty)$ and $\Theta_n(0)=c<\frac n{n-2}c$, there is a
unique $\xi_*>0$ such that
\begin{equation}\label{eq:s3-critical-period}
 \Theta_n(\xi_*)=\frac n{n-2}c,
 \qquad T_*=\frac{2\pi}{\xi_*}.
\end{equation}

For $\xi>0$, set $Z_\xi(t,\beta)=A_\xi(\beta)\cos(\xi t)$.
By~\eqref{eq:s3-he-angular-symbol-problem}, $Z_\xi$ satisfies
the interior equation $\mathcal A Z_\xi=0$. Using $A_\xi(L)=1$
and taking the real part in~\eqref{eq:s3-he-mode-conormal}
gives its trace and conormal derivative on
$\R\times\partial\Sph^{n-1}_+$,
\begin{equation}\label{eq:s3-he-periodic-multiplier}
 Z_\xi(t,L)=\cos(\xi t),\qquad
 \partial_\beta Z_\xi(t,L)=\Theta_n(\xi)\cos(\xi t).
\end{equation}
The boundary condition linearized at $C_n$, whose trace is
$1$, is $\partial_\beta Z(t,L)=\frac{n}{n-2}c\,Z(t,L)$.
Thus $Z_\xi$ satisfies this condition precisely when
$\Theta_n(\xi)=\frac{n}{n-2}c$. By~\eqref{eq:s3-critical-period},
this occurs at $\xi=\xi_*$. The Jacobi field $Z_{\xi_*}$ has
fundamental period $T_*=2\pi/\xi_*$, as is seen from its trace.
This is the period obtained in~\cite[Theorem~1.3]{DG} for
the linearized equation at order $1/2$ in boundary dimension
$n-1$ when the trace is independent of the spherical variable.

As $\eps\to1$, we compare the trace of $Z_{\xi_*}$ with
that of the difference $v_\eps-C_n$ after normalization at
an extremum where $|\tilde v_\eps-1|$ is largest. The bound
for $|\partial_t\widehat v(t,L)|/\widehat v(t,L)$ in~\eqref{flatbase:gradient}
gives the divergence of the period as $\eps\to0$.

\begin{proposition}\label{s3:endpoint-periods}
Let $n\geq3$ and let $(v_\eps,T_\eps)_{0<\eps<1}$ be the
family in Theorem~\ref{s3:global-family}.
\begin{enumerate}[label=\textup{(\roman*)}]
\item For every integer $k\geq0$ and every $0<\gamma<1$,
\begin{equation}\label{eq:s3-stationary-endpoint-limit}
 \lim_{\eps\to1}
 \|v_\eps-C_n\|_{C^{k,\gamma}_{\mathrm{uloc}}
                     (\R\times\Sph^{n-1}_+)}=0.
\end{equation}

\item Let $\xi_*>0$ be the unique solution of
$\Theta_n(\xi_*)=\frac{n}{n-2}c$ and set
$T_*=\frac{2\pi}{\xi_*}$. Then
\begin{equation}\label{eq:s3-period-endpoint-limits}
 \lim_{\eps\to1}T_\eps=T_*,
 \qquad
 \lim_{\eps\to0}T_\eps=+\infty.
\end{equation}
\end{enumerate}
\end{proposition}

\begin{proof}
Let $(\eps_j)$ be a sequence in $(0,1)$ converging to $1$,
and write $v_j=v_{\eps_j}$ and $T_j=T_{\eps_j}$. For all
sufficiently large $j$, Theorem~\ref{s3:global-family} gives
\begin{equation}\label{eq:s3-he-stationary-uniform-bounds}
 \begin{aligned}
 \tfrac12\leq\eps_j\leq\tilde v_j(t)&\leq M_n,\\
 \tfrac12 C_n(\theta)\leq v_j(t,\theta)&\leq M_n C_n(\theta)
 \end{aligned}
 \qquad\text{for }(t,\theta)\in\R\times\Sph^{n-1}_+.
\end{equation}
Set $a_j=\|\tilde v_j-1\|_{L^\infty(\R)}$. The minimum of
$\tilde v_j$ is $\eps_j$ at $t=0$, and its maximum is attained
at $t=T_j/2$. Thus
\[
 a_j=\max\{1-\eps_j,\ \tilde v_j(T_j/2)-1\}.
\]
Since $1-\eps_j\to0$, if $a_j$ did not tend to zero, there
would be a subsequence with $\tilde v_j(T_j/2)\geq1+\delta$
for some $\delta>0$. The functions $v_j(t+T_j/2,\theta)$
satisfy~\eqref{eq:s3-he-stationary-uniform-bounds}. By
Proposition~\ref{prop:flat-uniform-estimates}, a further
subsequence converges in $C^k([-R,R]\times\Sph^{n-1}_+)$
for every $k\geq0$ and $R>0$ to a solution
$v_\infty\in\mathscr B_n$.
The lower bound $\tilde v_j\geq\eps_j$ passes to the limit
as $\tilde v_\infty\geq1$, while the bound at the maxima
gives $\tilde v_\infty(0)\geq1+\delta$. Since the trace of
$C_n$ is $1$, Lemma~\ref{s3:two-critical} requires
$\tilde v_\infty$ to take a value below $1$. This contradiction
shows that $a_j\to0$.

For large $j$, the trace bounds
$0<1-a_j\leq\tilde v_j\leq1+a_j$ give, by
\eqref{eq:s3-canonical-comparison},
\begin{equation}\label{eq:s3-he-stationary-extension-comparison}
 |v_j(t,\theta)-C_n(\theta)|\leq a_j C_n(\theta)
 \qquad\text{on }\R\times\Sph^{n-1}_+.
\end{equation}
The difference $D_j=v_j-C_n$ satisfies $\mathcal AD_j=0$,
where $\mathcal A$ is defined in~\eqref{eq:jacobi-local-pair}.
Subtracting the boundary conditions for $v_j$ and $C_n$
gives
\[
 \partial_\beta D_j(t,L)=q_j(t)D_j(t,L),
 \qquad
 q_j(t)=\frac n{n-2}c\int_0^1
       \bigl(1-s+s\tilde v_j(t)\bigr)^{\frac2{n-2}}\,\dd s.
\]
For $0\leq s\leq1$, the convex combination
$1-s+s\tilde v_j(t)$ lies in $[1/2,M_n]$ by
\eqref{eq:s3-he-stationary-uniform-bounds}. The function
$x\mapsto x^{2/(n-2)}$ is smooth on this interval. The
uniform regularity bounds for $\tilde v_j$ in
Proposition~\ref{prop:flat-uniform-estimates} therefore bound
$q_j$ in $C^{k,\gamma}_{\mathrm{uloc}}(\R)$ for every $k$
and $\gamma$.

For each $\tau\in\R$, the local Schauder estimate of
Section~\ref{subsec:linearized-spaces} bounds the
$C^{k,\gamma}$ norm of $D_j$ on
$[\tau-1,\tau+1]\times\Sph^{n-1}_+$ by its $L^\infty$ norm
on $[\tau-2,\tau+2]\times\Sph^{n-1}_+$. The constant is
independent of $j$ and $\tau$ by the uniform bounds on $q_j$
and the invariance of $\mathcal A$ under translation in $t$.
Taking the supremum over $\tau$ and using
\eqref{eq:s3-he-stationary-extension-comparison}, we obtain
\begin{equation}\label{eq:s3-he-stationary-schauder}
 \begin{aligned}
 \|D_j\|_{C^{k,\gamma}_{\mathrm{uloc}}
                   (\R\times\Sph^{n-1}_+)}
 &\leq C_{k,\gamma}\|D_j\|_{L^\infty(\R\times\Sph^{n-1}_+)}\\
 &\leq C'_{k,\gamma}a_j\longrightarrow0.
 \end{aligned}
\end{equation}
This proves \textup{(i)}.

To determine the limit of $T_j$, set $s_j=0$ if
$1-\eps_j\geq\tilde v_j(T_j/2)-1$ and $s_j=T_j/2$ otherwise.
Then $|\tilde v_j(s_j)-1|=a_j>0$. Define
\begin{equation}\label{eq:s3-he-normalized-differences}
 Z_j(t,\theta)=
 \frac{v_j(t+s_j,\theta)-C_n(\theta)}{\tilde v_j(s_j)-1},
 \qquad h_j(t)=Z_j(t,L).
\end{equation}
Evenness gives $v_j(s_j-t,\theta)=v_j(t-s_j,\theta)$.
Since $2s_j$ is either $0$ or $T_j$, periodicity gives
$v_j(t-s_j,\theta)=v_j(t+s_j,\theta)$. Hence $Z_j$ is
even in $t$. By~\eqref{eq:s3-he-stationary-extension-comparison},
$|Z_j|\leq C_n$, and the normalization gives $h_j(0)=1$.
In particular, $h_j'(0)=0$. The equations for $D_j$ yield
\begin{equation}\label{eq:s3-he-normalized-linear-equation}
 \begin{cases}
 \mathcal A Z_j=0
     &\text{in }\R\times\Sph^{n-1}_+,\\
 \partial_\beta Z_j(t,L)=q_j(t+s_j)h_j(t)
     &\text{for }t\in\R.
 \end{cases}
\end{equation}

Dividing~\eqref{eq:s3-he-stationary-schauder} by $a_j$ gives
uniform bounds for $Z_j$ in
$C^{k,\gamma}_{\mathrm{uloc}}(\R\times\Sph^{n-1}_+)$ for
every $k$ and $\gamma$. By \textup{(i)} and the formula for
$q_j$, we also have
\[
 q_j(\cdot+s_j)\longrightarrow\frac n{n-2}c
 \qquad\text{in }C^{k,\gamma}_{\mathrm{uloc}}(\R)
\]
for every $k$ and $\gamma$, since the norm is invariant
under translation. Every subsequence of $(Z_j)$ therefore
has a further subsequence converging in
$C^k([-R,R]\times\Sph^{n-1}_+)$ for every $k\geq0$ and $R>0$.
The limit $Z$
is bounded, even and cylindrically symmetric, and
$Z(0,L)=1$. Passing to the limit
in~\eqref{eq:s3-he-normalized-linear-equation}, we obtain
\begin{equation}\label{eq:s3-he-linearized-trace-limit}
 \begin{cases}
 \mathcal A Z=0
     &\text{in }\R\times\Sph^{n-1}_+,\\
 \partial_\beta Z(t,L)=\dfrac n{n-2}c\,Z(t,L)
     &\text{for }t\in\R.
 \end{cases}
\end{equation}

The Jacobi field $Z_{\xi_*}(t,\beta)=A_{\xi_*}(\beta)
\cos(\xi_*t)$ constructed above is bounded, even and
cylindrically symmetric, and $Z_{\xi_*}(0,L)=1$. The
difference $Z-Z_{\xi_*}$ is therefore a bounded even
Jacobi field along $C_n$ with trace zero at $t=0$.
Proposition~\ref{s3:rigidity}\textup{(iii)} gives
$Z=Z_{\xi_*}$. This identifies every subsequential limit
and yields
\begin{equation}\label{eq:s3-he-cosine-limit}
 h_j\longrightarrow h,\qquad h(t)=\cos(\xi_*t),
\end{equation}
in $C^k([-R,R])$ for every $k\geq0$ and $R>0$.

Since $s_j$ is a trace extremum, the first positive critical
point of $h_j$ is $T_j/2$ by
Theorem~\ref{s3:global-family}\textup{(ii)}. Fix
$0<\eta<T_*/4$ and choose $r\in(0,T_*/2-\eta)$ so small that
$h''<0$ on $[0,r]$. For large $j$, convergence in $C^2$
gives $h_j''<0$ on $[0,r]$. As $h_j'(0)=0$, integration
gives $h_j'<0$ on $(0,r]$. On $[r,T_*/2-\eta]$, the
derivative $h'$ is bounded above by a strictly negative
constant, so uniform convergence gives $h_j'<0$ there as
well. Finally, $h'(T_*/2+\eta)>0$, and therefore
$h_j'(T_*/2+\eta)>0$ for large $j$. The first positive zero
of $h_j'$ consequently satisfies
\[
 T_*/2-\eta<T_j/2<T_*/2+\eta.
\]
Since $\eta$ and the sequence $(\eps_j)$ converging to $1$
were arbitrary, this proves $\lim_{\eps\to1}T_\eps=T_*$.

Finally, since $\tilde v_\eps(t)=c^{-\frac{n-2}{2}}\widehat v_\eps(t,L)$
by~\eqref{eq:flat-normalization-table}, the gradient bound
in~\eqref{flatbase:gradient} gives
$|\tilde v_\eps'|/\tilde v_\eps\leq(n-2)/2$. Integrating on
$[0,T_\eps/2]$ and using~\eqref{eq:s3-family-bounds}, we obtain
\[
 \log\frac1\eps
 <\log\frac{\tilde v_\eps(T_\eps/2)}{\tilde v_\eps(0)}
 \leq\frac{n-2}{4}T_\eps.
\]
Thus $T_\eps\to+\infty$ as $\eps\to0$.
\end{proof}

A prescribed period $T$ may contain several repetitions of
a solution with fundamental period $T_\eps$.

\begin{corollary}\label{s3:prescribed-parameters}
Let $\eps\in(0,1)$ and $T>0$. There exists a $T$-periodic
solution in $\mathscr B_n$ with necksize $\eps$ if and only
if $T=mT_\eps$ for some positive integer $m$. In this case,
the solutions are precisely $v_\eps(t+\tau,\theta)$ with
$\tau\in\R$.
\end{corollary}

\begin{proof}
By Theorem~\ref{s3:global-family}\textup{(iii)}, the solutions
with necksize $\eps$ are $v_\eps(t+\tau,\theta)$ with
$\tau\in\R$. They have the same periods as $v_\eps$, so
Theorem~\ref{s3:global-family}\textup{(ii)} gives the
criterion for $T$.
\end{proof}

For small $\eps$, split the increasing half-period at the
first point where $\tilde v_\eps=\delta$, with $\delta>0$
fixed and small. The first interval gives the logarithmic
term, and the second has uniformly bounded length.

\begin{corollary}\label{s3:logarithmic-period}
Let $n\geq3$ and let $T_\eps$ be the fundamental period of
$v_\eps$ in Theorem~\ref{s3:global-family}. There exist
$\eps_0\in(0,1)$ and $C>0$, depending only on $n$, such that
\[
 \left|T_\eps-\frac{4}{n-2}\log\frac{1}{\eps}\right|
 \leq C
 \qquad\text{for every }\eps\in(0,\eps_0).
\]
\end{corollary}

\begin{proof}
By~\eqref{eq:s3-family-bounds} and the boundary condition
in~\eqref{eq:boundary-cylinder},
\begin{equation}\label{eq:s3-he-zero-endpoint-bounds}
 \begin{aligned}
 0<v_\eps(t,\theta)&\leq M_n C_n(\theta),\\
 0<\tilde v_\eps(t)&\leq M_n,\\
 0<\partial_\beta v_\eps(t,L)
   &=c\,\tilde v_\eps(t)^{\frac n{n-2}}
     \leq cM_n^{\frac n{n-2}}.
 \end{aligned}
\end{equation}
The Harnack inequality in the proof of
Proposition~\ref{flatbase:selection}, applied to
$c^{\frac{n-2}{2}}v_\eps$ on $[t-1,t+1]\times\Sph^{n-1}_+$,
compares $v_\eps(t,\beta)$ with $\tilde v_\eps(t)$, with a
constant independent of $t$, $\beta$ and $\eps$ because
$c\,\tilde v_\eps(t)^{2/(n-2)}\leq cM_n^{2/(n-2)}$.
Put $a=(n-2)/2$ and set
\begin{equation}\label{eq:s3-logarithmic-period-J}
 J(t)=\int_0^L\omega(\beta)v_\eps(t,\beta)\,\dd\beta.
\end{equation}
Integrating against $\omega$ gives
\[
 K^{-1}\tilde v_\eps(t)\leq J(t)\leq K\tilde v_\eps(t)
 \qquad\text{for every }t\in\R,
\]
where $K\geq1$ depends only on $n$.

Integrating the interior equation in~\eqref{eq:boundary-cylinder}
against $\omega$ gives
\[
 J''(t)=a^2J(t)-[\omega(\beta)\partial_\beta
                           v_\eps(t,\beta)]_0^L.
\]
The term at $\beta=0$ vanishes since $\omega(0)=0$. As
$\omega(L)=1$, the boundary condition in~\eqref{eq:boundary-cylinder}
gives
\[
 J''=(a^2-q)J,
 \qquad
 q(t)=\frac{c\,\tilde v_\eps(t)^{n/(n-2)}}{J(t)}.
\]
Since $v_\eps$ is even in $t$, we have $J'(0)=0$.

Fix $\delta\in(0,1)$ so small that
$cK\delta^{2/(n-2)}\leq a^2/2$. For $0<\eps<\delta$,
the trace increases strictly from $\eps$ at $t=0$ to a
value greater than $1$ at $t=T_\eps/2$. Hence there is a
unique $s_\eps\in(0,T_\eps/2)$ with
$\tilde v_\eps(s_\eps)=\delta$. On $[0,s_\eps]$, the
lower bound for $J$ gives
\[
 0\leq q(t)\leq cK\tilde v_\eps(t)^{2/(n-2)}
             \leq cK\delta^{2/(n-2)}\leq\frac{a^2}{2}.
\]
Thus $J''\geq a^2J/2$ on this interval, and $J'(0)=0$
gives $J'\geq0$. Multiplying $J''\leq a^2J$ by $2J'$
and integrating from $0$ to $t$ also gives
\[
 J'(t)^2\leq a^2\bigl(J(t)^2-J(0)^2\bigr).
\]
In particular, $0\leq J'\leq aJ$ on $[0,s_\eps]$.

Set $F=J'+aJ$. On $[0,s_\eps]$ we have
\[
 aJ\leq F\leq2aJ,\qquad
 F'=aF-qJ,\qquad
 (\log F)'=a-\frac{qJ}{F}.
\]
Since $J/F\leq1/a$ and $q\leq a^2/2$, we have
$(\log F)'\geq a/2$. Integrating from $t$ to $s_\eps$
and using $J(s_\eps)\leq K\delta$ yields
\[
 J(t)\leq \frac{F(t)}{a}
 \leq \frac{F(s_\eps)}{a}
       e^{-\frac a2(s_\eps-t)}
 \leq 2K\delta\,e^{-\frac a2(s_\eps-t)}
 \qquad (0\leq t\leq s_\eps).
\]
The inequality $\tilde v_\eps\leq KJ$ also gives
$q\leq C J^{2/(n-2)}$. Consequently,
\[
 \int_0^{s_\eps}q(t)\,\dd t
 \leq C\delta^{2/(n-2)}
       \int_0^\infty e^{-\frac{a}{n-2}r}\,\dd r
 \leq C,
\]
with a bound independent of $\eps$. Integrating the identity
for $(\log F)'$ now gives
\[
 0\leq as_\eps-\log\frac{F(s_\eps)}{F(0)}
   =\int_0^{s_\eps}\frac{q(t)J(t)}{F(t)}\,\dd t
   \leq C.
\]
Since $F(0)=aJ(0)$ and $aJ(s_\eps)\leq F(s_\eps)
\leq2aJ(s_\eps)$, the trace values at $0$ and $s_\eps$
give
\[
 \frac{\delta}{K^2\eps}
 \leq\frac{F(s_\eps)}{F(0)}
 \leq\frac{2K^2\delta}{\eps}.
\]
As $\delta$ is fixed in terms of $n$, it follows that
\[
 \left|s_\eps-\frac1a\log\frac1\eps\right|\leq C
 \qquad\text{for }0<\eps<\delta.
\]

To bound the length of $[s_\eps,T_\eps/2]$, let $\sigma_\eps$
be the unique point in $(s_\eps,T_\eps/2)$ such that
$\tilde v_\eps(\sigma_\eps)=1$. Let $I$ be as
in~\eqref{eq:s3-moment} for $v=v_\eps$. The bound
$v_\eps\leq M_n C_n$ gives $0<I\leq C$, and symmetry at
$0$ and $T_\eps/2$ gives $I'(0)=I'(T_\eps/2)=0$.
By~\eqref{flatbase:gradient} and~\eqref{eq:s3-he-zero-endpoint-bounds},
$|\tilde v_\eps'|\leq(n-2)M_n/2$. Choose
$\ell=(1-\delta)/((n-2)M_n)$. Then
$\sigma_\eps-s_\eps\geq2\ell$ and
\[
 \delta\leq\tilde v_\eps(t)\leq\frac{1+\delta}{2}
 \qquad\text{for }s_\eps\leq t\leq s_\eps+\ell.
\]
By~\eqref{eq:s3-moment-equation}, $I''\geq0$ on
$[0,\sigma_\eps]$ and $I''\geq d>0$ on
$[s_\eps,s_\eps+\ell]$, with $d$ depending only on $n$ and
$\delta$. Since $I'(0)=0$, this gives $I'\geq d\ell$ on
$[s_\eps+\ell,\sigma_\eps]$. Integration and boundedness of
$I$ give
\[
 \sigma_\eps-s_\eps\leq\ell+\frac{C}{d\ell}.
\]

On $[\sigma_\eps,T_\eps/2]$, the trace is increasing and
at least $1$. Since $x\mapsto x-x^{n/(n-2)}$ is decreasing
on $[1,\infty)$, $I'$ is concave there by~\eqref{eq:s3-moment-equation}.
Its graph lies above the chord
joining $I'(\sigma_\eps)\geq d\ell$ and $I'(T_\eps/2)=0$.
Integrating the chord gives
\[
 \frac{d\ell}{2}\bigl(T_\eps/2-\sigma_\eps\bigr)
 \leq I(T_\eps/2)-I(\sigma_\eps)\leq C.
\]
Thus $T_\eps/2-s_\eps$ is bounded independently of $\eps$.
Together with the estimate for $s_\eps$, this gives the
asserted inequality.
\end{proof}

The coefficient $4/(n-2)$ also occurs in the period estimate
for the classical Yamabe family~\cite[Proposition~2]{MP}.

\section{Global classification}\label{sec:global-classification}

We use the half-space Liouville theorem to classify the solutions
with a removable singularity at $0$. For solutions with a
nonremovable singularity, we use the classification of $\mathscr B_n$
in Theorem~\ref{s3:global-family}.

\subsection{Global solutions and the Kelvin reduction}
\label{subsec:global-kelvin-reduction}

We first classify the solutions of~\eqref{eq:intro-flat-equation}
whose traces $u=U(\cdot,0)$ extend continuously across $0$.
We then use this classification and Kelvin inversion to obtain
global two-sided bounds for $|x|^{(n-2)/2}U(x)$ when the singularity
at $0$ is nonremovable. These bounds and cylindrical symmetry
imply that the normalized Emden--Fowler transform belongs to
$\mathscr B_n$.

We normalize the coefficient in the boundary condition to $1$
by setting
\[
 \widehat U=\left(\frac{n-2}{2}\right)^{\frac{n-2}{2}}U,
 \qquad \widehat u=\widehat U(\cdot,0).
\]
The boundary condition becomes
$-\partial_{x_n}\widehat U(\bar x,0)=\widehat u(\bar x)^{n/(n-2)}$
on $\R^{n-1}\setminus\{0\}$. Thus $\widehat U$ satisfies the
extension problem studied in~\cite{CJXS,JLX} at order $1/2$
and in boundary dimension $n-1$.

We first show that $U$ is a bubble if
$|\bar x|^{(n-2)/2}\widehat u(\bar x)\to0$ as $\bar x\to0$.
By the annular Harnack inequality~\cite[Lemma~3.2, (21)]{CJXS},
$\widehat U(r\theta)=o(r^{-(n-2)/2})$ uniformly for
$\theta\in\Sph^{n-1}_+$.

We use the barrier argument in the proof of~\cite[Proposition~3.4]{CJXS}.
With $x=r\theta$ as in Section~\ref{subsec:notation}, set
\[
 w_\varepsilon(r\theta)
 =\left(Cr^{-\frac{n-2}{4}}
          +\varepsilon r^{-\frac{n-2}{2}}\right)
   (1-\delta\theta_n),\qquad \varepsilon>0.
\]
At $\delta=0$, both radial powers have negative Laplacian
because their exponents lie in $(2-n,0)$. By continuity on
$\Sph^{n-1}_+$ and homogeneity of each term, $w_\varepsilon$
is positive and $\Delta w_\varepsilon<0$ for sufficiently small
$0<\delta<1$. On $\partial\R^n_+\setminus\{0\}$,
\[
 -\partial_{x_n}w_\varepsilon(\bar x,0)
 =\frac{\delta}{|\bar x|}w_\varepsilon(\bar x,0).
\]
Choose $\rho\in(0,1)$ so that
$\widehat u(\bar x)^{2/(n-2)}<\delta/|\bar x|$ for
$0<|\bar x|<\rho$. Take $C$ large enough that
$w_\varepsilon\geq\widehat U$ on $\{x\in\R^n_+:|x|=\rho\}$
for every $\varepsilon>0$. For fixed $\varepsilon$, the
Harnack estimate gives $\widehat U/w_\varepsilon\to0$ as
$x\to0$.

If $\widehat U/w_\varepsilon$ exceeds $1$ for $0<|x|<\rho$,
its maximum is attained either at an interior point or on
the flat boundary. At an interior maximum, the product rule
gives $0=\Delta\widehat U\leq
(\widehat U/w_\varepsilon)\Delta w_\varepsilon<0$.
At a maximum on the flat boundary, $-\partial_{x_n}$
of the quotient is nonnegative, whereas
\[
 -\partial_{x_n}\left(\frac{\widehat U}{w_\varepsilon}\right)
 =\left(\widehat u^{2/(n-2)}-\frac{\delta}{|\bar x|}\right)
   \frac{\widehat U}{w_\varepsilon}<0.
\]
Thus $\widehat U\leq w_\varepsilon$, and letting
$\varepsilon\to0^+$ gives
\[
 \widehat U(r\theta)\leq Cr^{-\frac{n-2}{4}}
 \qquad\text{for }0<r<\rho,\quad\theta\in\Sph^{n-1}_+.
\]

For $0<r<\rho/2$, the functions
$r^{(n-2)/4}\widehat U(r\,\cdot)$ are uniformly bounded on
$\{x\in\R^n_+:1/2<|x|<2\}$. Their boundary condition has
coefficient $r^{1/2}$ in front of the nonlinearity.
The local H\"older estimate used in the proof of
Proposition~\ref{prop:flat-uniform-estimates} therefore
gives uniform H\"older bounds in a neighborhood of $\Sph^{n-1}_+$.
Since $s\mapsto s^{n/(n-2)}$ is Lipschitz on bounded intervals
in $[0,\infty)$, the boundary condition gives uniform H\"older
bounds with the same exponent for the normal derivatives of
the rescaled functions.
The first Neumann Schauder estimate in that proof, together
with the interior gradient estimate, gives a uniform gradient
bound on $\Sph^{n-1}_+$. Returning to the original
variables gives
\[
 |\nabla\widehat U(r\theta)|\leq Cr^{-1-\frac{n-2}{4}}.
\]
The two bounds imply
$\widehat U\in H^1(\{x\in\R^n:x_n>0,\ |x|<\rho/2\})$
since $2+(n-2)/2<n$.

For a smooth test function supported in $|x|<\rho/2$, integrate
by parts on $\{x\in\R^n_+:x_n>0,\ r<|x|<\rho/2\}$.
The integral on $r\Sph^{n-1}_+$ is
$O(r^{3(n-2)/4})$ by the gradient bound. The flat boundary
term converges because $\widehat u^{n/(n-2)}
=O(|\bar x|^{-n/4})$ is locally integrable. Letting $r\to0$ gives the weak equation across $0$. Together
with smoothness away from $0$, this shows that $\widehat U$
is a nonzero nonnegative weak solution on the whole half-space
and belongs to $H^1_{\mathrm{loc}}$ up to the boundary.
We can therefore apply the half-space Liouville theorem
of Jin, Li and Xiong~\cite[Theorem~1.8]{JLX} with boundary
dimension $n-1$ and fractional order $1/2$. Returning to
the original normalization, we obtain
\[
 u(\bar x)=A\bigl(|\bar x-a|^2+\lambda^2\bigr)^{-\frac{n-2}{2}}
 \qquad
 \bigl(A>0,\ a\in\R^{n-1},\ \lambda>0\bigr),
\]
and $U$ is the Poisson extension of $u$. For $z=(a,-\lambda)$,
the function $A|x-z|^{2-n}$ is bounded and harmonic on the
half-space and has trace $u$. It therefore equals $U$.
In particular, the singularity is removable. Conversely, if
$u$ extends continuously across $0$, then
$|\bar x|^{(n-2)/2}\widehat u(\bar x)\to0$, and the same
argument gives this formula for $U$.

For $U(x)=A|x-z|^{2-n}$, the boundary condition
in~\eqref{eq:intro-flat-equation} gives
\[
 (n-2)\lambda A|(\bar x,0)-z|^{-n}
 =\frac{n-2}{2}A^{\frac n{n-2}}|(\bar x,0)-z|^{-n}.
\]
It follows that $A^{2/(n-2)}=2\lambda$, so $U$ has the
form~\eqref{eq:intro-flat-bubble}.

\begin{proposition}\label{s3:global-bilateral-reduction}
Let $n\geq3$ and let $U$ be a positive solution of
\eqref{eq:intro-flat-equation}, smooth on
$\R^n_+\setminus\{0\}$. Suppose that $u=U(\cdot,0)$ has
a nonremovable singularity at $0$. Then
there exist constants $0<m\leq M<\infty$ such that
\begin{equation}\label{eq:s3-global-bilateral-bounds}
 m|x|^{-\frac{n-2}{2}}\leq U(x)\leq M|x|^{-\frac{n-2}{2}}
 \qquad\text{for every }x\in \R^n_+\setminus\{0\}.
\end{equation}
In cylindrical variables, the normalized Emden--Fowler
transform of $U$, given by
\begin{equation}\label{eq:s3-global-cylinder-reduction}
 v(t,\theta)=b_n^{-1}e^{-\frac{n-2}{2}t}U(e^{-t}\theta),
 \qquad (t,\theta)\in\R\times\Sph^{n-1}_+,
\end{equation}
belongs to $\mathscr B_n$.
\end{proposition}

\begin{proof}
By~\cite[Proposition~3.1]{CJXS},
$|\bar x|^{(n-2)/2}\widehat u(\bar x)$ is bounded above near $0$.
If its lower limit at $0$ were zero,
\cite[Proposition~3.3]{CJXS} would give convergence to zero,
and the preceding removability argument would give a continuous
extension of $u$. Thus it is bounded away from zero near $0$.
For small $r>0$, the Harnack inequality
\cite[Lemma~3.2, (21)]{CJXS} applies on
\[
 \{x\in\R^n_+:x_n>0,\ r/2<|x|<2r\},
\]
with a constant independent of $r$. By continuity up to
$\partial\R^n_+$, it bounds $\widehat U(r\theta)$ above and
below by positive multiples of $\widehat u(r\vartheta)$,
uniformly for $\theta\in\Sph^{n-1}_+$ and
$\vartheta\in\Sph^{n-2}$. Since $\widehat U$ is a fixed positive
multiple of $U$, the bounds for
$|\bar x|^{(n-2)/2}\widehat u(\bar x)$ give
\eqref{eq:s3-global-bilateral-bounds} near $0$.

Take $a=0$ and $\rho=1$ in~\eqref{eq:flat-kelvin}. The function
\[
 U_{0,1}(x)=|x|^{2-n}U\!\left(\frac{x}{|x|^2}\right)
\]
satisfies~\eqref{eq:intro-flat-equation}. If $U_{0,1}(\cdot,0)$
extended continuously across $0$, then $U_{0,1}$ would be a
bubble of the form~\eqref{eq:intro-flat-bubble}. For a bubble
with parameters $(a,\lambda)$, put $z=(a,-\lambda)$. The unit
Kelvin transform gives the bubble with parameters
$(a/|z|^2,\lambda/|z|^2)$. Since $(U_{0,1})_{0,1}=U$, it would
follow that $U$ is a bubble, contrary to the nonremovability of
the singularity of $u$ at $0$.

The estimates near $0$ therefore apply to both $U$ and
$U_{0,1}$. By Section~\ref{subsec:notation}, $U_{0,1}$ corresponds
to $t\mapsto-t$. In terms of $r=e^{-t}$, this gives
\[
 r^{\frac{n-2}{2}}U(r\theta)
 =(r^{-1})^{\frac{n-2}{2}}U_{0,1}(r^{-1}\theta).
\]
Choose $r_0\in(0,1)$ so that the local estimates hold for
both functions when $0<|x|<r_0$. The identity gives the
bounds for $|x|>r_0^{-1}$. On the compact set
$\{x\in\R^n_+:r_0\leq |x|\leq r_0^{-1}\}$, the function
$|x|^{(n-2)/2}U(x)$ is continuous and positive. Its minimum
and maximum on this set complete the bounds in
\eqref{eq:s3-global-bilateral-bounds}.

By~\cite[Theorem~1.3]{CJXS}, $\widehat U$ is invariant under
rotations of $\bar x$, and so is $U$. Thus the function
$v$ in~\eqref{eq:s3-global-cylinder-reduction} is
cylindrically symmetric and satisfies
\eqref{eq:boundary-cylinder} by Section~\ref{subsec:notation}.
Moreover, \eqref{eq:s3-global-bilateral-bounds} gives
\[
 m b_n^{-1}\leq\tilde v(t)\leq M b_n^{-1}
 \qquad\text{for every }t\in\R.
\]
Therefore $v\in\mathscr B_n$.
\end{proof}

\subsection{Proof of Theorem B}
\label{subsec:global-classification-proof}

We recover the parameters of a bubble from $u=U(\cdot,0)$.
When the singularity at $0$ is nonremovable, we use
Theorem~\ref{s3:global-family} to identify the necksize and phase.
As in that theorem, we set $v_1=C_n$.

\begin{restatedtheorem}\label{s3:global-flat-classification}
Let $n\geq3$ and let $U$ be a positive solution of
\eqref{eq:intro-flat-equation}, smooth on
$\R^n_+\setminus\{0\}$. Exactly one of the following holds.

\begin{enumerate}[label=\textup{(\roman*)}]
\item The singularity at $0$ is removable, and $U$ is the
bubble
\begin{equation}\label{eq:s3-global-regular-classification}
 U(\bar x,x_n)=
 \left(\frac{2\lambda}
 {|\bar x-a|^2+(x_n+\lambda)^2}\right)^{\frac{n-2}{2}}
\end{equation}
for unique parameters $a\in\R^{n-1}$ and $\lambda>0$.

\item The singularity at $0$ is nonremovable. There exist
a unique necksize $\eps\in(0,1]$ and a phase $\tau\in\R$
such that, writing $x=r\theta$ with $r>0$ and
$\theta\in\Sph^{n-1}_+$, we have
\begin{equation}\label{eq:s3-global-singular-classification}
 U(r\theta)=b_nr^{-(n-2)/2}\,
             v_\eps\bigl(-\log r+\tau,\theta\bigr).
\end{equation}
For $0<\eps<1$, the phase is unique modulo the fundamental
period $T_\eps$. For $\eps=1$, all choices of $\tau$ give
the same cylindrical solution.
\end{enumerate}

Conversely, both formulas give positive solutions of
\eqref{eq:intro-flat-equation}. The bubbles in \textup{(i)}
extend smoothly across $0$, and the solutions in \textup{(ii)}
have a nonremovable singularity there.
\end{restatedtheorem}

\begin{proof}
If the singularity at $0$ is removable,
Section~\ref{subsec:global-kelvin-reduction} gives the bubble
in~\eqref{eq:s3-global-regular-classification}. Its trace is
\[
 u(\bar x)=U(\bar x,0)
 =\left(\frac{2\lambda}{|\bar x-a|^2+\lambda^2}
   \right)^{\frac{n-2}{2}}.
\]
The unique maximum of $u$ is attained at $\bar x=a$ and is
equal to $(2/\lambda)^{(n-2)/2}$. This determines $a$ and $\lambda$.

If the singularity at $0$ is nonremovable,
Proposition~\ref{s3:global-bilateral-reduction} shows that
the normalized transform $v$ defined in
\eqref{eq:s3-global-cylinder-reduction} belongs to $\mathscr B_n$.
By Theorem~\ref{s3:global-family}\textup{(iii)}, either
$v=C_n$ or $v(t,\theta)=v_\eps(t+\tau,\theta)$ for some
$0<\eps<1$ and $\tau\in\R$. The inverse transformation
gives~\eqref{eq:s3-global-singular-classification}, with
$\eps=1$ in the cylindrical case. Since $v_1=C_n$ is
independent of $t$, every choice of $\tau$ gives the same
solution in this case.

For any $\vartheta\in\Sph^{n-2}$, the normalized trace is
\[
 \tilde v(t)
 =b_n^{-1}e^{-\frac{n-2}{2}t}u(e^{-t}\vartheta)
 =\tilde v_\eps(t+\tau).
\]
Thus $\eps=\min_{t\in\R}\tilde v(t)$, including the
cylindrical case, where $\tilde v_1\equiv1$. Since the Emden--Fowler transformation is invertible,
Theorem~\ref{s3:global-family}\textup{(iii)} gives uniqueness
of $\tau$ modulo $T_\eps$ when $0<\eps<1$.

Conversely, the bubbles
in~\eqref{eq:s3-global-regular-classification} solve
\eqref{eq:intro-flat-equation} by the computation in
Section~\ref{subsec:global-kelvin-reduction} and are smooth
across $0$. By Section~\ref{subsec:notation}, the functions in
\eqref{eq:s3-global-singular-classification} solve
\eqref{eq:intro-flat-equation} because $v_\eps$ satisfies
\eqref{eq:boundary-cylinder}. For $r>0$ and
$\vartheta\in\Sph^{n-2}$, their traces are
\[
 u(r\vartheta)=b_nr^{-(n-2)/2}\,
       \tilde v_\eps(-\log r+\tau)
 \geq b_n\eps r^{-(n-2)/2}.
\]
The lower bound tends to $+\infty$ as $r\to0$, so the
singularity at $0$ is nonremovable.
\end{proof}

\appendix

\section{Estimates for the extension and its boundary trace}
\label{sec:flat-foundations}

We establish the estimates used in Proposition~\ref{s3:rigidity}.
For $n\geq3$, we bound $\widehat v$ on $\R\times\Sph^{n-1}_+$
in terms of two-sided bounds for $\widehat v(t,L)$ and estimate
$\partial_t\widehat v(t,L)$. For $n\geq5$, we also prove
second-order differential inequalities for $\widehat v(t,L)$
and use them to estimate $b(t)$ in~\eqref{flatbase:secant-coefficient},
the coefficient in the boundary condition for differences of
solutions and for Jacobi fields.

We then consider the restriction of $\mathcal A$ from
Section~\ref{subsec:linearized-spaces} to functions $f=f(\beta)$
with boundary condition $f'(L)=bf(L)$, where $b\geq0$.
For $n\geq5$, we estimate the first eigenvalue and its positive
eigenfunction and bound the gap between the first two eigenvalues
from below. We conclude with a resolvent identity valid for
all $n\geq3$.

\subsection{Bounds for the extension}
\label{subsec:flat-extension-bounds}

We use bounds for $\widehat v(t,L)$ to obtain a comparison
with $C_n$ on $\R\times\Sph^{n-1}_+$. We also prove uniqueness
for the bounded extension with trace $\widehat v(t,L)$ and
obtain uniform regularity for its translates in $t$.

By cylindrical symmetry, the normalization $\widehat v$
in~\eqref{eq:flat-normalization-table} depends only on $t$ and
$\beta$. With $\beta$, $L$ and $\omega$ as in
Section~\ref{subsec:notation}, \eqref{eq:s3-rp-auxiliary-pair}
takes the form
\begin{equation}\label{flatbase:aux-equation}
 \begin{cases}
 \widehat v_{tt}
   +\omega^{-1}(\omega\widehat v_\beta)_\beta
   -\dfrac{(n-2)^2}{4}\widehat v=0
       & \text{for } (t,\beta)\in\R\times(0,L),\\
 \widehat v_\beta(t,L)=\widehat v(t,L)^{\frac n{n-2}}
       & \text{for }t\in\R.
 \end{cases}
\end{equation}

\begin{proposition}\label{flatbase:selection}
Let $n\geq3$ and let $\widehat v$ be a positive smooth
cylindrically symmetric solution of~\eqref{flatbase:aux-equation}
on $\R\times\Sph^{n-1}_+$. Suppose that there are constants
$m,M>0$ such that
\[
 m\leq\widehat v(t,L)\leq M\qquad\text{for all }t\in\R.
\]
Then
\[
 mC_n(\theta)\leq\widehat v(t,\theta)\leq MC_n(\theta)
 \qquad\text{on }\R\times\Sph^{n-1}_+.
\]
Here $C_n$ is given by~\eqref{eq:s3-stationary-angular}.
The function $\widehat v$ is the unique bounded solution of
$\mathcal A w=0$ on $\R\times\Sph^{n-1}_+$ with
$w(t,(\vartheta,0))=\widehat v(t,L)$ for every
$(t,\vartheta)\in\R\times\Sph^{n-2}$, where $\mathcal A$ is
as in Section~\ref{subsec:linearized-spaces}. For every
integer $k\geq0$ and $0<\gamma<1$, $\widehat v$ belongs to
$C^{k,\gamma}_{\mathrm{uloc}}(\R\times\Sph^{n-1}_+)$, and the
family $\{\widehat v(t+\tau,\theta):\tau\in\R\}$ is precompact
in $C^{k,\gamma}_{\mathrm{loc}}(\R\times\Sph^{n-1}_+)$.
\end{proposition}

\begin{proof}
Write the boundary condition in~\eqref{flatbase:aux-equation}
as $\widehat v_\beta(t,L)=b(t)\widehat v(t,L)$, where
\[
 b(t)=\widehat v(t,L)^{\frac{2}{n-2}}
       \leq M^{\frac{2}{n-2}}.
\]
Set $\rho=L/4$. In the coordinates $(t,\beta,\vartheta)$ of
Section~\ref{subsec:notation}, the set of points in
$\R\times\Sph^{n-1}_+$ at distance less than $\rho$ from
$\R\times\partial\Sph^{n-1}_+$ is
\[
 \{(t,\beta,\vartheta):t\in\R,\ L-\rho<\beta\leq L,\ 
                              \vartheta\in\Sph^{n-2}\},
\]
since the distance to $\partial\Sph^{n-1}_+$ is $L-\beta$.
Choose a fixed smooth function $\chi(\beta)$
that is zero for $\beta\leq L-\rho$ and equal to $1$ for
$\beta\geq L-\rho/2$. Set $X_b=b(t)\chi(\beta)\partial_\beta$,
with $X_b=0$ where $\beta\leq L-\rho$. Both $X_b$ and
\[
 \operatorname{div}X_b
 =b(t)\bigl(\chi'(\beta)+(n-2)\cot\beta\,\chi(\beta)\bigr)
\]
are bounded in terms of $n$ and $M$. The divergence
only involves differentiation in $\beta$, since $X_b$ has no
component in the $t$ direction. At $\beta=L$,
$\langle X_b,\partial_\beta\rangle=b(t)$, so
$\langle\nabla\widehat v-\widehat vX_b,\partial_\beta\rangle=0$
by~\eqref{flatbase:aux-equation}. Applying the divergence theorem
to $\widehat v\phi X_b$ in the weak formulation of
\eqref{flatbase:aux-equation} gives
\[
 \int_{\R\times\Sph^{n-1}_+}
 \left[
 \begin{aligned}
 &\langle\nabla\widehat v-\widehat v X_b,\nabla\phi\rangle\\
 &-\phi\langle X_b,\nabla\widehat v\rangle
 +\left(\frac{(n-2)^2}{4}-\operatorname{div}X_b\right)\widehat v\phi
 \end{aligned}
 \right]\dd\mu=0,
\]
for every $\phi\in C_c^\infty(\R\times\Sph^{n-1}_+)$,
with support allowed to meet $\R\times\partial\Sph^{n-1}_+$.
Here the gradient and divergence are taken with respect to
the product metric of Section~\ref{subsec:notation}, with
volume element $\dd\mu$.

Let $\iota(t,\theta)=(t,(\theta_1,\ldots,\theta_{n-1},-\theta_n))$.
In the coordinate $s=\beta-L$ on $|s|<\rho$, the product
metric is $\dd t^2+\dd s^2+\cos^2s\,g_{\Sph^{n-2}}$ and
$\iota(t,s,\vartheta)=(t,-s,\vartheta)$. Extend $\widehat v$
and the bounded scalar $\operatorname{div}X_b$ to
$\iota(\R\times\Sph^{n-1}_+)$ by composition with $\iota$,
and extend $X_b$ to this set by the differential of $\iota$.
Write $\widehat v^{\,e}$, $(\operatorname{div}X_b)^e$ and
$X_b^{\,e}$ for the extensions. For $0<s<\rho$,
\[
 X_b^{\,e}(t,s,\vartheta)
   =-b(t)\chi(L-s)\partial_s.
\]
The two pieces defining $\widehat v^{\,e}$ have the same
trace on $\R\times\partial\Sph^{n-1}_+$, so
$\widehat v^{\,e}\in H^1_{\mathrm{loc}}(\R\times\Sph^{n-1})$.

Let $\psi\in C_c^\infty(\R\times\Sph^{n-1})$ and set
\[
 \phi=(\psi+\psi\circ\iota)|_{\R\times\Sph^{n-1}_+}.
\]
Split the integral over $\R\times\Sph^{n-1}$ into the integrals
over $\R\times\Sph^{n-1}_+$ and
$\iota(\R\times\Sph^{n-1}_+)$. On the second set, change
variables by $\iota$. Since $\iota$ preserves the metric and
volume element, we obtain
\[
 \begin{aligned}
 &\int_{\R\times\Sph^{n-1}}
 \left[
 \begin{aligned}
 &\langle\nabla\widehat v^{\,e}-\widehat v^{\,e}X_b^{\,e},
         \nabla\psi\rangle\\
 &-\psi\langle X_b^{\,e},\nabla\widehat v^{\,e}\rangle
 +\left(\frac{(n-2)^2}{4}-(\operatorname{div}X_b)^e\right)
                  \widehat v^{\,e}\psi
 \end{aligned}
 \right]\dd\mu\\
 &\quad=
 \int_{\R\times\Sph^{n-1}_+}
 \left[
 \begin{aligned}
 &\langle\nabla\widehat v-\widehat v X_b,\nabla\phi\rangle\\
 &-\phi\langle X_b,\nabla\widehat v\rangle
 +\left(\frac{(n-2)^2}{4}-\operatorname{div}X_b\right)\widehat v\phi
 \end{aligned}
 \right]\dd\mu
 =0.
 \end{aligned}
\]
Thus $\widehat v^{\,e}$ satisfies a homogeneous weak equation
on $\R\times\Sph^{n-1}$. In local coordinates, the product
metric gives a uniformly elliptic principal part in divergence
form and the lower-order coefficients are bounded in terms of
$n$ and $M$. Applying the interior Harnack inequality
\cite[Theorem~7]{Serrin} on a fixed chain of overlapping charts
in $(-2,2)\times\Sph^{n-1}$, translated by $\tau$ in $t$, gives
\[
 \begin{aligned}
 \sup_{[\tau-1,\tau+1]\times\Sph^{n-1}_+}\widehat v
 &\leq C\inf_{[\tau-1,\tau+1]\times\Sph^{n-1}_+}\widehat v\\
 &\leq C\widehat v(\tau,L)\leq CM
 \qquad(\tau\in\R).
 \end{aligned}
\]
Here $C=C(n,M)$ is independent of $\tau$ and $m$.

We now use the maximum principle for the interior operator
$\mathcal A$ of Section~\ref{subsec:linearized-spaces}. Let $w$
be a bounded solution of $\mathcal A w=0$ on
$\R\times\Sph^{n-1}_+$ with $w|_{\R\times\partial\Sph^{n-1}_+}
\leq0$. Fix $0<\delta<(n-2)/2$ and $\varepsilon>0$, and
set $w_\varepsilon=w-\varepsilon\cosh(\delta t)$. Then
\[
 \mathcal A w_\varepsilon
 =-\varepsilon\left(\frac{(n-2)^2}{4}-\delta^2\right)
                     \cosh(\delta t)<0.
\]
We have $w_\varepsilon\leq0$ on
$[-R,R]\times\partial\Sph^{n-1}_+$ and on
$\{-R,R\}\times\Sph^{n-1}_+$ whenever
$\varepsilon\cosh(\delta R)\geq
\|w\|_{L^\infty(\R\times\Sph^{n-1}_+)}$.
The maximum principle on $[-R,R]\times\Sph^{n-1}_+$ gives
$w_\varepsilon\leq0$ there for every sufficiently large $R$.
Hence $w_\varepsilon\leq0$ on the whole cylinder. Letting
$\varepsilon\to0^+$ gives $w\leq0$.

By~\eqref{eq:s3-stationary-angular}, $\mathcal A C_n=0$ and
$C_n(L)=1$. The Harnack bound gives boundedness of
$\widehat v-MC_n$ and $mC_n-\widehat v$, whose values at
$\beta=L$ are nonpositive. The comparison gives
$mC_n\leq\widehat v\leq MC_n$. If $w_1,w_2$ are bounded
solutions with $w_1|_{\R\times\partial\Sph^{n-1}_+}
=w_2|_{\R\times\partial\Sph^{n-1}_+}$, the same argument
applied to $w_1-w_2$ and $w_2-w_1$ gives $w_1=w_2$.

Since $C_n>0$ on $\Sph^{n-1}_+$, which is compact, the
comparison gives
\[
 0<m\min_{\Sph^{n-1}_+} C_n
 \leq \widehat v(t,\theta)
 \leq M\max_{\Sph^{n-1}_+} C_n
 \qquad\text{for all }(t,\theta)\in\R\times\Sph^{n-1}_+.
\]
For $\tau\in\R$, set $\widehat v_\tau(t,\theta)
=\widehat v(t+\tau,\theta)$. By translation invariance,
$\widehat v_\tau$ solves~\eqref{flatbase:aux-equation} and
satisfies the same bounds. Proposition~\ref{prop:flat-uniform-estimates}
for the boundary coefficient $1$ gives, for every integer $k\geq0$
and $0<\gamma<1$,
\[
 \sup_{\tau\in\R}
 \|\widehat v_\tau\|_{C^{k,\gamma}_{\mathrm{uloc}}
                       (\R\times\Sph^{n-1}_+)}
 \leq C(n,m,M,k,\gamma).
\]
By the compactness assertion of the same proposition, every
sequence $(\widehat v_{\tau_j})_{j\geq1}$ has a subsequence
converging in $C^{k,\gamma}_{\mathrm{loc}}(\R\times\Sph^{n-1}_+)$.
\end{proof}

\subsection{Bounds for the trace}
\label{subsec:flat-trace-bounds}

We prove a gradient bound on $\R\times(0,L]$ which controls
$\partial_t\widehat v(t,L)$ for $n\geq3$. For $n\geq5$, we
use the gradient bound to derive differential inequalities for
$\widehat v(t,L)$ and then estimate the boundary coefficient
in~\eqref{flatbase:secant-coefficient}. After changing the
normalization, the boundary estimate below and the first-order
bound for the secant coefficient also appear in
\cite[Proposition~2.2 and Lemma~2.3]{JinLiaoXiong}.

\begin{proposition}\label{flatbase:trace-gradient}
Let $n\geq3$ and let $\widehat v$ satisfy the hypotheses of
Proposition~\ref{flatbase:selection}. On $\R\times(0,L]$, consider
the metric
\[
 g_*=\sin^{-2}\beta\,(\dd t^2+\dd\beta^2).
\]
Then
\[
 \left|\nabla_{g_*}
  \log\bigl(\sin^{\frac{n-2}{2}}\beta\,\widehat v\bigr)
 \right|_{g_*}^2
 \leq\frac{(n-2)^2}{4}.
\]
In particular,
\begin{equation}\label{flatbase:gradient}
 \left(\frac{\partial_t\widehat v(t,L)}{\widehat v(t,L)}\right)^2
 +\widehat v(t,L)^{\frac4{n-2}}
 \leq\frac{(n-2)^2}{4}
 \qquad\text{for all }t\in\R.
\end{equation}
\end{proposition}

\begin{proof}
The interior equation in~\eqref{flatbase:aux-equation} gives
\[
 (\partial_t^2+\partial_\beta^2)
   \bigl(\sin^{\frac{n-2}{2}}\beta\,\widehat v\bigr)
 =\frac{(n-2)(n-4)}4\csc^2\beta\,
       \sin^{\frac{n-2}{2}}\beta\,\widehat v.
\]
Since $\Delta_{g_*}=\sin^2\beta(\partial_t^2+\partial_\beta^2)$,
we obtain
\begin{equation}\label{eq:flat-trace-log-laplacian}
 \Delta_{g_*}
  \log\bigl(\sin^{\frac{n-2}{2}}\beta\,\widehat v\bigr)
 =\frac{(n-2)(n-4)}4-
 \left|\nabla_{g_*}
  \log\bigl(\sin^{\frac{n-2}{2}}\beta\,\widehat v\bigr)
 \right|_{g_*}^2.
\end{equation}
Define
\[
 \Psi(t,\beta)=
 \frac{
 \left|\nabla_{g_*}
  \log\bigl(\sin^{\frac{n-2}{2}}\beta\,\widehat v\bigr)
 \right|_{g_*}^2-\frac{(n-2)^2}{4}}
 {\sin\beta\,\widehat v(t,\beta)^{\frac2{n-2}}}.
\]

In the coordinates of Section~\ref{subsec:notation}, cylindrical
symmetry and smoothness give an even smooth extension of
$\widehat v(t,\beta)$ across $\beta=0$. In particular,
$\partial_\beta\widehat v(t,0)=0$.
Using the even extension and the uniform estimates in
Proposition~\ref{flatbase:selection} in the definition of $\Psi$
gives
\[
 \Psi(t,\beta)=O(\beta)
 \qquad\text{as }\beta\to0,
\]
uniformly in $t$. The same estimates bound $\Psi$ on
$\R\times[\delta,L]$ for each $\delta>0$, so
$\sup_{\R\times(0,L]}\Psi$ is finite.

Suppose that $\sup_{\R\times(0,L]}\Psi>0$ and choose
$(t_j,\beta_j)$ such that
\[
 \Psi(t_j,\beta_j)\longrightarrow\sup_{\R\times(0,L]}\Psi.
\]
Since $\Psi(t,\beta)=O(\beta)$ uniformly in $t$, there is
$\delta_0>0$ such that $\beta_j\geq\delta_0$ for all sufficiently
large $j$. After passing to a subsequence,
$\beta_j\to\beta_\infty\in[\delta_0,L]$. By
Proposition~\ref{flatbase:selection}, a further subsequence of
$\widehat v(t+t_j,\theta)$ converges in
$C^k([-R,R]\times\Sph^{n-1}_+)$ for every $R>0$ and every
integer $k\geq0$. We denote the limit again by $\widehat v$.
For this limit, $\Psi$ attains a positive maximum at
$(0,\beta_\infty)$.

Suppose first that the maximum occurs with $0<\beta<L$.
Since $\nabla_{g_*}\Psi=0$ there, differentiation of the
definition of $\Psi$ gives
\begin{multline}\label{eq:flat-trace-maximum-gradient}
 \nabla_{g_*}
 \left|\nabla_{g_*}
  \log\bigl(\sin^{\frac{n-2}{2}}\beta\,\widehat v\bigr)
 \right|_{g_*}^2\\
 =\frac{2\Psi\sin\beta\,\widehat v^{\frac2{n-2}}}{n-2}
 \nabla_{g_*}
  \log\bigl(\sin^{\frac{n-2}{2}}\beta\,\widehat v\bigr).
\end{multline}
For every smooth function $f$ and tangent vector $Y$,
\[
 2\nabla_{g_*}^2 f(\nabla_{g_*}f,Y)
 =Y\bigl(|\nabla_{g_*}f|_{g_*}^2\bigr).
\]
The gradient on the right of~\eqref{eq:flat-trace-maximum-gradient}
is nonzero since $\Psi>0$. By~\eqref{eq:flat-trace-maximum-gradient}
and the Hessian identity above, it is an eigenvector of
$\nabla_{g_*}^2\log(\sin^{\frac{n-2}{2}}\beta\,\widehat v)$,
viewed as an endomorphism with respect to $g_*$. Its eigenvalue
and the other eigenvalue are, respectively,
\[
 \frac{\Psi\sin\beta\,\widehat v^{\frac2{n-2}}}{n-2},
 \qquad
 -\frac{n-2}{2}
 -\frac{n-1}{n-2}\Psi\sin\beta\,\widehat v^{\frac2{n-2}}.
\]
Their sum is the Laplacian in~\eqref{eq:flat-trace-log-laplacian}.

The curvature of $g_*$ is
$-\sin^2\beta\,\partial_\beta^2(-\log\sin\beta)=-1$.
Using the two Hessian eigenvalues in Bochner's formula and
applying the product rule to
$\Psi\sin\beta\,\widehat v^{2/(n-2)}$ at $\nabla_{g_*}\Psi=0$,
we obtain
\[
 \Delta_{g_*}\Psi
 =\Psi\left(n-2+\frac{2(n-1)}{n-2}
              \Psi\sin\beta\,\widehat v^{\frac2{n-2}}\right)>0.
\]
Here we used~\eqref{eq:flat-trace-log-laplacian}
and~\eqref{eq:flat-trace-maximum-gradient}.
This contradicts $\Delta_{g_*}\Psi\leq0$ at the maximum.

It remains to consider $\beta_\infty=L$. The boundary condition
in~\eqref{flatbase:aux-equation} gives
\[
 \left.\partial_\beta
  \log\bigl(\sin^{\frac{n-2}{2}}\beta\,\widehat v\bigr)
 \right|_{\beta=L}
 =\widehat v(t,L)^{\frac2{n-2}}.
\]
Consequently,
\[
 \Psi(t,L)
 =\widehat v(t,L)^{-\frac2{n-2}}
 \left[
   \left(\partial_t\log\widehat v(t,L)\right)^2
   +\widehat v(t,L)^{\frac4{n-2}}
   -\frac{(n-2)^2}{4}
 \right].
\]
Differentiating this identity gives
\[
 \begin{aligned}
 \partial_t\Psi(t,L)
 ={}&2\widehat v(t,L)^{-\frac2{n-2}}
       \partial_t\log\widehat v(t,L)\\
 &{}\times
 \left[
  \partial_t^2\log\widehat v(t,L)
  +\frac2{n-2}\widehat v(t,L)^{\frac4{n-2}}
  -\frac{\Psi(t,L)}{n-2}\widehat v(t,L)^{\frac2{n-2}}
 \right].
 \end{aligned}
\]
At the maximum, $\partial_t\Psi=0$ and $\partial_t^2\Psi\leq0$.
If $\partial_t\log\widehat v(t,L)\ne0$, the bracket vanishes.
If $\partial_t\log\widehat v(t,L)=0$, differentiating once
more gives
\[
 \partial_t^2\log\widehat v(t,L)
 \left[
  \partial_t^2\log\widehat v(t,L)
  +\frac{\widehat v(t,L)^{\frac4{n-2}}
                    +\frac{(n-2)^2}{4}}{n-2}
 \right]\leq0.
\]
Hence $\partial_t^2\log\widehat v(t,L)$ lies between the two
roots of this quadratic. In both cases,
\[
 \partial_t^2\log\widehat v(t,L)
 +\frac2{n-2}\widehat v(t,L)^{\frac4{n-2}}
 \geq
 \frac{\Psi(t,L)}{n-2}\widehat v(t,L)^{\frac2{n-2}}.
\]

Differentiating the boundary condition also gives
\[
 \left.\partial_t\partial_\beta
  \log\bigl(\sin^{\frac{n-2}{2}}\beta\,\widehat v\bigr)
 \right|_{\beta=L}
 =\frac2{n-2}\widehat v(t,L)^{\frac2{n-2}}
       \partial_t\log\widehat v(t,L).
\]
Using~\eqref{eq:flat-trace-log-laplacian} and
$\partial_\beta(\sin^2\beta)|_{\beta=L}=0$ in the definition
of $\Psi$, we obtain
\[
 \begin{aligned}
 \partial_\beta\Psi(t,L)
 =-2\bigg[
 &\left(1-\frac1{n-2}\right)
          \Psi(t,L)\widehat v(t,L)^{\frac2{n-2}}
 +\partial_t^2\log\widehat v(t,L)\\
 &+\frac2{n-2}\widehat v(t,L)^{\frac4{n-2}}
 \bigg]
 \leq-2\Psi(t,L)\widehat v(t,L)^{\frac2{n-2}}<0.
 \end{aligned}
\]
This contradicts $\partial_\beta\Psi(t,L)\geq0$ at a maximum
on $\R\times(0,L]$. Therefore $\Psi\leq0$ everywhere.
At $\beta=L$, the inequality $\Psi(t,L)\leq0$ gives
\eqref{flatbase:gradient}.
\end{proof}

We use Proposition~\ref{flatbase:trace-gradient} to compare
$\partial_\beta\widehat v/\sin\beta$ with
$\widehat v^{n/(n-2)}$ on $\R\times(0,L]$. Their
$\beta$ derivatives at $L$ give the upper bound for
$\partial_t^2\widehat v(t,L)$. For the other inequality, we use
the quotient of two positive solutions with the same linear
boundary condition.

\begin{proposition}\label{flatbase:pressure}
Let $n\geq5$ and let $\widehat v$ satisfy the hypotheses of
Proposition~\ref{flatbase:selection}. For every $t\in\R$,
\begin{align}
 \partial_t^2\widehat v(t,L)
 &\leq\frac{(n-2)^2}{4}\widehat v(t,L)
   -\frac n{n-2}\widehat v(t,L)^{\frac{n+2}{n-2}},
 \label{flatbase:upper}\\
 \partial_t^2\bigl(\widehat v(t,L)^{-\frac2{n-2}}\bigr)
 &<\widehat v(t,L)^{-\frac2{n-2}}.
 \label{flatbase:lower}
\end{align}
\end{proposition}

\begin{proof}
Expanding the gradient bound in
Proposition~\ref{flatbase:trace-gradient} gives
\[
 \widehat v_t^2+\widehat v_\beta^2
 +(n-2)\cot\beta\,\widehat v\widehat v_\beta
 \leq\frac{(n-2)^2}{4}\widehat v^2.
\]
Indeed, the bound is
\[
 \sin^2\beta\left[
 \left(\frac{\widehat v_t}{\widehat v}\right)^2+
 \left(\frac{\widehat v_\beta}{\widehat v}
              +\frac{n-2}{2}\cot\beta\right)^2\right]
 \leq\frac{(n-2)^2}{4},
\]
and $\csc^2\beta-\cot^2\beta=1$.
Differentiating the interior equation in
\eqref{flatbase:aux-equation} in $\beta$ and applying the
chain rule to $\widehat v^{n/(n-2)}$ gives
\[
 \begin{aligned}
 &\left(-\partial_t^2-\partial_\beta^2
       -n\cot\beta\,\partial_\beta+\frac{n^2}{4}\right)
 \left(\frac{\widehat v_\beta}{\sin\beta}
                     -\widehat v^{\frac n{n-2}}\right)\\
 &\quad=\frac{2n}{(n-2)^2}
 \widehat v^{\frac{4-n}{n-2}}
 \left[
 \widehat v_t^2+\widehat v_\beta^2
 +(n-2)\cot\beta\,\widehat v\widehat v_\beta
 -\frac{(n-2)^2}{4}\widehat v^2
 \right]\leq0.
 \end{aligned}
\]
The identity
\[
 \frac{\widehat v_\beta(t,\beta)}{\sin\beta}
 =\frac{\beta}{\sin\beta}
   \int_0^1\partial_\beta^2\widehat v(t,s\beta)\,\dd s
\]
and evenness at $\beta=0$ give a smooth even extension of the
quotient. The function
$\widehat v_\beta/\sin\beta-\widehat v^{n/(n-2)}$
can consequently be regarded as a smooth cylindrically symmetric
function on $\R\times\Sph^{n+1}_+$ with product metric
$\dd t^2+g_{\Sph^{n+1}}$, where $g_{\Sph^{n+1}}$ is the
round metric of radius $1$. It is bounded by
Proposition~\ref{flatbase:selection}. Since
$\Delta_{\Sph^{n+1}_+}=\partial_\beta^2+n\cot\beta\,\partial_\beta$
on functions of $\beta$, we obtain
\begin{equation}\label{flatbase:upper-comparison-equation}
 \left(-\partial_t^2-\Delta_{\Sph^{n+1}_+}+\frac{n^2}{4}\right)
 \left(\frac{\widehat v_\beta}{\sin\beta}
                     -\widehat v^{\frac n{n-2}}\right)\leq0.
\end{equation}
By the boundary condition in~\eqref{flatbase:aux-equation}
and $\sin L=1$,
\[
 \left.
 \left(\frac{\partial_\beta\widehat v}{\sin\beta}
                   -\widehat v^{\frac n{n-2}}\right)
 \right|_{\beta=L}
 =\partial_\beta\widehat v(t,L)
                   -\widehat v(t,L)^{\frac n{n-2}}=0.
\]
For $\varepsilon>0$, the operator in
\eqref{flatbase:upper-comparison-equation} applied to
$\widehat v_\beta/\sin\beta-\widehat v^{n/(n-2)}
-\varepsilon\cosh t$ is at most
$-\varepsilon(n^2/4-1)\cosh t<0$. The function is nonpositive
on $[-R,R]\times\partial\Sph^{n+1}_+$ by its value at $\beta=L$
and on $\{-R,R\}\times\Sph^{n+1}_+$ for sufficiently large $R$
by boundedness of $\widehat v_\beta/\sin\beta-\widehat v^{n/(n-2)}$.
As in Proposition~\ref{flatbase:selection}, the maximum principle
on $[-R,R]\times\Sph^{n+1}_+$ and the limit
$\varepsilon\to0^+$ give
\[
 \frac{\widehat v_\beta(t,\beta)}{\sin\beta}
       -\widehat v(t,\beta)^{\frac n{n-2}}\leq0
 \qquad\text{for }(t,\beta)\in\R\times[0,L].
\]
Since equality holds at $\beta=L$,
\[
 \begin{aligned}
 0&\leq
 \left.\partial_\beta
 \left(\frac{\widehat v_\beta}{\sin\beta}
                       -\widehat v^{\frac n{n-2}}\right)
 \right|_{\beta=L}\\
 &=\frac{(n-2)^2}{4}\widehat v(t,L)
   -\partial_t^2\widehat v(t,L)
   -\frac n{n-2}\widehat v(t,L)^{\frac{n+2}{n-2}}.
 \end{aligned}
\]
Here we used $\cot L=\cos L=0$, $\sin L=1$ and both equations
in~\eqref{flatbase:aux-equation}. This proves~\eqref{flatbase:upper}.

To prove~\eqref{flatbase:lower}, define
\[
 \begin{aligned}
 h_\pm(t,\beta)
 &=\frac{n-2}{2}\widehat v(t,\beta)
   -\cot\beta\,\widehat v_\beta(t,\beta)
   \pm\widehat v_t(t,\beta),\\
 f_\pm(t,\beta)&=e^{\pm t}\sin\beta\,h_\pm(t,\beta).
 \end{aligned}
\]
We derive their equations directly from~\eqref{flatbase:aux-equation}.
Differentiating this equation in $\beta$ and using
$(\cot\beta)'=-\csc^2\beta$ and
$(\cot\beta)''=2\cot\beta\csc^2\beta$ gives
\[
 \begin{aligned}
 &\left(\partial_t^2+\partial_\beta^2
          +n\cot\beta\,\partial_\beta-\frac{n^2-4}{4}\right)
       \left(\frac{n-2}{2}\widehat v
                         -\cot\beta\,\widehat v_\beta\right)\\
 &\quad=2\widehat v_{\beta\beta}
       +2(n-2)\cot\beta\,\widehat v_\beta
       -\frac{(n-2)^2}{2}\widehat v
       =-2\widehat v_{tt}.
 \end{aligned}
\]
Differentiating~\eqref{flatbase:aux-equation} in $t$ also gives
\[
 \begin{aligned}
 &\left(\partial_t^2+\partial_\beta^2
          +n\cot\beta\,\partial_\beta-\frac{n^2-4}{4}\right)
                  \widehat v_t\\
 &\quad=2\cot\beta\,\widehat v_{t\beta}-(n-2)\widehat v_t\\
 &\quad=-2\partial_t\left(\frac{n-2}{2}\widehat v
                            -\cot\beta\,\widehat v_\beta\right).
 \end{aligned}
\]
Adding and subtracting these identities yields, on $\R\times(0,L)$,
\begin{equation}\label{eq:flat-pressure-h-interior}
 h_{\pm,tt}\pm2h_{\pm,t}+h_{\pm,\beta\beta}
 +n\cot\beta\,h_{\pm,\beta}
 -\frac{n^2-4}{4}h_\pm=0.
\end{equation}
For $f_\pm=e^{\pm t}\sin\beta\,h_\pm$, the product rule gives
\[
 \begin{aligned}
 f_{\pm,tt}
 &=e^{\pm t}\sin\beta\,
       (h_{\pm,tt}\pm2h_{\pm,t}+h_\pm),\\
 f_{\pm,\beta\beta}+(n-2)\cot\beta\,f_{\pm,\beta}
 &=e^{\pm t}\sin\beta\,
       \bigl[h_{\pm,\beta\beta}+n\cot\beta\,h_{\pm,\beta}\\
 & \hspace{2em}
       +\bigl((n-2)\cot^2\beta-1\bigr)h_\pm\bigr].
 \end{aligned}
\]
Since
\[
 (n-2)\cot^2\beta-
 \left(\frac{(n-2)^2}{4}+(n-2)\csc^2\beta\right)
 =-\frac{n^2-4}{4},
\]
equation~\eqref{eq:flat-pressure-h-interior} gives, on
$\R\times(0,L)$,
\begin{equation}\label{eq:flat-pressure-f-interior}
 \begin{aligned}
 f_{\pm,tt}+f_{\pm,\beta\beta}
 &+(n-2)\cot\beta\,f_{\pm,\beta}\\
 &-\left(\frac{(n-2)^2}{4}+(n-2)\csc^2\beta\right)f_\pm=0.
 \end{aligned}
\end{equation}

At $\beta=L$,
\[
 h_\pm(t,L)=\frac{n-2}{2}\widehat v(t,L)
                         \pm\partial_t\widehat v(t,L).
\]
Differentiating the definition of $h_\pm$ in $\beta$ and
using~\eqref{flatbase:aux-equation} gives
\[
 \begin{aligned}
 h_{\pm,\beta}(t,L)
 &=\frac n2\widehat v(t,L)^{\frac n{n-2}}
   \pm\frac n{n-2}\widehat v(t,L)^{\frac2{n-2}}
                         \partial_t\widehat v(t,L)\\
 &=\frac n{n-2}\widehat v(t,L)^{\frac2{n-2}}h_\pm(t,L).
 \end{aligned}
\]
Since $\cos L=0$,
\begin{equation}\label{eq:flat-pressure-f-boundary}
 f_{\pm,\beta}(t,L)
 =\frac n{n-2}\widehat v(t,L)^{\frac2{n-2}}f_\pm(t,L)
 \qquad(t\in\R).
\end{equation}

Since $\cot\beta\,\widehat v_\beta
=\cos\beta\,(\widehat v_\beta/\sin\beta)$, this term extends
smoothly and evenly across $\beta=0$.
Thus $h_\pm$ are smooth on $\R\times\Sph^{n+1}_+$ and
\eqref{eq:flat-pressure-h-interior} holds there with
$\Delta_{\Sph^{n+1}_+}$ in place of
$\partial_\beta^2+n\cot\beta\,\partial_\beta$.
The integral formula for $\widehat v_\beta/\sin\beta$ and
Proposition~\ref{flatbase:selection} give uniform bounds for
$h_\pm$ and their derivatives on $\R\times\Sph^{n+1}_+$.

Let $m>0$ be a lower bound for
$\widehat v(t,L)$ on $\R$. By~\eqref{flatbase:gradient},
\[
 \begin{aligned}
 h_\pm(t,L)
 &\geq\widehat v(t,L)
       \left(\frac{n-2}{2}
             -\left|\frac{\partial_t\widehat v(t,L)}
                            {\widehat v(t,L)}\right|\right)\\
 &\geq m\left(\frac{n-2}{2}
       -\sqrt{\frac{(n-2)^2}{4}-m^{\frac4{n-2}}}\right)>0.
 \end{aligned}
\]
This bound is independent of $t$ and of the choice of sign.
Suppose that one of these functions has
$\inf_{\R\times[0,L]}h_\pm\leq0$. Choose $(t_j,\beta_j)$ such that
\[
 h_\pm(t_j,\beta_j)\longrightarrow
                \inf_{\R\times[0,L]}h_\pm\leq0.
\]
By Proposition~\ref{flatbase:selection} and the definition
of $h_\pm$, we may pass to a subsequence for which
$\beta_j\to\beta_\infty\in[0,L]$ and $h_\pm(t+t_j,\beta)$
converges in $C^2([-R,R]\times\Sph^{n+1}_+)$ for every $R>0$.
Here convergence through $\beta=0$ follows from the integral
formula for $\widehat v_\beta/\sin\beta$. We denote the limit
again by $h_\pm$.
It attains its minimum at $(0,\beta_\infty)$, where
$h_\pm(0,\beta_\infty)\leq0$. The positive lower bound for
$h_\pm(t,L)$ gives $\beta_\infty<L$. At this minimum,
\[
 \partial_t h_\pm=0,\qquad
 (\partial_t^2+\Delta_{\Sph^{n+1}})h_\pm\geq0.
\]
If $h_\pm(0,\beta_\infty)<0$, then
\[
 (\partial_t^2+\Delta_{\Sph^{n+1}})h_\pm
 \pm2\partial_t h_\pm-\frac{n^2-4}{4}h_\pm>0
\]
there, contradicting~\eqref{eq:flat-pressure-h-interior}.
If $h_\pm(0,\beta_\infty)=0$, the strong minimum principle
applied to~\eqref{eq:flat-pressure-h-interior} gives
$h_\pm\equiv0$, contradicting
the positive bound on $\R\times\partial\Sph^{n+1}_+$.
Consequently
\[
 0<\inf_{\R\times[0,L]}h_\pm
 \leq\sup_{\R\times[0,L]}h_\pm<\infty.
\]

Multiplying~\eqref{eq:flat-pressure-f-interior} for $f_\pm$
by $f_\mp$ and subtracting cancels the zeroth-order terms.
Since $\omega=\sin^{n-2}\beta$, we obtain
\[
 \partial_t\!\left[\omega(f_-f_{+,t}-f_+f_{-,t})\right]
 +\partial_\beta\!\left[\omega(f_-f_{+,\beta}
                                      -f_+f_{-,\beta})\right]=0.
\]
Using
\[
 f_+f_-=\sin^2\beta\,h_+h_-,
 \qquad
 \log(f_+/f_-)=2t+\log(h_+/h_-),
\]
and dividing by $2$, we obtain
\begin{equation}\label{eq:flat-pressure-flux}
 \begin{aligned}
 &\partial_t\!\left[\sin^n\beta\,h_+h_-
                 \left(1+\frac12\partial_t
                                  \log\frac{h_+}{h_-}\right)\right]\\
 &\quad+\frac12\partial_\beta\!\left[
       \sin^n\beta\,h_+h_-\,
       \partial_\beta\log\frac{h_+}{h_-}\right]=0.
 \end{aligned}
\end{equation}
The boundary condition~\eqref{eq:flat-pressure-f-boundary} gives
\[
 \partial_\beta\log\frac{h_+}{h_-}(t,L)=0.
\]
The uniform positive bounds and derivative bounds for $h_\pm$
give constants $c_1,c_2>0$ such that
\[
 c_1\beta^n
 \leq\sin^n\beta\,h_+(t,\beta)h_-(t,\beta)
 \leq c_2\beta^n
 \qquad(t\in\R,\ 0<\beta\leq L).
\]
They also give
\[
 \partial_t\log(h_+/h_-)=O(1),\qquad
 \partial_\beta\log(h_+/h_-)=O(\beta)
 \qquad\text{as }\beta\to0,
\]
uniformly in $t$. The second estimate follows from evenness
at $\beta=0$.

For $(t,\beta)\in\R\times[0,L]$, define
\[
 k(t,\beta)=\int_0^\beta\sin^n s\,h_+(t,s)h_-(t,s)
 \left(1+\frac12\partial_t
                  \log\frac{h_+(t,s)}{h_-(t,s)}\right)\dd s.
\]
Then $k=O(\beta^{n+1})$ uniformly in $t$, and
\[
 \begin{aligned}
 k_\beta(t,\beta)
 &=\sin^n\beta\,h_+h_-
       \left(1+\frac12\partial_t\log\frac{h_+}{h_-}\right),\\
 k_t(t,\beta)
 &=-\frac12\sin^n\beta\,h_+h_-\,
                          \partial_\beta\log\frac{h_+}{h_-}.
 \end{aligned}
\]
The second identity follows by differentiating the integral
and using~\eqref{eq:flat-pressure-flux}. The term at $s=0$
vanishes because it is $O(s^{n+1})$. Dividing the two identities
by $\sin^n\beta\,h_+h_-$ and differentiating gives
\begin{equation}\label{eq:flat-pressure-k-equation}
 \partial_t\left(\frac{k_t}{\sin^n\beta\,h_+h_-}\right)
 +\partial_\beta\left(\frac{k_\beta}{\sin^n\beta\,h_+h_-}\right)=0
 \quad\text{on }\R\times(0,L).
\end{equation}
We have $k(t,0)=0$. Since~\eqref{eq:flat-pressure-f-boundary}
gives $\partial_\beta\log(h_+/h_-)(t,L)=0$, the formula for
$k_t$ gives $k_t(t,L)=0$. Hence $k(t,L)=J$ for a real
constant $J$.

We claim that $\min(0,J)\leq k\leq\max(0,J)$. The formula
for $k_\beta$ and the uniform bounds for $h_\pm$ and their
derivatives give $|k_\beta|\leq C$ on $\R\times[0,L]$. Together
with $k(t,0)=0$ and $k(t,L)=J$, this gives
\[
 |k(t,\beta)|\leq C\beta,\qquad
 |k(t,\beta)-J|\leq C(L-\beta).
\]
Suppose that $\sup_{\R\times[0,L]}k>\max(0,J)$ and choose
$(t_j,\beta_j)$ such that $k(t_j,\beta_j)\to\sup_{\R\times[0,L]}k$.
The two bounds give $\beta_j\in[\delta,L-\delta]$ for some
$\delta>0$ and all sufficiently large $j$. After passing to a
subsequence, $\beta_j\to\beta_\infty\in[\delta,L-\delta]$.
By Proposition~\ref{flatbase:selection} and the formulas for
$h_\pm$ and $k$, the translates $h_\pm(t+t_j,\beta)$ and
$k(t+t_j,\beta)$ converge in $C^2([-R,R]\times[0,L])$ for every
$R>0$ along a further subsequence. At $\beta=0$, this convergence
follows from the integral formula for $\widehat v_\beta/\sin\beta$.
The limiting $h_\pm$ are bounded below by positive constants.
The limit of $k(t+t_j,\beta)$ therefore satisfies
\eqref{eq:flat-pressure-k-equation} with the limiting
coefficients and has the same values $0$ and $J$ at $\beta=0$
and $\beta=L$. Its maximum is attained at $(0,\beta_\infty)$.
The strong maximum principle on $\R\times(0,L)$ makes the limit
constant, which contradicts its value $0$ at $\beta=0$.
The same argument applied to $-k$, with boundary value $-J$
at $\beta=L$, proves the lower bound.

If $J<0$, the function $k$ is nonconstant because its boundary
values are $0$ and $J$. The strong minimum principle applied
to~\eqref{eq:flat-pressure-k-equation} gives $k>J$ for
$0<\beta<L$. The coefficient in this equation is smooth and
positive near each point $(t,L)$. The boundary point lemma at
the minimum $k(t,L)=J$ therefore gives
\[
 k_\beta(t,L)<0.
\]
If $J=0$, the bound $\min(0,J)\leq k\leq\max(0,J)$ gives
$k\equiv0$ and $k_\beta(t,L)=0$. In either case,
\[
 \begin{aligned}
 0&\geq\frac{k_\beta(t,L)}{h_+(t,L)h_-(t,L)}\\
 &=\frac{\dd}{\dd t}\left[
 t+\frac12\log
 \frac{\frac{n-2}{2}\widehat v(t,L)+\partial_t\widehat v(t,L)}
      {\frac{n-2}{2}\widehat v(t,L)-\partial_t\widehat v(t,L)}
 \right].
 \end{aligned}
\]
The uniform bounds for $h_\pm(t,L)$ bound the quotient
above and below by positive constants. Hence
$t+\frac12\log(h_+(t,L)/h_-(t,L))=t+O(1)$ tends to
$+\infty$ as $t\to+\infty$, contradicting its nonpositive
derivative. Hence $J>0$.

For $J>0$, the strong maximum principle for
\eqref{eq:flat-pressure-k-equation} gives $k<J$ on
$\R\times(0,L)$, since $k(t,0)=0$. The boundary point lemma
at the maximum $k(t,L)=J$ gives
$k_\beta(t,L)>0$. Dividing by $h_+(t,L)h_-(t,L)>0$ and
differentiating the logarithm gives
\[
 1+\frac{
 \frac{n-2}{2}
 \left[\widehat v(t,L)\,\partial_t^2\widehat v(t,L)
                 -(\partial_t\widehat v(t,L))^2\right]}
 {\frac{(n-2)^2}{4}\widehat v(t,L)^2
                 -(\partial_t\widehat v(t,L))^2}>0.
\]
The denominator equals $h_+(t,L)h_-(t,L)>0$. Multiplying by
this product and dividing by
$\frac{n-2}{2}\widehat v(t,L)^2$ yields
\[
 \frac{\partial_t^2\widehat v(t,L)}{\widehat v(t,L)}
 >\frac n{n-2}
       \left(\frac{\partial_t\widehat v(t,L)}
                        {\widehat v(t,L)}\right)^2
 -\frac{n-2}{2}.
\]
Finally,
\[
 \begin{aligned}
 \partial_t^2\bigl(\widehat v(t,L)^{-\frac2{n-2}}\bigr)
 &=\frac2{n-2}\widehat v(t,L)^{-\frac2{n-2}}
 \left[
 \frac n{n-2}
       \left(\frac{\partial_t\widehat v(t,L)}
                        {\widehat v(t,L)}\right)^2
 -\frac{\partial_t^2\widehat v(t,L)}{\widehat v(t,L)}
 \right]\\
 &<\widehat v(t,L)^{-\frac2{n-2}},
 \end{aligned}
\]
which proves~\eqref{flatbase:lower}.
\end{proof}

\begin{corollary}\label{flatbase:secant}
Let $n\geq5$, and let $\widehat v_1,\widehat v_2$ satisfy the
hypotheses of Proposition~\ref{flatbase:selection}, with
constants $m,M>0$ such that
\[
 m\leq\widehat v_i(t,L)\leq M
 \qquad(t\in\R,\ i=1,2).
\]
Suppose either that both traces are even in $t$, or that, for
some $t_0\in\R$ and some solution $\widehat v$ satisfying the
same hypotheses,
\[
 \widehat v_1(t,L)=\widehat v(t_0+t,L),\qquad
 \widehat v_2(t,L)=\widehat v(t_0-t,L).
\]
The coefficient
\begin{equation}\label{flatbase:secant-coefficient}
 b(t)=\frac n{n-2}\int_0^1
 \bigl[(1-\lambda)\widehat v_1(t,L)
                 +\lambda\widehat v_2(t,L)\bigr]^{\frac2{n-2}}
 \,\dd\lambda
\end{equation}
is smooth, positive and even. It satisfies
\begin{align}
 0<b\leq\frac n2,\qquad
 b''&\geq2\frac{(b')^2}{b}-b,
 \label{flatbase:secant-lower}\\
 \frac{b''}{b}
 &\leq\frac{n-2}{2}-\frac{2b^2}{n}
          -\frac{n-4}{2}\left(\frac{b'}b\right)^2,
 \label{flatbase:secant-upper}\\
 \left|\frac{b'}b\right|^2
 &\leq1-\frac{4b^2}{n^2},\notag\\
 -\tanh t\leq\frac{b'}b
 &\leq\tanh\left(\frac{n-2}{2}t\right)
 \qquad(t\geq0).
 \label{flatbase:history}
\end{align}
Taking $\widehat v_1=\widehat v_2=\widehat v$, with $\widehat v$
even in $t$, gives the same estimates for the coefficient
$\frac n{n-2}\widehat v(t,L)^{2/(n-2)}$ in the linearized
boundary condition.
\end{corollary}

\begin{proof}
In the first case, the integrand in~\eqref{flatbase:secant-coefficient}
is even in $t$ for every $\lambda$. In the second case, reflection
in $t$ interchanges $\widehat v_1(t,L)$ and $\widehat v_2(t,L)$.
Changing the integration variable from $\lambda$ to $1-\lambda$
in~\eqref{flatbase:secant-coefficient} gives $b(-t)=b(t)$.

Since $\widehat v_i(t,L)\geq m>0$ for $i=1,2$, we have
\[
 (1-\lambda)\widehat v_1(t,L)+\lambda\widehat v_2(t,L)\geq m
 \qquad(0\leq\lambda\leq1).
\]
This bound and the smoothness of $\widehat v_i(t,L)$ allow
differentiation under the integral in~\eqref{flatbase:secant-coefficient}.
The same formula gives $b(t)\geq\frac n{n-2}m^{\frac2{n-2}}>0$.
By the boundary condition in~\eqref{flatbase:aux-equation}
and the fundamental theorem of calculus in $\lambda$,
\[
 \partial_\beta(\widehat v_2-\widehat v_1)(t,L)
 =b(t)(\widehat v_2(t,L)-\widehat v_1(t,L)).
\]

By~\eqref{flatbase:lower}, each $\widehat v_i(t,L)$ satisfies
\[
 \partial_t^2\widehat v_i(t,L)
 \geq\frac n{n-2}
       \frac{(\partial_t\widehat v_i(t,L))^2}{\widehat v_i(t,L)}
       -\frac{n-2}{2}\widehat v_i(t,L)
 \qquad(i=1,2).
\]
For $0\leq\lambda\leq1$, weighted Cauchy--Schwarz gives
\[
 \begin{aligned}
 &\frac{\bigl[(1-\lambda)\partial_t\widehat v_1(t,L)
                    +\lambda\partial_t\widehat v_2(t,L)\bigr]^2}
 {(1-\lambda)\widehat v_1(t,L)+\lambda\widehat v_2(t,L)}\\
 &\quad\leq
 (1-\lambda)\frac{(\partial_t\widehat v_1(t,L))^2}
                  {\widehat v_1(t,L)}
 +\lambda\frac{(\partial_t\widehat v_2(t,L))^2}
                  {\widehat v_2(t,L)}.
 \end{aligned}
\]
Multiplying the trace inequalities by $1-\lambda$ and $\lambda$,
respectively, and adding gives
\[
 \begin{aligned}
 &\partial_t^2
 \bigl[(1-\lambda)\widehat v_1(t,L)+\lambda\widehat v_2(t,L)\bigr]\\
 &\quad\geq\frac n{n-2}
 \frac{\bigl[(1-\lambda)\partial_t\widehat v_1(t,L)
                  +\lambda\partial_t\widehat v_2(t,L)\bigr]^2}
 {(1-\lambda)\widehat v_1(t,L)+\lambda\widehat v_2(t,L)}\\
 &\qquad-\frac{n-2}{2}
 \bigl[(1-\lambda)\widehat v_1(t,L)+\lambda\widehat v_2(t,L)\bigr].
 \end{aligned}
\]
Differentiating the power $2/(n-2)$ therefore gives
\[
 \begin{aligned}
 &\partial_t^2
 \bigl[(1-\lambda)\widehat v_1(t,L)+\lambda\widehat v_2(t,L)\bigr]^{
                                                        \frac2{n-2}}\\
 &\quad\geq
 2\frac{\left\{\partial_t
 \bigl[(1-\lambda)\widehat v_1(t,L)+\lambda\widehat v_2(t,L)\bigr]^{
                                                        \frac2{n-2}}
 \right\}^2}
 {\bigl[(1-\lambda)\widehat v_1(t,L)
                   +\lambda\widehat v_2(t,L)\bigr]^{\frac2{n-2}}}\\
 &\qquad-
 \bigl[(1-\lambda)\widehat v_1(t,L)+\lambda\widehat v_2(t,L)\bigr]^{
                                                        \frac2{n-2}}.
 \end{aligned}
\]
Cauchy--Schwarz in the integral defining $b'$ gives
\[
 \frac n{n-2}\int_0^1
 \frac{\left\{\partial_t
 \bigl[(1-\lambda)\widehat v_1(t,L)+\lambda\widehat v_2(t,L)\bigr]^{
                                                        \frac2{n-2}}
 \right\}^2}
 {\bigl[(1-\lambda)\widehat v_1(t,L)
                   +\lambda\widehat v_2(t,L)\bigr]^{\frac2{n-2}}}
 \,\dd\lambda
 \geq\frac{(b')^2}{b}.
\]
Integrating the bound for the second derivative of the
integrand in~\eqref{flatbase:secant-coefficient} and using the
Cauchy--Schwarz estimate for $b'$ gives the differential
inequality in~\eqref{flatbase:secant-lower}.

For the upper estimate, use~\eqref{flatbase:upper} and convexity
of the power $(n+2)/(n-2)$. For every $\lambda\in[0,1]$,
\[
 \begin{aligned}
 &\partial_t^2
 \bigl[(1-\lambda)\widehat v_1(t,L)+\lambda\widehat v_2(t,L)\bigr]\\
 &\quad\leq\frac{(n-2)^2}{4}
 \bigl[(1-\lambda)\widehat v_1(t,L)+\lambda\widehat v_2(t,L)\bigr]\\
 &\qquad-\frac n{n-2}
 \bigl[(1-\lambda)\widehat v_1(t,L)+\lambda\widehat v_2(t,L)\bigr]^{
                                                     \frac{n+2}{n-2}}.
 \end{aligned}
\]
Differentiating the power $2/(n-2)$ gives
\[
 \begin{aligned}
 &\partial_t^2
 \bigl[(1-\lambda)\widehat v_1(t,L)+\lambda\widehat v_2(t,L)\bigr]^{
                                                        \frac2{n-2}}\\
 &\quad\leq\frac{n-2}{2}
 \bigl[(1-\lambda)\widehat v_1(t,L)+\lambda\widehat v_2(t,L)\bigr]^{
                                                        \frac2{n-2}}\\
 &\qquad-\frac{2n}{(n-2)^2}
 \bigl[(1-\lambda)\widehat v_1(t,L)+\lambda\widehat v_2(t,L)\bigr]^{
                                                        \frac6{n-2}}\\
 &\qquad-\frac{n-4}{2}
 \frac{\left\{\partial_t
 \bigl[(1-\lambda)\widehat v_1(t,L)+\lambda\widehat v_2(t,L)\bigr]^{
                                                        \frac2{n-2}}
 \right\}^2}
 {\bigl[(1-\lambda)\widehat v_1(t,L)
                   +\lambda\widehat v_2(t,L)\bigr]^{\frac2{n-2}}}.
 \end{aligned}
\]
Jensen's inequality for the cube and
\eqref{flatbase:secant-coefficient} give
\[
 \begin{aligned}
 &\frac n{n-2}\int_0^1
 \bigl[(1-\lambda)\widehat v_1(t,L)+\lambda\widehat v_2(t,L)\bigr]^{
                                                        \frac6{n-2}}
 \,\dd\lambda\\
 &\quad\geq\frac n{n-2}
 \left\{\int_0^1
 \bigl[(1-\lambda)\widehat v_1(t,L)+\lambda\widehat v_2(t,L)\bigr]^{
                                                        \frac2{n-2}}
 \,\dd\lambda\right\}^3\\
 &\quad=\frac{(n-2)^2}{n^2}b^3.
 \end{aligned}
\]
Integrate the upper bound for the second derivative of the
integrand in~\eqref{flatbase:secant-coefficient}. Use the
lower bound for the integral of the power $6/(n-2)$ and the
Cauchy--Schwarz bound for $b'$, which enters with coefficient
$-(n-4)/2<0$. This gives
\[
 b''\leq\frac{n-2}{2}b-\frac2n b^3
                    -\frac{n-4}{2}\frac{(b')^2}{b}.
\]
Dividing by $b>0$ proves~\eqref{flatbase:secant-upper}.

Combining~\eqref{flatbase:secant-lower}
and~\eqref{flatbase:secant-upper} gives
\[
 2\left(\frac{b'}b\right)^2-1
 \leq\frac{b''}{b}
 \leq\frac{n-2}{2}-\frac{2b^2}{n}
               -\frac{n-4}{2}\left(\frac{b'}b\right)^2.
\]
Rearranging yields
\[
 \left(\frac{b'}b\right)^2+\frac{4b^2}{n^2}\leq1.
\]
Hence $b\leq n/2$ and $|b'/b|<1$. Subtracting $(b'/b)^2$
from the two bounds for $b''/b$ gives
\[
 \left(\frac{b'}b\right)^2-1
 \leq\left(\frac{b'}b\right)'
 \leq\frac{n-2}{2}
       \left[1-\left(\frac{b'}b\right)^2\right].
\]
Consequently
\[
 -1\leq
 \frac{\dd}{\dd t}\operatorname{arctanh}\left(\frac{b'}b\right)
 \leq\frac{n-2}{2}.
\]
Since $b$ is even, $b'(0)=0$. Integrating from $0$ to $t\geq0$
and applying the increasing function $\tanh$ proves
\eqref{flatbase:history}. For $\widehat v_1=\widehat v_2
=\widehat v$, \eqref{flatbase:secant-coefficient} reduces to
$b(t)=\frac n{n-2}\widehat v(t,L)^{2/(n-2)}$, the coefficient
in the boundary condition of~\eqref{eq:jacobi-local-pair}
in the normalization~\eqref{eq:flat-normalization-table}.
\end{proof}

\subsection{Spectral estimates}
\label{subsec:flat-spectral-estimates}

We use $\beta$, $L$ and $\omega$ as in
Section~\ref{subsec:notation}. For cylindrically symmetric
functions $f=f(\beta)$, we use the differential expression of
$\mathcal A$ from Section~\ref{subsec:linearized-spaces} and
impose the boundary condition $f'(L)=bf(L)$ for $b\geq0$.
For $n\geq5$, we estimate the first eigenvalue and the
corresponding positive eigenfunction, and give a lower bound
for the gap between the first two eigenvalues. We then prove
a resolvent identity for all $n\geq3$.

For $n\geq5$, let
\[
 \mathcal Q=\{z\in H^1(0,L):z(0)=0\}.
\]
We use real-valued functions and write $\|\cdot\|_2$ for the norm
in $L^2(0,L)$. For $z\in\mathcal Q$ and $b\in\R$, set
\begin{equation}\label{eq:flat-spectral-unweighted-form}
 q_b[z]=\int_0^L
 \left(|z'|^2+\frac{(n-2)(n-4)}4\csc^2\beta\,|z|^2\right)
 \dd\beta-b|z(L)|^2.
\end{equation}
We write $q_b[\cdot,\cdot]$ for the associated bilinear form.
Hardy's inequality and $\sin\beta\geq2\beta/\pi$ on $[0,L]$
give
\[
 \int_0^L\csc^2\beta\,|z|^2\,\dd\beta
 \leq\frac{\pi^2}{4}\int_0^L\frac{|z|^2}{\beta^2}\,\dd\beta
 \leq\pi^2\|z'\|_2^2.
\]
For the trace $z(L)$, the condition $z(0)=0$ gives, for
every $\varepsilon>0$,
\begin{equation}\label{eq:flat-spectral-unweighted-trace}
 |z(L)|^2=2\int_0^L zz'\,\dd\beta
 \leq\varepsilon\|z'\|_2^2+\varepsilon^{-1}\|z\|_2^2.
\end{equation}
For fixed $b$, choose $\varepsilon>0$ with $|b|\varepsilon<1$.
Adding a sufficiently large multiple of $\|z\|_2^2$ to
$q_b[z]$ gives a squared norm equivalent to
$\|z\|_{H^1(0,L)}^2$ on $\mathcal Q$. Thus $q_b$ is closed
and bounded below. The inclusion
$\mathcal Q\hookrightarrow L^2(0,L)$ is compact by the
compact embedding $H^1(0,L)\hookrightarrow L^2(0,L)$.

Let $A_b$ be the self-adjoint operator associated with $q_b$.
Its differential expression on $(0,L)$ is
\begin{equation}\label{eq:flat-spectral-unweighted-expression}
 -\partial_\beta^2+\frac{(n-2)(n-4)}4\csc^2\beta.
\end{equation}
The domain of $A_b$ consists of $z\in\mathcal Q$ for which
\eqref{eq:flat-spectral-unweighted-expression} applied to $z$
belongs to $L^2(0,L)$ in the sense of distributions and
$z'(L)=bz(L)$. The potential in
\eqref{eq:flat-spectral-unweighted-expression} is bounded near
$L$, so such a $z$ is $H^2$ near $L$ and $z'(L)$ is defined.
Integration by parts in $q_b[z,\phi]$ for
$\phi\in\mathcal Q$ vanishing near $0$ gives the boundary
term $(z'(L)-bz(L))\phi(L)$. This gives the condition at $L$.
Conversely, if this term vanishes, the integration by parts
identity extends to all $\phi\in\mathcal Q$ by density of
the tests vanishing near $0$, proved below.

Write $\mu_j(b)$, $j\geq0$, for the eigenvalues of $A_b$ in
increasing order. By the min--max principle, with $E$ ranging
over linear subspaces of $\mathcal Q$,
\begin{equation}\label{eq:flat-spectral-min-max}
 \mu_j(b)=
 \min_{\substack{E\subset\mathcal Q\\ \dim E=j+1}}
 \max_{z\in E\setminus\{0\}}\frac{q_b[z]}{\|z\|_2^2}.
\end{equation}
In derivatives with respect to $b$, we keep $n$ and the
spatial variable fixed.

For $n\geq3$, identify the functions in
$H^1_{\mathrm{sym}}(\Sph^{n-1}_+)$ with functions of $\beta$,
where the subscript $\mathrm{sym}$ denotes cylindrical symmetry
as in Section~\ref{subsec:linearized-spaces}. On such functions,
the spherical $L^2$ inner product equals $|\Sph^{n-2}|$ times
the inner product in $L^2((0,L),\omega\,\dd\beta)$. For $b\geq0$,
let $\mathcal A_b$ be the self-adjoint operator in this weighted
space associated with the form on $H^1_{\mathrm{sym}}(\Sph^{n-1}_+)$
given by
\begin{equation}\label{eq:flat-spectral-weighted-form}
 a_b[f]=\int_0^L\omega(\beta)
 \left(|f'|^2+\frac{(n-2)^2}{4}|f|^2\right)\dd\beta
 -b|f(L)|^2.
\end{equation}
We write $a_b[\cdot,\cdot]$ for the associated bilinear form.
For $0<\delta<L/2$, Cauchy--Schwarz in
$f(L)=f(s)+\int_s^L f'(\beta)\,\dd\beta$ and averaging over
$L-\delta\leq s\leq L$ give
\[
 |f(L)|^2
 \leq\frac2\delta\int_{L-\delta}^L|f|^2\,\dd\beta
      +2\delta\int_{L-\delta}^L|f'|^2\,\dd\beta.
\]
Since $\omega$ has a positive lower bound on $[L/2,L]$,
choosing $\delta$ sufficiently small gives, for each
$\varepsilon>0$,
\begin{equation}\label{eq:flat-spectral-weighted-trace}
 |f(L)|^2
 \leq\varepsilon\int_0^L\omega|f'|^2\,\dd\beta
      +C_\varepsilon\int_0^L\omega|f|^2\,\dd\beta.
\end{equation}
For fixed $b$, take $\varepsilon$ so that $b\varepsilon<1$.
Adding a sufficiently large multiple of the squared
$L^2((0,L),\omega\,\dd\beta)$ norm to $a_b[f]$ gives a squared
norm equivalent to $\|f\|_{H^1(\Sph^{n-1}_+)}^2$. Thus $a_b$
is closed and bounded below. The compact embedding
$H^1(\Sph^{n-1}_+)\hookrightarrow L^2(\Sph^{n-1}_+)$ gives
compactness of
$H^1_{\mathrm{sym}}(\Sph^{n-1}_+)\hookrightarrow
L^2((0,L),\omega\,\dd\beta)$ under the identification above.

Integration by parts in~\eqref{eq:flat-spectral-weighted-form}
gives the interior expression
\begin{equation}\label{eq:flat-spectral-weighted-expression}
 \mathcal A_b f=-\omega^{-1}(\omega f')'
                     +\frac{(n-2)^2}{4}f
 \quad\text{on }(0,L).
\end{equation}
Since $\omega(L)=1$, integration by parts against tests
supported near $L$ gives $f'(L)=bf(L)$. If $\mathcal A_b f=\lambda f$,
the identity
\[
 a_b[f,\phi]=\lambda\int_0^L\omega f\phi\,\dd\beta
 \qquad\bigl(\phi\in H^1_{\mathrm{sym}}(\Sph^{n-1}_+)\bigr)
\]
is the weak eigenvalue equation on $\Sph^{n-1}_+$ for
$-\Delta_{\Sph^{n-1}_+}+(n-2)^2/4$ with boundary condition
$\partial_\nu f=bf$. It holds in the spherical formulation for
every smooth test function on $\Sph^{n-1}_+$, since averaging
a test function over rotations of the first $n-1$ coordinates
leaves its pairings with $f$ unchanged. Elliptic regularity for
this equation and boundary condition gives smoothness on
$\Sph^{n-1}_+$.

For $n\geq5$, multiplication by $\sin^{(n-2)/2}\beta$ is
unitary from the weighted space $L^2((0,L),\omega\,\dd\beta)$
to $L^2(0,L)$. For smooth cylindrically symmetric $f$, write
$z=\sin^{(n-2)/2}\beta\,f$. Then $z\in\mathcal Q$ and
$z(L)=f(L)$. By~\eqref{eq:flat-spectral-weighted-form},
\[
 a_b[f]=\int_0^L
 \left(\left|z'-\frac{n-2}{2}\cot\beta\,z\right|^2
                   +\frac{(n-2)^2}{4}|z|^2\right)\dd\beta
                   -b|z(L)|^2.
\]
The boundary terms from integration by parts vanish because
$\cot L=0$ and $\cot\beta\,|z|^2=O(\beta^{n-3})$ near $0$.
The resulting coefficient of $|z|^2$ is
\[
 \frac{(n-2)^2}{4}(1+\cot^2\beta)
 -\frac{n-2}{2}\csc^2\beta
 =\frac{(n-2)(n-4)}4\csc^2\beta.
\]
Hence
\begin{equation}\label{eq:flat-spectral-form-identification}
 a_b[f]=q_b[z],
 \qquad z=\sin^{\frac{n-2}{2}}\!\beta\,f.
\end{equation}

Smooth cylindrically symmetric functions are dense in
$H^1_{\mathrm{sym}}(\Sph^{n-1}_+)$ by smoothing and averaging
over rotations. Their images under multiplication by
$\sin^{(n-2)/2}\beta$ include every smooth function on $[0,L]$
that vanishes near $0$, since division by $\sin^{(n-2)/2}\beta$
gives a smooth cylindrically symmetric preimage. These functions
are dense in $\mathcal Q$. Indeed, cutting off $z\in\mathcal Q$ on
$[0,\varepsilon]$ with derivative bounded by $C/\varepsilon$
on $[\varepsilon,2\varepsilon]$ gives convergence in $H^1(0,L)$
by Hardy's inequality and
\[
 \int_0^{2\varepsilon}\frac{|z|^2}{\beta^2}\,\dd\beta
 \longrightarrow0.
\]
We smooth the cutoff functions to obtain the density in
$\mathcal Q$. Therefore~\eqref{eq:flat-spectral-form-identification}
holds on the form domains and gives the unitary equivalence of
$\mathcal A_b$ and $A_b$ for $n\geq5$.

\begin{proposition}\label{flatbase:angular}
Let $n\geq5$. For each $b\geq0$, let $e_b$ be the eigenfunction
of $A_b$ corresponding to $\mu_0(b)$, chosen positive on
$(0,L]$ and normalized by $\|e_b\|_2=1$. Set
\[
 C(b)=e_b(L)^2,\qquad r(b)=\frac{(n-2)^2}{4}-\mu_0(b).
\]
The map $b\mapsto(\mu_0(b),e_b)$ is smooth with values in
$\R\times\mathcal Q$, and
\begin{equation}\label{flatbase:angular-bounds}
 -\frac{\dd\mu_0}{\dd b}(b)=C(b),\qquad
 b^2\leq r(b)\leq b^2+\frac{n-2}{2},\qquad
 C(b)\geq2b.
\end{equation}
Moreover,
\begin{equation}\label{flatbase:angular-gap}
 bC(b)\leq2r(b),\qquad \mu_1(b)-\mu_0(b)>r(b)+n-1.
\end{equation}
The form $q_b$ restricted to $\{z\in\mathcal Q:z(L)=0\}$
has first eigenvalue $n^2/4$, with eigenfunction
$\sin^{\frac{n-2}{2}}\beta\cos\beta$. For every $b\geq0$,
\[
 \mu_0(b)<\frac{n^2}{4}<\mu_1(b).
\]
\end{proposition}

\begin{proof}
The case $j=0$ of~\eqref{eq:flat-spectral-min-max} allows us
to choose a nonnegative first eigenfunction, since taking the
absolute value preserves $q_b[z]/\|z\|_2^2$. At a zero in
$(0,L)$ its derivative would vanish, and at a zero at $L$ the
same conclusion follows from $e_b'(L)=be_b(L)$. Uniqueness
for $A_b e_b=\mu_0(b)e_b$ on $(0,L)$, using
\eqref{eq:flat-spectral-unweighted-expression}, gives $e_b>0$
on $(0,L]$ and simplicity, since $e_b(L)$ fixes both data
at $L$.

Fix $b_0\geq0$ and choose $M>0$ such that
$q_{b_0}[z]+M\|z\|_2^2$ is equivalent to $\|z\|_{H^1(0,L)}^2$
on $\mathcal Q$. By~\eqref{eq:flat-spectral-unweighted-form}
and~\eqref{eq:flat-spectral-unweighted-trace}, the map
\[
 (b,\mu,e)\longmapsto
 \bigl(q_b[e,\cdot]-\mu\langle e,\cdot\rangle_{L^2(0,L)},
       \|e\|_2^2-1\bigr)
\]
is smooth from $\R\times\R\times\mathcal Q$ to
$\mathcal Q^*\times\R$. At $(b_0,\mu_0(b_0),e_{b_0})$,
testing the first component of the derivative in $(\mu,e)$
against $e_{b_0}$ determines the variation of $\mu$. The second
component determines the projection of the variation of $e$
onto $e_{b_0}$. On $\mathcal Q\cap e_{b_0}^\perp$, the
eigenfunction expansion for $A_{b_0}$ gives
\[
 q_{b_0}[h]-\mu_0(b_0)\|h\|_2^2
 \geq \bigl(\mu_1(b_0)-\mu_0(b_0)\bigr)\|h\|_2^2.
\]
The gap is positive by simplicity. Adding
$(M+\mu_0(b_0))\|h\|_2^2$ to the left-hand side gives
$q_{b_0}[h]+M\|h\|_2^2$. Thus $q_{b_0}-\mu_0(b_0)\|\cdot\|_2^2$
is coercive on $\mathcal Q\cap e_{b_0}^\perp$, and the
Lax--Milgram theorem determines the remaining component
of the variation of $e$. The derivative is therefore invertible,
and the implicit function theorem gives a smooth local branch
of eigenvalues and $L^2(0,L)$-normalized eigenfunctions through
$(\mu_0(b_0),e_{b_0})$.

To identify this branch, use the trace bound
\eqref{eq:flat-spectral-unweighted-trace} in the difference
of the forms at $b$ and $b_0$ to obtain
\[
 |q_b[z]-q_{b_0}[z]|=|b-b_0||z(L)|^2
 \leq C|b-b_0|\bigl(q_{b_0}[z]+M\|z\|_2^2\bigr).
\]
Dividing by $\|z\|_2^2$ and applying
\eqref{eq:flat-spectral-min-max} to the shifted forms, for
$C|b-b_0|<1$ we obtain
\[
 |\mu_j(b)-\mu_j(b_0)|
 \leq C|b-b_0|(\mu_j(b_0)+M),\qquad j=0,1.
\]
Since $\mu_0(b_0)<\mu_1(b_0)$, the eigenvalue of the smooth
branch stays below $\mu_1(b)$ for $b$ near $b_0$, by continuity.
It therefore equals $\mu_0(b)$.
Continuity of $e\mapsto e(L)$ on $\mathcal Q$ follows from
\eqref{eq:flat-spectral-unweighted-trace}. The normalized
eigenfunction on the branch has positive value at $L$ for
$b$ near $b_0$, so it equals $e_b$. This proves smooth
dependence in $\R\times\mathcal Q$, including at $b_0=0$.

Differentiating
$q_b[e_b,\phi]=\mu_0(b)\langle e_b,\phi\rangle_{L^2(0,L)}$
in $b$ with $\phi\in\mathcal Q$ fixed and then taking $\phi=e_b$
gives
\[
 \frac{\dd\mu_0}{\dd b}(b)=-e_b(L)^2=-C(b).
\]
The terms containing $\partial_b e_b$ cancel by the same
form identity with $\phi=\partial_b e_b$.

The condition $e_b(0)=0$ and the normalization
$\|e_b\|_2=1$ give
\[
 0\leq\int_0^L(\partial_\beta e_b-be_b)^2\,\dd\beta
 =\int_0^L|\partial_\beta e_b|^2\,\dd\beta+b^2-bC(b).
\]
Since $\csc^2\beta\geq1$, \eqref{eq:flat-spectral-unweighted-form}
gives
\[
 \mu_0(b)\geq\frac{(n-2)(n-4)}{4}-b^2.
\]
For $j=0$, the min--max formula
\eqref{eq:flat-spectral-min-max} and the identification of forms
\eqref{eq:flat-spectral-form-identification} give
\[
 \mu_0(b)=
 \min_{\substack{f\in H^1_{\mathrm{sym}}(\Sph^{n-1}_+)\\ f\ne0}}
 \frac{a_b[f]}{\int_0^L\omega|f|^2\,\dd\beta}.
\]
For $b=0$, the minimum is $(n-2)^2/4$, attained at $f=1$
by~\eqref{eq:flat-spectral-weighted-form}. For $b>0$,
$f(\beta)=e^{b(\beta-L)}$ is Lipschitz in the spherical distance
$\beta$ and belongs to $H^1_{\mathrm{sym}}(\Sph^{n-1}_+)$.
Using this function in the quotient and setting $x=L-\beta$
in its squared $L^2((0,L),\omega\,\dd\beta)$ norm gives
\begin{equation}\label{eq:flat-spectral-rayleigh}
 \begin{gathered}
 I(b)=\int_0^L e^{-2bx}\cos^{n-2}x\,\dd x
      <\int_0^\infty e^{-2bx}\,\dd x=\frac1{2b},\\
 \mu_0(b)\leq\frac{a_b[e^{b(\beta-L)}]}{I(b)}
 =\frac{(n-2)^2}{4}+b^2-\frac b{I(b)}.
 \end{gathered}
\end{equation}
Together with $\mu_0(b)\geq(n-2)(n-4)/4-b^2$,
this gives $b^2\leq r(b)\leq b^2+(n-2)/2$.

By~\eqref{eq:flat-spectral-form-identification}, the function
$f=e_b/\sin^{(n-2)/2}\beta$ is a smooth cylindrically symmetric
eigenfunction of $\mathcal A_b$ on $\Sph^{n-1}_+$. It is
positive for $0<\beta\leq L$. Since
$\Delta_{\Sph^{n-1}_+}f-r(b)f=0$ and $r(b)\geq0$, the strong
minimum principle gives $f(0)>0$. Smoothness and evenness at
$\beta=0$ give
\begin{equation}\label{eq:flat-spectral-regular-expansion}
 \begin{aligned}
 e_b(\beta)&=\beta^{\frac{n-2}{2}}\bigl(f(0)+O(\beta^2)\bigr),\\
 (\log e_b)''&=-\frac{n-2}{2}\csc^2\beta+O(1)
                          \qquad\text{as }\beta\downarrow0.
 \end{aligned}
\end{equation}
From $A_b e_b=\mu_0(b)e_b$ and
\eqref{eq:flat-spectral-unweighted-expression},
\begin{equation}\label{eq:flat-spectral-riccati}
 (\log e_b)''
 =\frac{(n-2)(n-4)}4\csc^2\beta-\mu_0(b)
                  -\bigl((\log e_b)'\bigr)^2.
\end{equation}
At any zero of $(\log e_b)''$ in $(0,L)$, differentiation of
\eqref{eq:flat-spectral-riccati} gives
\[
 (\log e_b)^{\prime\prime\prime}=-\frac{(n-2)(n-4)}2\csc^2\beta\cot\beta<0.
\]
Thus $(\log e_b)''$, which is negative near $0$
by~\eqref{eq:flat-spectral-regular-expansion}, has no first zero
in $(0,L)$. Hence $\log e_b$ is concave on $(0,L]$.
The tangent bound at $L$, where $(\log e_b)'(L)=b$, gives
$e_b(\beta)^2\leq C(b)e^{2b(\beta-L)}$. Integrating on $(0,L)$
and using $\|e_b\|_2=1$ gives $C(b)\geq2b$ for $b>0$.
For $b=0$ this bound follows from $C(0)>0$.

Equation~\eqref{eq:flat-spectral-weighted-expression} gives
$(\omega f')'=r(b)\omega f$. Since $\omega(0)f'(0)=0$ and
$f'(L)=bf(L)$, integration on $(0,L)$ gives
\[
 bC(b)=r(b)e_b(L)\int_0^L
                 \sin^{\frac{n-2}{2}}\beta\,e_b(\beta)\,\dd\beta.
\]
Concavity of $\log e_b$ gives
$e_b(L)e_b(\beta)\leq e_b((L+\beta)/2)^2$. Consequently
\[
 bC(b)\leq r(b)\int_0^L e_b\left(\frac{L+\beta}{2}\right)^2
                   \dd\beta
 =2r(b)\int_{L/2}^L e_b(\beta)^2\,\dd\beta\leq2r(b).
\]

The identification~\eqref{eq:flat-spectral-form-identification}
takes the form domain with $f(L)=0$ onto
$\{z\in\mathcal Q:z(L)=0\}$. The Dirichlet minimum in the proof
of Lemma~\ref{flatbase:robin-resolvent}, attained at $\cos\beta$,
therefore gives the claimed eigenvalue and eigenfunction for
this restriction of $q_b$. The strict inequalities
$\mu_0(b)<n^2/4<\mu_1(b)$ follow from the same lemma, and yield
\[
 \mu_1(b)-\mu_0(b)>\frac{n^2}{4}-\mu_0(b)=r(b)+n-1.
\]
\end{proof}

We now prove the resolvent identity for $\mathcal A_b$ in all
dimensions $n\geq3$.

\begin{lemma}\label{flatbase:robin-resolvent}
Let $n\geq3$ and $b\geq0$. Write $\lambda_0<\lambda_1\leq\cdots$
for the eigenvalues of $\mathcal A_b$, and let $e_j$ be
corresponding real orthonormal eigenfunctions in
$L^2((0,L),\omega\,\dd\beta)$. Then
\[
 \lambda_0<\frac{n^2}{4}<\lambda_1.
\]
Let $\mathcal A_{b,\perp}$ be the restriction of $\mathcal A_b$
to $e_0^\perp$ with respect to the inner product
$\ip{f}{\phi}=\int_0^L\omega f\phi\,\dd\beta$. Write $\ell(f)=f(L)$
for the trace functional on $H^1_{\mathrm{sym}}(\Sph^{n-1}_+)$ and
for its restriction to $H^1_{\mathrm{sym}}(\Sph^{n-1}_+)\cap e_0^\perp$.
For $\lambda_0<\xi\leq n^2/4$,
\begin{equation}\label{flatbase:robin-resolvent-identity}
 \begin{aligned}
 \ell\bigl((\mathcal A_{b,\perp}-\tfrac{n^2}{4})^{-1}\ell\bigr)
 &=\sum_{j\geq1}\frac{e_j(L)^2}{\lambda_j-n^2/4}
   =\frac{e_0(L)^2}{n^2/4-\lambda_0},\\
 \ell\bigl((\mathcal A_{b,\perp}-\xi)^{-1}\ell\bigr)
 &\leq\frac{e_0(L)^2}{n^2/4-\lambda_0}
       \leq\frac{e_0(L)^2}{\xi-\lambda_0}.
 \end{aligned}
\end{equation}
For each $\xi$ in this range, the inverse is defined by the
bilinear form $a_b[\cdot,\cdot]-\xi\ip{\cdot}{\cdot}$ on
$H^1_{\mathrm{sym}}(\Sph^{n-1}_+)\cap e_0^\perp$ and maps
the dual of this subspace into the subspace itself.
\end{lemma}

\begin{proof}
The interior expression in~\eqref{eq:flat-spectral-weighted-expression}
applied to $-\cos\beta$ gives $-(n^2/4)\cos\beta$, since
$(\omega\sin\beta)'=(n-1)\omega\cos\beta$. At $L$, the function
$-\cos\beta$ has value $0$ and derivative $1$. Integration by
parts in~\eqref{eq:flat-spectral-weighted-form} therefore gives
\begin{equation}\label{eq:flat-spectral-green-cos}
 a_b[-\cos\beta,\phi]
 -\frac{n^2}{4}\ip{-\cos\beta}{\phi}
 =[\omega\sin\beta\,\phi]_0^L=\phi(L)
\end{equation}
for every smooth cylindrically symmetric $\phi$. The term at
$0$ vanishes since $\omega\sin\beta=\sin^{n-1}\beta$, and at $L$
we used $\omega(L)=\sin L=1$. By density of smooth cylindrically
symmetric functions and the trace bound
\eqref{eq:flat-spectral-weighted-trace}, this identity holds
for every $\phi\in H^1_{\mathrm{sym}}(\Sph^{n-1}_+)$.

For tests with $\phi(L)=0$, \eqref{eq:flat-spectral-green-cos}
is the weak eigenvalue equation for $-\cos\beta$ on the form
domain $\{f\in H^1_{\mathrm{sym}}(\Sph^{n-1}_+): f(L)=0\}$.
This domain is closed by~\eqref{eq:flat-spectral-weighted-trace}
and embeds compactly into $L^2((0,L),\omega\,\dd\beta)$. Hence
the minimum of $a_b[f]/\ip{f}{f}$ over its nonzero functions is
attained. The minimizer can be chosen nonnegative because taking
the absolute value preserves this quotient. Its inner product
with $\cos\beta$ is positive since $\cos\beta>0$ for $0\leq\beta<L$. Orthogonality of
eigenfunctions for distinct eigenvalues shows that the
Dirichlet minimum is $n^2/4$.

For the unrestricted form $a_b$ on
$H^1_{\mathrm{sym}}(\Sph^{n-1}_+)$, the first eigenvalue
$\lambda_0$ is the minimum of $a_b[f]/\ip{f}{f}$ over $f\ne0$.
Taking $f=1$ gives
\[
 \lambda_0\leq\frac{a_b[1]}{\int_0^L\omega\,\dd\beta}
 =\frac{(n-2)^2}{4}-\frac{b}{\int_0^L\omega\,\dd\beta}
 \leq\frac{(n-2)^2}{4}<\frac{n^2}{4}.
\]
The condition $e_j'(L)=be_j(L)$ and uniqueness for
$\mathcal A_b e_j=\lambda_j e_j$ with
\eqref{eq:flat-spectral-weighted-expression} give
simplicity and $e_j(L)\ne0$. The function
$e_1(L)e_0-e_0(L)e_1$ vanishes at $L$ and has nonzero
components along both $e_0$ and $e_1$. Testing the Dirichlet Rayleigh
quotient with this function gives
\[
 \frac{n^2}{4}
 \leq
 \frac{\lambda_0e_1(L)^2+\lambda_1e_0(L)^2}
      {e_1(L)^2+e_0(L)^2}
 <\lambda_1.
\]
The quotient follows from $\ip{e_i}{e_j}=\delta_{ij}$ and
$a_b[e_i,e_j]=\lambda_i\ip{e_i}{e_j}$. The strict inequality uses
$\lambda_0<\lambda_1$ and $e_1(L)\ne0$.

The compact embedding of $H^1_{\mathrm{sym}}(\Sph^{n-1}_+)$
into $L^2((0,L),\omega\,\dd\beta)$ gives compact resolvent for
$\mathcal A_b$. By the spectral theorem, the squared form norm is
\[
 a_b[f]+(1-\lambda_0)\ip{f}{f}
 =\sum_{j\geq0}(\lambda_j-\lambda_0+1)|\ip{f}{e_j}|^2.
\]
This norm controls $\ip{f}{f}^{1/2}$ because
$\lambda_j-\lambda_0+1\geq1$, and is equivalent to the $H^1$ norm
by the shifted form bound for $a_b$.
For $j\geq1$ and $\lambda_0<\xi\leq n^2/4$,
\[
 0<\frac{\lambda_1-n^2/4}{\lambda_1-\lambda_0+1}
 \leq\frac{\lambda_j-\xi}{\lambda_j-\lambda_0+1}\leq1.
\]
Indeed, the middle ratio is
$1-(\xi-\lambda_0+1)/(\lambda_j-\lambda_0+1)$.
Thus $a_b[\cdot,\cdot]-\xi\ip{\cdot}{\cdot}$ is bounded and
coercive on $H^1_{\mathrm{sym}}(\Sph^{n-1}_+)\cap e_0^\perp$
in the form norm. The Lax--Milgram theorem gives the inverse
in~\eqref{flatbase:robin-resolvent-identity}.

Taking $\phi=e_0$ in~\eqref{eq:flat-spectral-green-cos} gives
$(\lambda_0-n^2/4)\ip{-\cos\beta}{e_0}=e_0(L)$.
On $H^1_{\mathrm{sym}}(\Sph^{n-1}_+)\cap e_0^\perp$ the
component along $e_0$ contributes zero to
$a_b[\cdot,\cdot]-(n^2/4)\ip{\cdot}{\cdot}$. Restricting
\eqref{eq:flat-spectral-green-cos} to test functions in this
subspace gives
\[
 -\cos\beta
 =-\frac{e_0(L)}{n^2/4-\lambda_0}e_0
   +(\mathcal A_{b,\perp}-\tfrac{n^2}{4})^{-1}\ell.
\]
Applying $\ell(f)=f(L)$ and using $\cos L=0$ gives
\[
 \ell\bigl((\mathcal A_{b,\perp}-\tfrac{n^2}{4})^{-1}\ell\bigr)
 =\frac{e_0(L)^2}{n^2/4-\lambda_0}.
\]

By~\eqref{eq:flat-spectral-weighted-trace}, $\ell$ is continuous
in the form norm. Applying $\ell$ to the finite sum with
coefficients $e_j(L)/(\lambda_j-\lambda_0+1)$ in the basis
$(e_j)$ gives its squared form norm. Continuity of $\ell$
bounds that norm independently of the number of terms, so
\[
 \sum_{j\geq0}\frac{e_j(L)^2}{\lambda_j-\lambda_0+1}<\infty.
\]
For $u=(\mathcal A_{b,\perp}-\xi)^{-1}\ell$, taking $\phi=e_j$
in $a_b[u,\phi]-\xi\ip{u}{\phi}=\ell(\phi)$ gives
$(\lambda_j-\xi)\ip{u}{e_j}=e_j(L)$ for $j\geq1$. Hence
\[
 (\mathcal A_{b,\perp}-\xi)^{-1}\ell
 =\sum_{j\geq1}\frac{e_j(L)}{\lambda_j-\xi}e_j.
\]
The lower bound for $(\lambda_j-\xi)/(\lambda_j-\lambda_0+1)$
and the summability of $e_j(L)^2/(\lambda_j-\lambda_0+1)$
give convergence in the form norm. Applying $\ell$ term by term
gives the sum of $e_j(L)^2/(\lambda_j-\xi)$. Each term increases with $\xi$.
At $\xi=n^2/4$ the sum gives the first line of
\eqref{flatbase:robin-resolvent-identity}. Comparison with that
value and $0<\xi-\lambda_0\leq n^2/4-\lambda_0$ give the two
inequalities.
\end{proof}

For $n\geq5$, the unitary map in
\eqref{eq:flat-spectral-form-identification} takes
$H^1_{\mathrm{sym}}(\Sph^{n-1}_+)\cap e_0^\perp$ onto
$\mathcal Q\cap e_b^\perp$ and preserves the trace
$\ell(f)=f(L)$, since $\sin L=1$. By polarization of
\eqref{eq:flat-spectral-form-identification}, the form
$a_b[\cdot,\cdot]-\xi\ip{\cdot}{\cdot}$ corresponds to
$q_b[\cdot,\cdot]-\xi\langle\cdot,\cdot\rangle_{L^2(0,L)}$
on $\mathcal Q\cap e_b^\perp$.
Thus~\eqref{flatbase:robin-resolvent-identity} holds also for
$A_b$ and its restriction to $e_b^\perp$ in $L^2(0,L)$.

\section{Uniqueness in dimensions five and higher}
\label{flatref:section}

We prove the uniqueness statements used in Proposition~\ref{s3:rigidity}
for $n\geq5$. We first derive a differential inequality for the trace
$Y(t,L)$ of a solution $Y$ of~\eqref{flatref:linear-strip}. We then
integrate this inequality over $[0,\infty)$ and use an expansion in
the eigenfunctions for the coefficient $b(0)$ in
\eqref{flatref:linear-strip} to prove uniqueness for odd and
even solutions.

\subsection{An estimate for the trace}
\label{app:weighted-identity}

Use $\beta$, $L$ and $\omega$ as in Section~\ref{subsec:notation}.
Let $b$ be smooth and even on $\R$, with $0<b\leq n/2$ and
$|b'|\leq b$. Let $Y$ be smooth and cylindrically symmetric on
$\R\times\Sph^{n-1}_+$ and satisfy
\begin{equation}\label{flatref:linear-strip}
 \begin{cases}
  \mathcal AY=0 & \text{in }\R\times\Sph^{n-1}_+,\\
  \partial_\beta Y(t,L)=b(t)Y(t,L) & \text{for }t\in\R,
 \end{cases}
\end{equation}
where $\mathcal A$ is defined in Section~\ref{subsec:linearized-spaces}.
Set
\[
 D_Y(t,\beta)=\sin^{\frac{n-2}{2}}\!\beta\,Y(t,\beta).
\]
By~\eqref{eq:flat-spectral-form-identification},
$t\mapsto D_Y(t,\cdot)$ is smooth with values in
$\mathcal Q=\{z\in H^1(0,L):z(0)=0\}$. Also $D_Y(t,L)=Y(t,L)$
since $\sin L=1$. Let $\phi$ be smooth and cylindrically symmetric on $\Sph^{n-1}_+$.
By the interior equation in~\eqref{flatref:linear-strip},
integration by parts in~\eqref{eq:flat-spectral-weighted-form} gives
\[
 a_{b(t)}[Y(t,\cdot),\phi]
 =\int_0^L\omega\,\partial_t^2Y(t,\beta)\phi(\beta)\,\dd\beta.
\]
Indeed, $\omega Y_\beta\phi$ vanishes at $\beta=0$ by smoothness
and $\omega(0)=0$. At $\beta=L$, it equals $b(t)Y(t,L)\phi(L)$
by the boundary condition in~\eqref{flatref:linear-strip} and
$\omega(L)=1$.
By density and the bilinear identity associated
with~\eqref{eq:flat-spectral-form-identification},
\begin{equation}\label{flatref:weak-difference}
 \langle (D_Y)_{tt},z\rangle=q_{b(t)}[D_Y,z]
 \qquad\text{for every }z\in\mathcal Q.
\end{equation}
We write $\langle\cdot,\cdot\rangle$ and $\|\cdot\|_2$ for the
inner product and norm in $L^2(0,L)$.

For each $t\in\R$, let $\mu=\mu_0(b(t))$ and $e=e_{b(t)}$
be the first eigenvalue of $A_{b(t)}$ and its positive normalized
eigenfunction, as in Proposition~\ref{flatbase:angular}.
In this subsection, $q_b$, $C$ and $r$ depend on $t$ through $b(t)$. Write
\[
 D_Y=ue+H,\qquad (D_Y)_t=we+J,
\]
where $H,J\in\mathcal Q$ and $H,J\perp e$ in $L^2(0,L)$.
Since $\langle e_t,e\rangle=0$ by normalization, differentiating
the projections and using~\eqref{flatref:weak-difference} gives
\begin{equation}\label{flatref:projections}
 \begin{aligned}
 u'&=w+\langle H,e_t\rangle,\\
 w'&=\mu u+\langle J,e_t\rangle,\\
 \frac{\dd}{\dd t}\langle H,J\rangle
 &=q_b[H]+\|J\|_2^2-u\langle e_t,J\rangle-w\langle e_t,H\rangle.
 \end{aligned}
\end{equation}
The last identity follows by pairing
$H_t=J+(w-u')e-u e_t$ with $J$ and
$J_t=(D_Y)_{tt}-w'e-w e_t$ with $H$.

We use the weight $e^{-(n-1)t}$ and set
\[
 \begin{gathered}
 F=\langle H,J\rangle-\frac{n-1}{4}u^2-\frac{w^2}{n-1},\\
 x_\pm=\frac{n-1}{2}u\pm w,\qquad
 P=\frac14-\frac{\mu}{(n-1)^2},\qquad Q=1-P,\\
 \mathcal E=q_b[H]+\|J\|_2^2-(n-1)\langle H,J\rangle .
 \end{gathered}
\]
By~\eqref{flatbase:angular-bounds} and $b\leq n/2$,
\[
 \begin{aligned}
 (n-1)^2P
 &=r+\frac{2n-3}{4}>0,\\
 (n-1)^2Q
 &=\frac{3(n-1)^2}{4}+\mu\\
 &\geq\frac{3(n-1)^2-6(n-1)+2}{4}>0.
 \end{aligned}
\]
Since $H\perp e$, Proposition~\ref{flatbase:angular} gives
$q_b[H]\geq\mu_1(b)\|H\|_2^2$ with $\mu_1(b)>n^2/4$.
In particular, $\mathcal E\geq0$ by the identity
\begin{equation}\label{flatref:angular-energy}
 \mathcal E=q_b[H]-\frac{(n-1)^2}{4}\|H\|_2^2
       +\left\|J-\frac{n-1}{2}H\right\|_2^2.
\end{equation}
By~\eqref{flatref:projections},
\begin{equation}\label{flatref:exact-dissipation}
 F'-(n-1)F
 =P x_+^2+Q x_-^2+\mathcal E
       -x_+\left\langle e_t,H+\frac{2J}{n-1}\right\rangle.
\end{equation}

Fix $t\in\R$ and let $(e_j)_{j\geq0}$ be an orthonormal basis of
eigenfunctions of $A_{b(t)}$, with $e_0=e$. Let $\mu_j=\mu_j(b(t))$
and $\delta=n^2/4-\mu=r+n-1$. The smooth dependence in
Proposition~\ref{flatbase:angular} allows us to differentiate
the weak eigenvalue equation for $e$, holding $e_j$ fixed. For $j\geq1$, using
$\partial_b q_b[z,\phi]=-z(L)\phi(L)$
from~\eqref{eq:flat-spectral-unweighted-form} gives
\begin{equation}\label{flatref:ground-derivative}
 \langle e_t,e_j\rangle
   =\frac{b'\sqrt C\,e_j(L)}{\mu_j-\mu}.
\end{equation}
Write $H+2J/(n-1)=2H+2(J-(n-1)H/2)/(n-1)$. Applying
Cauchy--Schwarz in the eigenfunction expansion with the two
terms of~\eqref{flatref:angular-energy} gives
\[
 \left\langle e_t,H+\frac{2J}{n-1}\right\rangle^2
 \leq\frac{4(b')^2C}{(n-1)^2}\,\mathcal E
 \sum_{j\geq1}
 \frac{e_j(L)^2}{(\mu_j-\mu)^2}
 \frac{\mu_j+3(n-1)^2/4}{\mu_j-(n-1)^2/4}.
\]
Since $Q>0$, the ratio
$(\mu_j+3(n-1)^2/4)/(\mu_j-\mu)$ is at most its value with
$\mu_j$ replaced by $n^2/4$. Also $(2n-1)/(n-1)^2\leq9/16$.
Thus
\[
 \frac{4}{(n-1)^2}\frac{\mu_j+3(n-1)^2/4}{\mu_j-\mu}
 \leq\frac{73}{16\delta}.
\]
For $z=\mu_j-n^2/4>0$, the arithmetic--geometric mean inequality
gives
\[
 \mu_j-\mu\geq2\sqrt{z\delta},\qquad
 \mu_j-\frac{(n-1)^2}{4}\geq\sqrt{z(2n-1)}.
\]
By~\eqref{eq:flat-spectral-form-identification}, the resolvent
identity~\eqref{flatbase:robin-resolvent-identity} applies to
$(e_j)_{j\geq0}$. It follows that
\[
 \left\langle e_t,H+\frac{2J}{n-1}\right\rangle^2
 \leq\frac{73(b'C)^2}{32\delta^{5/2}\sqrt{2n-1}}\,\mathcal E .
\]
To compare the coefficient with $P$, use $|b'|C\leq bC\leq2r$
by~\eqref{flatbase:angular-gap}. With $s=r/(n-1)$, the ratio of
the coefficient to $P$ is at most
\[
 \frac{73\sqrt2}{16}
 \frac{s^2}{(1+s)^{5/2}(s+7/16)}.
\]
Here we used $2n-1\geq2(n-1)$ and
$(2n-3)/4\geq7(n-1)/16$. Since $s+7/16\geq(\sqrt7/2)\sqrt s$
and $s^{3/2}/(1+s)^{5/2}$ is largest at $s=3/2$, the square
of this upper bound is at most $143883/175000<1$.
Applying this estimate in~\eqref{flatref:exact-dissipation} gives
\begin{equation}\label{flatref:coercive}
 F'-(n-1)F\geq \frac12 P x_+^2+Q x_-^2+\frac12\mathcal E.
\end{equation}

By Cauchy--Schwarz in the eigenfunction expansion of $H$ and
\eqref{flatbase:robin-resolvent-identity},
\[
 \frac{\delta}{C}|H(t,L)|^2+\frac{n^2}{4}\|H\|_2^2\leq q_b[H].
\]
By~\eqref{flatref:angular-energy}, $|H(t,L)|^2\leq C\mathcal E/\delta$.
Write $Y(t,L)=\sqrt C\,(x_++x_-)/(n-1)+H(t,L)$ and apply
Cauchy--Schwarz using the three terms on the right-hand side
of~\eqref{flatref:coercive}. The coefficient of $F'-(n-1)F$
in the bound for $b|Y(t,L)|^2/2$ is at most
\[
 \begin{aligned}
 &\frac{bC}{(n-1)^2P}+\frac{bC}{2(n-1)^2Q}+\frac{bC}{\delta}\\
 &\qquad\leq
 \frac{2r}{r+(2n-3)/4}+\frac{2r}{r+n-1}
       +\frac{r}{3(n-1)^2/4+(n-2)^2/4-r}.
 \end{aligned}
\]
Here we used $bC\leq2r$ from~\eqref{flatbase:angular-gap}.
Since $r\leq(n^2+2n-4)/4$ and $n-1\geq4$, we have
\[
 \begin{aligned}
 2-\frac{2r}{r+(2n-3)/4}
 &\geq\frac{2(2n-3)}{(n-1)^2+6(n-1)-2}
 \geq\frac{7}{5(n-1)},\\
 2-\frac{2r}{r+n-1}
 &\geq\frac{8(n-1)}{(n-1)^2+8(n-1)-1}
 \geq\frac{8}{3(n-1)}.
 \end{aligned}
\]
The third fraction is at most
\[
 \frac{(n-1)^2+4(n-1)-1}{3(n-1)^2-6(n-1)+2}
 \leq\frac13+\frac4{n-1},
\]
since $3(n-1)^2-6(n-1)+2\geq3(n-1)^2/2$. The three fractions
therefore sum to less than $13/3$. This proves
\begin{equation}\label{flatref:weighted-trace}
 b|Y(t,L)|^2\leq\frac{26}{3}\bigl(F'-(n-1)F\bigr).
\end{equation}

\subsection{Uniqueness for odd and even solutions}
\label{app:even-jacobi}

\begin{lemma}\label{flatref:parity-uniqueness}
\label{flatref:odd-endpoint}\label{flatneck:even}
Let $n\geq5$ and let $b$ be smooth and even on $\R$, with
$0<b\leq n/2$ and $|b'|\leq b$. Let $Y$ be a smooth
cylindrically symmetric solution of~\eqref{flatref:linear-strip}
on $\R\times\Sph^{n-1}_+$. Suppose that $Y$ and $\partial_tY$
have polynomial growth in the sense of~\eqref{eq:jacobi-polynomial-space}.
Then $Y\equiv0$ in either of the following cases.
\begin{enumerate}[label=\textup{(\roman*)}]
 \item $Y$ is odd in $t$, $\partial_tY(0,L)=0$, and
 $|b'(t)|/b(t)\leq(n-2)t/2$ for $t\geq0$.
 \item $Y$ is even in $t$ and $Y(0,L)=0$.
\end{enumerate}
\end{lemma}

\begin{proof}
The growth assumption gives $e^{-(n-1)t}F(t)\to0$ as $t\to\infty$.
If $Y$ is odd, then $D_Y(0,\cdot)=0$ and $F(0)=-w(0)^2/(n-1)$.
If $Y$ is even, then $(D_Y)_t(0,\cdot)=0$ and $F(0)=-(n-1)u(0)^2/4$.
Integrating~\eqref{flatref:weighted-trace} with weight $e^{-(n-1)t}$
gives
\begin{equation}\label{flatref:integrated-trace}
 I:=\int_0^\infty e^{-(n-1)t}b(t)|Y(t,L)|^2\,\dd t
 \leq
 \begin{cases}
  \dfrac{26}{3(n-1)}w(0)^2 & \text{if }Y\text{ is odd},\\[4pt]
  \dfrac{13(n-1)}6 u(0)^2 & \text{if }Y\text{ is even}.
 \end{cases}
\end{equation}

For the rest of the proof, let $\mu_j=\mu_j(b(0))$ and fix
an orthonormal basis $(e_j)_{j\geq0}$ of eigenfunctions of
$A_{b(0)}$, with $e_0=e_{b(0)}$. We evaluate $C$, $r$ and
$\delta=n^2/4-\mu_0$ at $b(0)$. Let $V=(D_Y)_t(0,\cdot)$
if $Y$ is odd and $V=D_Y(0,\cdot)$ if $Y$ is even. Set
$V_j=\langle V,e_j\rangle$. Thus $V_0=w(0)$ in the odd case
and $V_0=u(0)$ in the even case.

In both cases, $V(L)=0$ by the hypothesis at $(0,L)$.
The variational characterization of $n^2/4$ on
$\{z\in\mathcal Q:z(L)=0\}$ in Proposition~\ref{flatbase:angular}
gives $q_{b(0)}[V]\geq(n^2/4)\|V\|_2^2$. Expanding in the fixed
basis gives
\begin{equation}\label{flatref:initial}
 \delta |V_0|^2\leq
 \sum_{j\geq1}(\mu_j-n^2/4)|V_j|^2.
\end{equation}

For $j\geq1$, set $\kappa_j=\sqrt{\mu_j}>n/2$ and
$Y_j(t)=\langle D_Y(t,\cdot),e_j\rangle$. By
\eqref{flatref:weak-difference} with $z=e_j$,
\begin{equation}\label{flatref:fixed-equation}
 -Y_j''+\kappa_j^2Y_j=e_j(L)\bigl(b(t)-b(0)\bigr)Y(t,L).
\end{equation}
The right-hand side is $q_{b(0)}[D_Y,e_j]-q_{b(t)}[D_Y,e_j]$
by~\eqref{eq:flat-spectral-unweighted-form}. By parity, $Y_j(0)=0$
if $Y$ is odd and $Y_j'(0)=0$ if $Y$ is even.
Multiply~\eqref{flatref:fixed-equation} by $e^{-\kappa_jt}$ and
integrate by parts on $[0,\infty)$. The polynomial growth of $Y$
and $\partial_tY$ gives
\begin{equation}\label{flatref:fixed-datum}
 Y_j'(0)+\kappa_jY_j(0)=
 e_j(L)\int_0^\infty e^{-\kappa_jt}\bigl(b(t)-b(0)\bigr)Y(t,L)\,\dd t.
\end{equation}
The left-hand side is $V_j$ if $Y$ is odd and $\kappa_jV_j$ if
$Y$ is even.

Integrating $|b'|/b\leq1$ on $[0,t]$ gives
$|\log(b(t)/b(0))|\leq t$. Since $\sinh^2$ is even and
increasing on $[0,\infty)$, we have
\begin{equation}\label{flatref:coefficient-variation}
 \frac{|b(t)-b(0)|^2}{b(t)}
 =4b(0)\sinh^2\frac{\log b(t)-\log b(0)}{2}
 \leq b(0)e^t(1-e^{-t})^2.
\end{equation}
In case~\textup{(i)}, integrating $|b'(t)|/b(t)\leq(n-2)t/2$
gives $|\log(b(t)/b(0))|\leq(n-2)t^2/4$. By convexity,
$\sinh(\lambda t/2)\leq\lambda\sinh(t/2)$ for $0\leq\lambda\leq1$.
Thus in the odd case the right-hand side
of~\eqref{flatref:coefficient-variation} can be multiplied by
$\min\{1,(n-2)^2t^2/16\}$.

Applying Cauchy--Schwarz in~\eqref{flatref:fixed-datum} with the
integral $I$ in~\eqref{flatref:integrated-trace} gives
\[
 |Y_j'(0)+\kappa_jY_j(0)|^2
 \leq e_j(L)^2 I\int_0^\infty e^{-(2\kappa_j-n+1)t}
                    \frac{|b(t)-b(0)|^2}{b(t)}\,\dd t .
\]
Put $x=\kappa_j-n/2>0$, so that $\mu_j-n^2/4=x(x+n)$.
Using~\eqref{flatref:coefficient-variation}, with the additional
factor $\min\{1,(n-2)^2t^2/16\}$ in the odd case, reduces the
estimates for odd and even $Y$ to the first and second inequalities
below, respectively.
\begin{equation}\label{flatref:frozen-moments}
 \begin{aligned}
 x^2(x+n)^2\int_0^\infty e^{-2xt}(1-e^{-t})^2
       \min\{1,(n-2)^2t^2/16\}\,\dd t&\leq \frac{15n^2}{112},\\
 \frac{x^2(x+n)^2}{(x+n/2)^2}
        \int_0^\infty e^{-2xt}(1-e^{-t})^2\,\dd t&\leq\frac38.
 \end{aligned}
\end{equation}
Expanding the square in the integral without the minimum gives
\[
 \int_0^\infty e^{-2xt}(1-e^{-t})^2\,\dd t
 =\frac{1}{2x(x+1)(2x+1)}.
\]
For $0<x\leq5/2$, we omit the minimum, divide the bound for the
first line of~\eqref{flatref:frozen-moments} by $n^2$, and use
$n\geq5$. This gives at most
$\frac{x(x+5)^2}{50(x+1)(2x+1)}$, which is at most $15/112$ since
\[
 750(x+1)(2x+1)-112x(x+5)^2
   =(5-2x)(56x^2-50x+150)\geq0.
\]
For $5/2\leq x\leq n/2$, use $1+x/n\leq3/2$ and the decrease
of $\frac{x}{(x+1)(2x+1)}$ on $[5/2,\infty)$ to obtain the same bound.
For $x\geq n/2$, use $1-e^{-t}\leq t$ and integrate
$t^4e^{-2xt}$. The left-hand side of the first line, divided by $n^2$,
is then at most
\[
 \frac{3(n-2)^2(x+n)^2}{64n^2x^3}
 \leq\frac{27(n-2)^2}{32n^3}\leq\frac1{16}.
\]
The first inequality uses $x+n\leq3x$ and $x\geq n/2$.
The derivative of $(n-2)^2/n^3$ is $(n-2)(6-n)/n^4$,
so its maximum on $[5,\infty)$ occurs at $n=6$.
For the second line of~\eqref{flatref:frozen-moments}, we use
$(x+n)/(x+n/2)\leq2$ and
$\frac{2x}{(x+1)(2x+1)}\leq3/8$, which follows from $6x^2-7x+3>0$.

The inequalities~\eqref{flatref:frozen-moments} bound the sum
in~\eqref{flatref:initial} by
\[
 \begin{cases}
 \dfrac{15n^2}{112}\dfrac{b(0)C}{\delta} I
        &\text{if }Y\text{ is odd},\\[5pt]
 \dfrac38\dfrac{b(0)C}{\delta} I
        &\text{if }Y\text{ is even}.
 \end{cases}
\]
Indeed, in both cases the spectral sum is
$\sum_{j\geq1}e_j(L)^2/(\mu_j-n^2/4)=C/\delta$ by
\eqref{flatbase:robin-resolvent-identity}. Using
$b(0)C\leq2r$ from~\eqref{flatbase:angular-gap} gives
\[
 \frac{b(0)C}{\delta^2}
 \leq\frac{2r}{(r+n-1)^2}\leq\frac{1}{2(n-1)}.
\]
Consequently, by~\eqref{flatref:integrated-trace},
\[
 \sum_{j\geq1}(\mu_j-n^2/4)|V_j|^2
 \leq
 \begin{cases}
  \dfrac{1625}{1792}\delta |V_0|^2 &\text{if }Y\text{ is odd},\\[4pt]
  \dfrac{13}{32}\delta |V_0|^2 &\text{if }Y\text{ is even}.
 \end{cases}
\]
In the odd case we used $n/(n-1)\leq5/4$. Both constants are
less than $1$, so~\eqref{flatref:initial} gives $V_0=0$ and
hence $F(0)=0$. Integrating~\eqref{flatref:coercive} with weight
$e^{-(n-1)t}$ gives $x_+=x_-=0$ and $\mathcal E=0$ on $[0,\infty)$.
Thus $u=w=0$, and~\eqref{flatref:angular-energy} gives $H=J=0$
since $\mu_1(b(t))>n^2/4$. By parity, $Y\equiv0$ on
$\R\times\Sph^{n-1}_+$.
\end{proof}

The growth assumption in Lemma~\ref{flatref:parity-uniqueness} can be
replaced by
\[
 \liminf_{t\to\infty}e^{-(n-1)t/2}
 \bigl(\|D_Y(t,\cdot)\|_2+\|(D_Y)_t(t,\cdot)\|_2\bigr)=0 .
\]
Choose $t_i\to\infty$ along which the expression in the limit
tends to zero. Then $e^{-(n-1)t_i}F(t_i)\to0$, and the
nonnegativity in~\eqref{flatref:coercive} gives
\eqref{flatref:integrated-trace}. Since $\kappa_j>(n-1)/2$,
the boundary term in the integration leading
to~\eqref{flatref:fixed-datum} vanishes along the same sequence.
The Cauchy--Schwarz estimate using~\eqref{flatref:integrated-trace}
and~\eqref{flatref:frozen-moments} gives absolute convergence of the
integral in~\eqref{flatref:fixed-datum}.

\section{Uniqueness in dimensions three and four}
\label{nfour:section}

We prove the uniqueness statements used in Proposition~\ref{s3:rigidity}
for $n=3,4$. We expand $Y(t,\cdot)$ in the eigenfunctions of
$\mathcal A_{n/4}$, defined by~\eqref{eq:flat-spectral-weighted-form}.
For each parity, the condition at $(0,L)$ determines an integral
equation for $Y(t,L)$. We estimate its kernel separately in
dimensions three and four.

\begin{theorem}\label{nfour:linear-uniqueness}
Let $n\in\{3,4\}$ and let $b$ be smooth and even on $\R$, with
$0\leq b\leq n/2$. Let $Y$ be smooth and cylindrically symmetric on
$\R\times\Sph^{n-1}_+$, with polynomial growth in the sense
of~\eqref{eq:jacobi-polynomial-space}. Suppose that
\begin{equation}\label{nfour:linear-equation}
 \begin{cases}
  \mathcal AY=0 & \text{in }\R\times\Sph^{n-1}_+,\\
  \partial_\beta Y(t,L)=b(t)Y(t,L) & \text{for }t\in\R .
 \end{cases}
\end{equation}
Here $\mathcal A$ is as in Section~\ref{subsec:linearized-spaces}.
Then $Y\equiv0$ on $\R\times\Sph^{n-1}_+$ in either of the following cases.
\begin{enumerate}[label=\textup{(\roman*)}]
 \item $Y$ is odd in $t$ and $\partial_tY(0,L)=0$.
 \item $Y$ is even in $t$ and $Y(0,L)=0$.
\end{enumerate}
\end{theorem}

\begin{proof}
Use $\beta$, $L$ and $\omega$ as in Section~\ref{subsec:notation}.
Set $b_*=n/4$. Let $\mathcal A_{b_*}$ be the operator associated
with $a_{b_*}$ in~\eqref{eq:flat-spectral-weighted-form}, with form
domain $H^1_{\mathrm{sym}}(\Sph^{n-1}_+)$. The inner product
$\ip{\cdot}{\cdot}$ is taken in $L^2((0,L),\omega\,\dd\beta)$.

Testing~\eqref{eq:flat-spectral-weighted-form} with the constant
function $1$ gives
\[
 a_{b_*}[1]=
 \begin{cases}
  \dfrac14\displaystyle\int_0^L\sin\beta\,\dd\beta-\dfrac34
      =-\dfrac12 & (n=3),\\[3pt]
  \displaystyle\int_0^L\sin^2\beta\,\dd\beta-1
      =\dfrac L2-1<0 & (n=4).
 \end{cases}
\]

Let $\lambda_j$ and $e_j$, $j\geq0$, be the eigenvalues and real
orthonormal eigenfunctions of $\mathcal A_{b_*}$. By the Rayleigh
quotient and Lemma~\ref{flatbase:robin-resolvent},
$\lambda_0<0$ and $\lambda_1>n^2/4$. We choose $e_0$ positive
on $(0,L]$ and the signs of the other eigenfunctions so that
$e_j(L)>0$. The absolute value of a Rayleigh minimizer is also a minimizer.
If this nonnegative eigenfunction vanishes at some $\beta\in(0,L)$,
then its derivative vanishes there and uniqueness for the
eigenfunction equation gives $e_0=0$. At $\beta=L$, the
condition $e_j'(L)=b_*e_j(L)$ and the same uniqueness argument
give $e_j(L)\ne0$ for every $j\geq0$.
Write $\kappa_0=\sqrt{-\lambda_0}$, $k_j=\sqrt{\lambda_j}$ for
$j\geq1$, and $C_j=e_j(L)^2$. By~\eqref{flatbase:robin-resolvent-identity},
\begin{equation}\label{nfour:sum-residues}
 A_*:=\sum_{j\geq1}\frac{C_j}{k_j^2-n^2/4}
      =\frac{C_0}{n^2/4+\kappa_0^2}.
\end{equation}

Let $F(t)=(b(t)-b_*)Y(t,L)$ and $Y_j(t)=\ip{Y(t,\cdot)}{e_j}$.
Testing~\eqref{nfour:linear-equation} against $e_j$ and integrating
by parts in $\beta$ gives
\begin{equation}\label{nfour:modal-equation}
 -Y_j''+\lambda_jY_j=e_j(L)F .
\end{equation}
At $\beta=0$ the boundary term vanishes by smoothness and $\omega(0)=0$.
At $\beta=L$ the difference of the boundary terms is
$F(t)e_j(L)$ by~\eqref{nfour:linear-equation} and
$e_j'(L)=b_*e_j(L)$. For $j\geq1$, the conditions at $t=0$ are
$Y_j(0)=0$ in the odd case and $Y_j'(0)=0$ in the even case.
Thus
\begin{equation}\label{nfour:high-representation}
 \begin{aligned}
 Y_j(t)&=e_j(L)\int_0^\infty G_{k_j}^{D/N}(t,\tau)F(\tau)\,\dd\tau,\\
 G_{k_j}^{D/N}(t,\tau)&=
 \frac{e^{-k_j|t-\tau|}\mp e^{-k_j(t+\tau)}}{2k_j}.
 \end{aligned}
\end{equation}
The minus sign gives the Dirichlet condition and the plus sign
gives the Neumann condition at $t=0$. The jump of $\partial_tG_{k_j}^{D/N}$
at $t=\tau$ is $-1$, so
$(-\partial_t^2+k_j^2)G_{k_j}^{D/N}(t,\tau)=\delta_\tau$
in distributions on $(0,\infty)$, where $\delta_\tau$ is the Dirac
measure at $\tau$. Both kernels are nonnegative
since $|t-\tau|\leq t+\tau$ for $t,\tau\geq0$.
The integral in~\eqref{nfour:high-representation} has polynomial
growth, since $F$ has polynomial growth and the kernel decays
exponentially in $|t-\tau|$. Applying $-\partial_t^2+k_j^2$ to
the difference between $Y_j$ and the first line of
\eqref{nfour:high-representation} gives zero. Polynomial growth
excludes the term in $e^{k_jt}$, and the condition at $0$ excludes
the term in $e^{-k_jt}$.

For $\lambda>n/2$, use $G_\lambda^D$ and $G_\lambda^N$ from
\eqref{nfour:high-representation} with $k_j$ replaced by $\lambda$.
Then
\begin{equation}\label{nfour:green-moments}
 \begin{aligned}
 \int_0^\infty G_\lambda^D(t,\tau)\sinh(n\tau/2)\,\dd\tau
 &=\frac{\sinh(nt/2)}{\lambda^2-n^2/4},\\
 \int_0^\infty G_\lambda^N(t,\tau)\cosh(n\tau/2)\,\dd\tau
 &=\frac{\cosh(nt/2)}{\lambda^2-n^2/4}.
 \end{aligned}
\end{equation}
Applying $-\partial_t^2+\lambda^2$ to both sides gives
$\sinh(nt/2)$ in the first line and $\cosh(nt/2)$ in the second.
Both sides of the first line vanish at $t=0$, and both sides
of the second have zero $t$-derivative there. The growth bound
$O(e^{nt/2})$
and $\lambda>n/2$ give uniqueness. Let $h(t)=\sinh(nt/2)$
if $Y$ is odd and $h(t)=\cosh(nt/2)$ if $Y$ is even.
Smoothness, parity and polynomial growth give $|F|\leq Mh$
on $[0,\infty)$ for some $M>0$. By
\eqref{nfour:green-moments}, the $j$th term $e_j(L)Y_j(t)$
for $j\geq1$ is bounded in absolute value by
$Mh(t)C_j/(k_j^2-n^2/4)$. The series converges locally uniformly
in $t$ by~\eqref{nfour:sum-residues}. Its sum, including the term
$j=0$, is $Y(t,L)$ because the expansion of $Y(t,\cdot)$ converges
in the shifted form norm from the proof of
Lemma~\ref{flatbase:robin-resolvent}. The map $f\mapsto f(L)$
is continuous in that norm by~\eqref{eq:flat-spectral-weighted-trace}.

For $j\geq1$, multiply~\eqref{nfour:modal-equation} by
$e^{-k_jt}$ and integrate by parts on $[0,\infty)$ to obtain
\[
 Y_j'(0)+k_jY_j(0)
 =e_j(L)\int_0^\infty e^{-k_j\tau}F(\tau)\,\dd\tau .
\]
The term at infinity vanishes since $Y_j$ and $F$ have polynomial
growth and~\eqref{nfour:modal-equation} gives polynomial growth of
$Y_j''$ and hence of $Y_j'$. Expand $\partial_tY(0,\cdot)$ in the odd case and $Y(0,\cdot)$
in the even case in the same shifted form norm. Taking the trace
at $\beta=L$ and using parity gives
\begin{equation}\label{nfour:critical-traces}
 \begin{aligned}
 \partial_tY(0,L)
 &=e_0(L)Y_0'(0)+\sum_{j\geq1}C_j
   \int_0^\infty e^{-k_j\tau}F(\tau)\,\dd\tau
   &&\text{if }Y\text{ is odd},\\
 Y(0,L)
 &=e_0(L)Y_0(0)+\sum_{j\geq1}\frac{C_j}{k_j}
   \int_0^\infty e^{-k_j\tau}F(\tau)\,\dd\tau
   &&\text{if }Y\text{ is even}.
 \end{aligned}
\end{equation}

Let $S(t)=\sin(\kappa_0t)/\kappa_0$ and $P=\pi/\kappa_0$.
By~\eqref{nfour:modal-equation}, $Y_0''+\kappa_0^2Y_0=-e_0(L)F$.
Its solution has homogeneous term $Y_0'(0)S(t)$ in the odd case
and $Y_0(0)\cos(\kappa_0t)$ in the even case. In either case
the term with forcing is $-e_0(L)\int_0^t S(t-\tau)F(\tau)\,\dd\tau$.
The hypotheses at $(0,L)$ and~\eqref{nfour:critical-traces}
determine the homogeneous term. Using it in the expansion of
$Y(t,L)$ gives, for $t\geq0$,
\begin{equation}\label{nfour:trace-representation}
 \begin{aligned}
 Y(t,L)&=\int_0^\infty K^D(t,\tau)F(\tau)\,\dd\tau
   &&\text{if }Y\text{ is odd},\\
 Y(t,L)&=\int_0^\infty K^N(t,\tau)F(\tau)\,\dd\tau
   &&\text{if }Y\text{ is even},
 \end{aligned}
\end{equation}
where
\begin{align}
 K^D(t,\tau)&=-C_0S(t-\tau)\mathbf1_{\{\tau<t\}}\nonumber\\
 &\quad+\sum_{j\geq1} C_j
     \bigl(G_{k_j}^D(t,\tau)-S(t)e^{-k_j\tau}\bigr),
      \label{nfour:odd-kernel}\\
 K^N(t,\tau)&=-C_0S(t-\tau)\mathbf1_{\{\tau<t\}}\nonumber\\
 &\quad+\sum_{j\geq1} C_j
     \left(G_{k_j}^N(t,\tau)-\frac{\cos(\kappa_0t)}{k_j}e^{-k_j\tau}\right).
      \label{nfour:even-kernel}
\end{align}
For each parity, the corresponding series converges in
$L^1((0,\infty),h(\tau)\,\dd\tau)$
by~\eqref{nfour:sum-residues}, \eqref{nfour:green-moments} and
\[
 \begin{aligned}
 \int_0^\infty e^{-k_j\tau}\sinh(n\tau/2)\,\dd\tau
   &=\frac{n/2}{k_j^2-n^2/4},\\
 \int_0^\infty e^{-k_j\tau}\cosh(n\tau/2)\,\dd\tau
   &=\frac{k_j}{k_j^2-n^2/4}.
 \end{aligned}
\]
Since $|F|\leq Mh$, these estimates also justify integration
of the series in~\eqref{nfour:trace-representation}.

\medskip
\noindent\emph{Bounds for the kernels.}\enspace
Define
\[
 \begin{aligned}
 Q^D(t)&=\frac{\displaystyle\int_0^\infty |K^D(t,\tau)|
                         \sinh(n\tau/2)\,\dd\tau}{\sinh(nt/2)}
                         &&(t>0),\\
 Q^N(t)&=\frac{\displaystyle\int_0^\infty |K^N(t,\tau)|
                         \cosh(n\tau/2)\,\dd\tau}{\cosh(nt/2)}
                         &&(t\geq0).
 \end{aligned}
\]
We bound $b_*Q^D$ on $(0,\infty)$ and $b_*Q^N$ on $[0,\infty)$.
By~\eqref{nfour:sum-residues} and~\eqref{nfour:green-moments},
the series in
\eqref{nfour:odd-kernel} and~\eqref{nfour:even-kernel}
contribute at most $A_*\sinh(nt/2)+(n/2)A_*|S(t)|$ and
$A_*\cosh(nt/2)+A_*|\cos(\kappa_0t)|$ to the respective integrals.
For the remaining term, set
\[
 \begin{aligned}
 I^D(t)&=\int_0^t|S(t-\tau)|\sinh(n\tau/2)\,\dd\tau,\\
 I^N(t)&=\int_0^t|S(t-\tau)|\cosh(n\tau/2)\,\dd\tau .
 \end{aligned}
\]
Since $S\geq0$ on $[0,P]$, applying $\partial_t^2+\kappa_0^2$
to $I^D$ and $I^N$ gives $\sinh(nt/2)$ and $\cosh(nt/2)$,
respectively, on this interval. Both integrals and their first
derivatives vanish at $0$. Hence on $[0,P]$,
\begin{equation}\label{nfour:low-first}
 \begin{aligned}
 (n^2/4+\kappa_0^2)I^D(t)&=\sinh(nt/2)-\frac n2 S(t),\\
 (n^2/4+\kappa_0^2)I^N(t)&=\cosh(nt/2)-\cos(\kappa_0t).
 \end{aligned}
\end{equation}

For $t=mP+\delta$, where $m\geq0$ is an integer and
$0\leq\delta<P$, the identity
\[
 |S(r)|=S(r)+2\sum_{1\leq j\leq r/P} S(r-jP),\qquad r\geq0,
\]
follows from $S(r-jP)=(-1)^jS(r)$. The convolutions with $S$ in place of $|S|$ satisfy
\eqref{nfour:low-first} on $[0,\infty)$. Applying the identity
for $|S|$ then gives
\[
 \begin{aligned}
 (n^2/4+\kappa_0^2)I^D(t)
 &=\sinh(nt/2)-\frac n2|S(t)|+2\sum_{j=1}^{m}\sinh(n(t-jP)/2),\\
 (n^2/4+\kappa_0^2)I^N(t)
 &=\cosh(nt/2)-\cos(\kappa_0\delta)+2\sum_{j=1}^{m}\cosh(n(t-jP)/2).
 \end{aligned}
\]
For $1\leq j\leq m$,
\[
 \begin{aligned}
 \frac{\sinh(n(t-jP)/2)}{\sinh(nt/2)}&\leq e^{-njP/2},\\
 \frac{\cosh(n(t-jP)/2)}{\cosh(nt/2)}&\leq2e^{-njP/2}.
 \end{aligned}
\]
Summing the geometric series gives, for $t\geq P$,
\begin{equation}\label{nfour:low-global}
 \begin{aligned}
 (n^2/4+\kappa_0^2)I^D(t)
 &\leq \sinh(nt/2)+\frac{2\sinh(nt/2)}{e^{nP/2}-1}-\frac n2|S(t)|,\\
 (n^2/4+\kappa_0^2)I^N(t)
 &\leq \cosh(nt/2)+\frac{4\cosh(nt/2)}{e^{nP/2}-1}+1 .
 \end{aligned}
\end{equation}
Combining these bounds with~\eqref{nfour:low-first} and the
estimates of the two series gives
\begin{align}
 \frac{Q^D(t)}{2A_*}&\leq1
   &&(0<t\leq P),\notag\\
 \frac{Q^D(t)}{2A_*}&\leq1+\frac1{e^{nP/2}-1}
   &&(t\geq P),\notag\displaybreak[3]\\
 \frac{Q^N(t)}{2A_*}&\leq1
   &&(0\leq t\leq P/2),\notag\\
 \frac{Q^N(t)}{2A_*}&\leq1+\frac{4\kappa_0}{n}e^{-1-nP/4}
   &&(P/2\leq t\leq P),\notag\\
 \frac{Q^N(t)}{2A_*}&\leq1+\frac2{e^{nP/2}-1}+\frac1{\cosh(nP/2)}
   &&(t\geq P).\label{nfour:preliminary-table}
\end{align}
In the fourth line, use $|\cos(\kappa_0t)|\leq\kappa_0(t-P/2)$,
$1/\cosh(nt/2)\leq2e^{-nt/2}$, and
$\sup_{s\geq0}se^{-ns/2}=2/(ne)$.

\medskip
\noindent\emph{Dimension three.}\enspace
Here $b_*=3/4$. In the variable $x=\cos\beta$ the expression
\eqref{eq:flat-spectral-weighted-expression} is
$T=-\partial_x((1-x^2)\partial_x)+1/4$ on $(0,1)$.
Consider $f(x)=e^{-9x/32}(1+x)^{-15/32}$. Since
$f'(0)=-3f(0)/4$, the function $f(\cos\beta)$ belongs to
$\operatorname{Dom}(\mathcal A_{b_*})$. Differentiation gives
\[
 \frac{(T+25/32)f}{f}
 =\frac{x(1-x)}{1+x}
   \left(\frac9{64}-\frac{81x}{1024}\right)>0
 \qquad(0<x<1).
\]
Pairing with $e_0$ and using self-adjointness gives
\[
 0<\ip{(\mathcal A_{b_*}+25/32)f(\cos\beta)}{e_0}
   =(25/32-\kappa_0^2)\ip{f(\cos\beta)}{e_0}.
\]
Thus $\kappa_0^2<25/32$. The eigenfunction equation integrated
on $(0,L)$ gives
\[
 \frac34 e_0(L)=(\kappa_0^2+1/4)\int_0^L\sin\beta\,e_0\,\dd\beta.
\]
Since $\|e_0\|_{L^2((0,L),\omega\,\dd\beta)}=1$ and
$\int_0^L\sin\beta\,\dd\beta=1$, the integral is at most $1$.
By~\eqref{nfour:sum-residues},
\[
 2b_*A_*\leq\frac{8(\kappa_0^2+1/4)^2}{3(\kappa_0^2+9/4)}
 <\frac{363}{388}<\frac{15}{16}.
\]
The second inequality uses $\kappa_0^2<25/32$ and the increase
of $(x+1/4)^2/(x+9/4)$ for $x\geq0$.

Also $\kappa_0<9/10$ gives $3P/2>5$. The fourth and fifth
bounds in~\eqref{nfour:preliminary-table} are at most
$1+(6/5)e^{-7/2}<1+1/20$ and $1+2/99+1/50<1+1/20$,
respectively. Here $e^5>100$ and
$e^{7/2}>24$ follow from $e>8/3$. These estimates give
$b_*Q^D(t)<63/64$ for $t>0$ and $b_*Q^N(t)<63/64$ for $t\geq0$.

\medskip
\noindent\emph{Dimension four.}\enspace
Here $b_*=1$. By~\eqref{eq:flat-spectral-weighted-expression},
$\psi_j=\sin\beta\,e_j$ satisfies
\[
 -\psi_j''=\lambda_j\psi_j\quad\text{on }(0,L),\qquad
 \psi_j(0)=0,\qquad \psi_j'(L)=\psi_j(L).
\]
Multiplication by $\sin\beta$ preserves the norm from
$L^2((0,L),\omega\,\dd\beta)$ to $L^2(0,L)$. Conversely, the
sine and hyperbolic sine solutions with $\psi_j(0)=0$ give
smooth cylindrically symmetric functions after division by
$\sin\beta$, by their expansions at $\beta=0$. Thus
\begin{equation}\label{nfour:roots}
 \kappa_0\coth(L\kappa_0)=1,\qquad
 k_j\cot(Lk_j)=1,\qquad 2j<k_j<2j+1\quad(j\geq1).
\end{equation}
Indeed, the derivative of $x\cot(Lx)$ is
$\bigl(\sin(2Lx)/2-Lx\bigr)/\sin^2(Lx)<0$ where it is defined.
The function decreases from $+\infty$ to $0$ on $(2j,2j+1)$
for $j\geq1$. It is at most $1/L<1$ on $(0,1)$ and nonpositive
wherever the cotangent is nonpositive.
Its values at $5/2$ and $17/6$ are $5/2>1$ and
$17(2-\sqrt3)/6<1$, respectively. Hence $5/2<k_1<17/6$.

The eigenfunctions $\psi_j$ are multiples of
$\sinh(\kappa_0\beta)$ for $j=0$ and $\sin(k_j\beta)$ for $j\geq1$.
Integrating their squares and using~\eqref{nfour:roots} gives
\[
 \begin{aligned}
 \int_0^L\sinh^2(\kappa_0\beta)\,\dd\beta
   &=\frac{\sinh^2(L\kappa_0)}{2\kappa_0^2}-\frac L2,\\
 \int_0^L\sin^2(k_j\beta)\,\dd\beta
   &=\frac L2-\frac{\sin^2(Lk_j)}{2k_j^2}.
 \end{aligned}
\]
Since $\sin L=1$, division of the squared values at $L$ by
these integrals gives
\begin{equation}\label{nfour:residues}
 C_0=\frac{2\kappa_0^2}{1-L+L\kappa_0^2},\qquad
 C_j=\frac{2}{L+(L-1)/k_j^2}\quad(j\geq1).
\end{equation}

We claim that
\begin{equation}\label{nfour:constants}
 \begin{gathered}
 \frac{441}{500}<\kappa_0<\frac{89}{100},\qquad P>\frac72,\\
 C_0>2,\qquad C_1>1,\qquad 2A_*<1.
 \end{gathered}
\end{equation}
By~\eqref{nfour:roots}, $\tanh(L\kappa_0)=\kappa_0$.
The strict concavity of $\tanh$ on $(0,\infty)$ implies
that $\tanh(Ls)/s$ decreases strictly for $s>0$.

The exponential series gives $1957/720<e<11/4$ and, for
$s=23/100$,
\[
 e^s\leq 1+s+\frac{s^2}{2}+\frac{s^3}{6}
          +\frac{s^4}{24(1-s/5)}<\frac{1259}{1000}.
\]
Indeed, the terms of degree at least four have successive
ratios at most $s/5$. Consequently
\[
 e^{277/100}>
 \left(\frac{1957}{720}\right)^3\frac{1000}{1259}
 >\frac{941}{59}.
\]
Since $\pi>223/71$, we have $2L(441/500)>277/100$. The
ratio $(1+441/500)/(1-441/500)$ equals $941/59$. Hence
$\tanh(L(441/500))>441/500$ and $\kappa_0>441/500$.
For the upper bound, $L<11/7$ gives $2L(89/100)<3-1/5$.
Since
\[
 e^{3-1/5}<\frac{(11/4)^3}{1+1/5+1/50}<\frac{189}{11}
          =\frac{1+89/100}{1-89/100},
\]
we have $\tanh(L(89/100))<89/100$. Thus $\kappa_0<89/100$
and $P>7/2$.

The first denominator in~\eqref{nfour:residues} is positive
by the integral formula for $C_0$. It is less than
$\kappa_0^2$ since $L>1$ and $\kappa_0<1$, so $C_0>2$.
Also $k_1>5/2$ and $L<11/7$ give $C_1>350/291>1$.

To prove $2A_*<1$, write $x=\kappa_0^2$. Equations
\eqref{nfour:sum-residues} and~\eqref{nfour:residues} give
\[
 \frac{4}{2A_*}=(4+x)\frac{1-L+Lx}{x}.
\]
By~\eqref{nfour:roots},
$L=\operatorname{arctanh}(\sqrt{x})/\sqrt{x}$. Thus the
right-hand side is the value at $x=\kappa_0^2$ of the
strictly increasing function on $(0,1)$ given by
\[
 (4+x)\sum_{j=0}^\infty\frac{2x^j}{(2j+1)(2j+3)}.
\]
It suffices to prove that this function is greater than $4$
at $x=(441/500)^2$. For $s=441/500$, the inequality
$e^{277/100}>941/59$ gives
$2\operatorname{arctanh}(s)<277/100$. Since
\[
 \frac{277}{100}<\frac{2s(4-3s^2)}{(4+s^2)(1-s^2)},
\]
we have
\[
 (4+s^2)\frac{1-(1-s^2)\operatorname{arctanh}(s)/s}{s^2}>4.
\]
This proves~\eqref{nfour:constants}.

Raising the following inequalities to the fifth power reduces
them to comparisons with $(11/4)^2$ and $(11/4)^3$, respectively.
Thus $e<11/4$ gives
\begin{equation}\label{nfour:taylor-bounds}
 e^{2/5}<\frac32,\qquad e^{3/5}<2 .
\end{equation}
Since $P>7/2$ and $(1957/720)^7>1000$, we have $e^{2P}>1000$.
In the fourth line of~\eqref{nfour:preliminary-table},
$\kappa_0e^{-1-P}<1/80$ since $\kappa_0<1$ and
$e^{9/2}>80$. The last line is at most $1+2/999+1/500<81/80$.
Together with $2A_*<1$ these estimates give
\begin{equation}\label{nfour:coarse}
 Q^D(t)<\frac{1000}{999}\quad(t>0),\qquad
 Q^N(t)<\frac{81}{80}\quad(t\geq0).
\end{equation}
For $0\leq t\leq3/2$, with $t>0$ in the odd case, the first
and third lines of~\eqref{nfour:preliminary-table} give the
stronger bound $2A_*<1$ for both quantities since $3/2<P/2$.

Let $t\geq3/2$. We use the terms with $j=0$ and $j=1$ in
\eqref{nfour:odd-kernel} and~\eqref{nfour:even-kernel}
on the interval $1/20\leq t-\tau\leq1/5$. On this interval
the $j=0$ term is negative. The inequality
$S(t-\tau)\geq(t-\tau)(1-\kappa_0^2(t-\tau)^2/6)$ and
$C_0>2$ give $C_0S(t-\tau)>149/1500$.
Also $\tau\geq13/10$ and $k_1>5/2$, so $k_1\tau>3$ and
$e^{-k_1\tau}<1/20$. Since $k_1<17/6$, we have
$k_1(t-\tau)<3/5$ and by~\eqref{nfour:taylor-bounds},
\[
 \frac{e^{-k_1(t-\tau)}}{2k_1}>\frac3{34}.
\]
Equation~\eqref{nfour:high-representation} now gives
$G_{k_1}^D(t,\tau)>(399/400)(3/34)$ and
$G_{k_1}^N(t,\tau)>3/34$. Since $|S(t)|<8/7$,
$k_1^{-1}<2/5$ and $C_1>1$, the $j=1$ terms of
\eqref{nfour:odd-kernel} and~\eqref{nfour:even-kernel}
are bounded below by
\[
 \frac{399}{400}\frac3{34}-\frac2{35}>\frac3{100},
 \qquad
 \frac3{34}-\frac1{50}>\frac1{15},
\]
respectively. On the same interval,
\[
 \begin{aligned}
 \frac{\sinh(2\tau)}{\sinh(2t)}
 &\geq e^{-2(t-\tau)}(1-e^{-4\tau})>\frac{13}{20},\\
 \frac{\cosh(2\tau)}{\cosh(2t)}
 &\geq e^{-2(t-\tau)}>\frac23 ,
 \end{aligned}
\]
by~\eqref{nfour:taylor-bounds} and $e^{4\tau}>40$.

For real numbers $u,w$ of opposite sign,
$|u|+|w|-|u+w|=2\min(|u|,|w|)$. The terms with $j=0$ and
$j=1$ have opposite signs throughout the interval
$1/20\leq t-\tau\leq1/5$. Its length is $3/20$, so
combining these terms before taking the absolute value
decreases the bound for $Q^D(t)$ by at least
$2(3/100)(13/20)(3/20)=117/20000$. It decreases the bound for
$Q^N(t)$ by at least $2(1/15)(2/3)(3/20)=1/75$.
In particular, for $t\geq3/2$,
\[
 \begin{aligned}
 Q^D(t)&<\frac{1000}{999}-\frac{117}{20000}<1,\\
 Q^N(t)&<\frac{81}{80}-\frac1{75}=\frac{1199}{1200}<1.
 \end{aligned}
\]

\medskip
For the weight $h$ chosen above, set
\[
 \|Y(\cdot,L)\|_h=\sup_{\substack{t\geq0\\h(t)>0}}
                       \frac{|Y(t,L)|}{h(t)} .
\]
This is finite by smoothness, parity and polynomial growth.
Since $0\leq b\leq2b_*$, we have $|F(t)|\leq b_*|Y(t,L)|$.
By~\eqref{nfour:trace-representation}, the ratio $|Y(t,L)|/h(t)$
is bounded by $b_*Q^D(t)\|Y(\cdot,L)\|_h$ in the odd case and
by $b_*Q^N(t)\|Y(\cdot,L)\|_h$ in the even case. Taking the
supremum in $t$ and using the bounds on $b_*Q^D$ and $b_*Q^N$
gives
\begin{equation}\label{nfour:contraction}
 \begin{aligned}
 \|Y(\cdot,L)\|_h&\leq\frac{63}{64}\|Y(\cdot,L)\|_h &&(n=3),\\
 \|Y(\cdot,L)\|_h&\leq q\|Y(\cdot,L)\|_h &&(n=4),
 \end{aligned}
\end{equation}
where
\[
 q=\max\left\{2A_*,\frac{1000}{999}-\frac{117}{20000},
                    \frac{1199}{1200}\right\}<1.
\]
Hence $Y(t,L)=0$ and $F(t)=0$ on $[0,\infty)$.
By~\eqref{nfour:high-representation}, $Y_j=0$ for $j\geq1$.
The expansion at $L$ then gives $e_0(L)Y_0(t)=Y(t,L)=0$.
Since $e_0(L)>0$, completeness and parity give $Y\equiv0$
on $\R\times\Sph^{n-1}_+$.
\end{proof}


\begin{thebibliography}{99}

\bibitem{ADN}
S. Agmon, A. Douglis and L. Nirenberg,
\href{https://doi.org/10.1002/cpa.3160120405}
  {\emph{Estimates near the boundary for solutions of elliptic
  partial differential equations satisfying general boundary
  conditions. I}},
Comm. Pure Appl. Math. \textbf{12} (1959), 623--727.

\bibitem{Andrade}
J. H. Andrade, A. DelaTorre, J. M. do \'O, J. Ratzkin and J. Wei,
\href{https://arxiv.org/abs/2511.05225v2}
  {\emph{Classification of fractional, singular Yamabe metrics
  on a twice-punctured sphere I}},
\href{https://arxiv.org/abs/2511.05225v2}{arXiv:2511.05225v2} (2025).

\bibitem{AoDGW}
W. Ao, A. DelaTorre, M. del Mar Gonz\'alez and J. Wei,
\href{https://doi.org/10.1515/crelle-2018-0032}
  {\emph{A gluing approach for the fractional Yamabe problem
  with isolated singularities}},
J. Reine Angew. Math. \textbf{763} (2020), 25--78.

\bibitem{PartII}
F.~G.~Bogo,
\emph{Jacobi Fields and Delaunay Asymptotics for the Singular
Escobar-Yamabe Equation},
companion manuscript.

\bibitem{Brezis}
H. Brezis,
\href{https://doi.org/10.1007/978-0-387-70914-7}
  {\emph{Functional Analysis, Sobolev Spaces and Partial
  Differential Equations}},
Universitext, Springer, New York, 2011.

\bibitem{CGS}
L. A. Caffarelli, B. Gidas and J. Spruck,
\href{https://math.jhu.edu/~js/Math646/asymmetry.pdf}
  {\emph{Asymptotic symmetry and local behavior of semilinear
  elliptic equations with critical Sobolev growth}},
Comm. Pure Appl. Math. \textbf{42} (1989), 271--297.

\bibitem{CJXS}
L. A. Caffarelli, T. Jin, Y. Sire and J. Xiong,
\href{https://arxiv.org/abs/1309.2021}
  {\emph{Local analysis of solutions of fractional semi-linear
  elliptic equations with isolated singularities}},
Arch. Ration. Mech. Anal. \textbf{213} (2014), 245--268.

\bibitem{CS}
L. A. Caffarelli and L. Silvestre,
\href{https://arxiv.org/abs/math/0608640}
  {\emph{An extension problem related to the fractional Laplacian}},
Comm. Partial Differential Equations \textbf{32} (2007), 1245--1260.

\bibitem{CaseChang}
J. S. Case and S.-Y. A. Chang,
\href{https://doi.org/10.1002/cpa.21564}
  {\emph{On fractional GJMS operators}},
Comm. Pure Appl. Math. \textbf{69} (2016), 1017--1061.

\bibitem{ChangGonzalez}
S.-Y. A. Chang and M. del Mar Gonz\'alez,
\href{https://doi.org/10.1016/j.aim.2010.07.016}
  {\emph{Fractional Laplacian in conformal geometry}},
Adv. Math. \textbf{226} (2011), 1410--1432.

\bibitem{CDR}
S. Cruz-Bl\'azquez, A. DelaTorre and D. Ruiz,
\href{https://doi.org/10.3934/dcds.2023077}
  {\emph{Qualitative properties of singular solutions to the
  fractional Yamabe problem}},
Discrete Contin. Dyn. Syst. \textbf{43} (2023), 4027--4041.

\bibitem{DDPGW}
A. DelaTorre, M. del Pino, M. del Mar Gonz\'alez and J. Wei,
\href{https://arxiv.org/abs/1510.08504}
  {\emph{Delaunay-type singular solutions for the fractional
  Yamabe problem}},
Math. Ann. \textbf{369} (2017), 597--626.

\bibitem{DG}
A. DelaTorre and M. del Mar Gonz\'alez,
\href{https://arxiv.org/abs/1504.03493}
  {\emph{Isolated singularities for a semilinear equation for
  the fractional Laplacian arising in conformal geometry}},
Rev. Mat. Iberoam. \textbf{34} (2018), 1645--1678.

\bibitem{Delaunay}
Ch. Delaunay,
\href{https://www.numdam.org/item/JMPA_1841_1_6__309_0/}
  {\emph{Sur la surface de r\'evolution dont la courbure
  moyenne est constante}},
J. Math. Pures Appl. (1) \textbf{6} (1841), 309--315.

\bibitem{Escobar}
J. F. Escobar,
\href{https://annals.math.princeton.edu/1992/136-1/p01}
  {\emph{Conformal deformation of a Riemannian metric to a scalar
  flat metric with constant mean curvature on the boundary}},
Ann. of Math. (2) \textbf{136} (1992), 1--50.

\bibitem{JLX}
T. Jin, Y. Y. Li and J. Xiong,
\href{https://doi.org/10.4171/JEMS/456}
  {\emph{On a fractional Nirenberg problem, part I: blow up
  analysis and compactness of solutions}},
J. Eur. Math. Soc. \textbf{16} (2014), 1111--1171.

\bibitem{JinLiaoXiong}
T. Jin, Y. Liao and J. Xiong,
\href{https://arxiv.org/abs/2611.07636}
  {\emph{Classification and asymptotics for the boundary Yamabe
  equation with an isolated singularity}},
arXiv:2611.07636 (2026).

\bibitem{KMPS}
N. Korevaar, R. Mazzeo, F. Pacard and R. Schoen,
\href{https://arxiv.org/abs/math/9807038}
  {\emph{Refined asymptotics for constant scalar curvature metrics
  with isolated singularities}},
Invent. Math. \textbf{135} (1999), 233--272.

\bibitem{LiZhu}
Y. Y. Li and M. Zhu,
\href{https://doi.org/10.1215/S0012-7094-95-08016-8}
  {\emph{Uniqueness theorems through the method of moving spheres}},
Duke Math. J. \textbf{80} (1995), 383--417.

\bibitem{Marques}
F. C. Marques,
\href{https://doi.org/10.1007/s00526-007-0144-3}
  {\emph{Isolated singularities of solutions to the Yamabe equation}},
Calc. Var. Partial Differential Equations \textbf{32} (2008), 349--371.

\bibitem{MP}
R. Mazzeo and F. Pacard,
\href{https://perso.pages.math.cnrs.fr/users/frank.pacard/Publications/26.pdf}
  {\emph{Constant scalar curvature metrics with isolated singularities}},
Duke Math. J. \textbf{99} (1999), 353--418.

\bibitem{MPU}
R. Mazzeo, D. Pollack and K. Uhlenbeck,
\href{https://arxiv.org/abs/dg-ga/9406004}
  {\emph{Moduli spaces of singular Yamabe metrics}},
J. Amer. Math. Soc. \textbf{9} (1996), 303--344.

\bibitem{McLean}
W. McLean,
\href{https://assets.cambridge.org/97805216/63328/frontmatter/9780521663328_frontmatter.pdf}
  {\emph{Strongly Elliptic Systems and Boundary Integral Equations}},
Cambridge University Press, Cambridge, 2000.

\bibitem{Serrin}
J. Serrin,
\href{https://doi.org/10.1007/BF02391014}
  {\emph{Local behavior of solutions of quasi-linear equations}},
Acta Math. \textbf{111} (1964), 247--302.

\bibitem{SilvaSantos}
A. Silva Santos,
\href{https://doi.org/10.1007/s00023-010-0024-9}
  {\emph{A construction of constant scalar curvature manifolds
  with Delaunay-type ends}},
Ann. Henri Poincar\'e \textbf{10} (2010), 1487--1535.

\end{thebibliography}
\end{document}